\documentclass{amsart}

\usepackage{amssymb, stmaryrd, ifthen, comment, mathtools, array, xspace, tikz, booktabs}
	\usetikzlibrary{matrix,fit,positioning,arrows.meta}
	\mathtoolsset{centercolon}

\usepackage[shortlabels]{enumitem}
	\setenumerate[0]{label={\rm (\roman*)}, leftmargin=*}

\usepackage[T1]{fontenc}
\usepackage[english]{babel}

\usepackage{hyperref}\hypersetup{bookmarksopen=true}

\usepackage[capitalize]{cleveref}
	\crefname{enumi}{}{}
	\Crefname{enumi}{Item}{Items}
	\crefname{equation}{}{}
	\Crefname{equation}{Equation}{Equations}

\numberwithin{equation}{section}
\renewcommand{\thesubsection}{\thesection.\Alph{subsection}}

\newtheorem{theorem}[equation]{Theorem}

\newtheorem{lemma}[equation]{Lemma}
\newtheorem{proposition}[equation]{Proposition}
\newtheorem{corollary}[equation]{Corollary}

\theoremstyle{definition}
\newtheorem{remark}[equation]{Remark}
\newtheorem*{remark*}{Remark}
\newtheorem{definition}[equation]{Definition}
\newtheorem{notation}[equation]{Notation}
\newtheorem{construction}[equation]{Construction}
\newtheorem{example}[equation]{Example}

\newtheorem*{construction*}{Construction}

\newcommand{\fG}{\mathfrak{g}} 
\newcommand{\fH}{\mathfrak{h}} 
\newcommand{\fu}{\mathfrak{u}} 
\newcommand{\ssl}{\mathfrak{sl}} 
\newcommand{\instrl}{\mathfrak{instrl}} 
\newcommand{\TKK}{K(\A)}

\newcommand{\Z}{\mathbb{Z}} 

\newcommand{\A}{\mathcal{A}} 
\newcommand{\B}{\mathcal{B}} 
\newcommand{\HH}{\mathcal{H}} 
\renewcommand{\SS}{\mathcal{S}} 

\newcommand{\id}{\mathrm{id}}
\newcommand{\op}{\mathrm{op}}
\newcommand{\dash}{\nobreakdash-\hspace{0pt}}
\newcommand{\ov}{\overline}

\newcommand{\subhalf}{{\leavevmode \raise.5ex\hbox{\the\scriptfont0 1}\kern-.13em /\kern-.07em\lower.25ex\hbox{\the\scriptfont0 2}}}
\newcommand{\idempotentmatrix}{\begin{psmallmatrix} 1 & 0 \\ 0 & 0 \end{psmallmatrix}}
\newcommand{\into}{\hookrightarrow}

\DeclareMathOperator{\Aut}{Aut}

\DeclareMathOperator{\Char}{char}

\DeclareMathOperator{\End}{End}

\DeclareMathOperator{\Fix}{Fix}

\DeclareMathOperator{\Trd}{Trd}

\DeclareMathOperator{\ad}{ad}

\DeclareMathOperator{\Zorn}{Zorn}
\DeclareMathOperator{\pr}{pr}

\allowdisplaybreaks
\begin{document}

\title{Lie algebras graded by the root system $AG_2$}
\author{Tom De Medts}
\address{Tom De Medts \\ Ghent University \\ Department of Mathematics, Computer Science and Statistics \\ Krijgslaan 299, S9 \\ B-9000 Gent \\ Belgium}
\email{tom.demedts@ugent.be}

\author{Michiel Smet}
\address{Michiel Smet \\ Ghent University \\ Department of Electronics and Information Systems \\ Krijgslaan 297, S8 \\ B-9000 Gent \\ Belgium}
\email{michiel.smet@ugent.be}

\dedicatory{Dedicated to Richard Weiss, friend and inexhaustible source of inspiration, \\ on the occasion of his 80th birthday}
\date{November 9, 2026}

\begin{abstract}
	We study Lie algebras graded by the generalized root system $AG_2$, which is obtained from $G_2$ by adjoining the doubles of its short roots, over fields of characteristic different from $2$ and $3$.
	We show that every such Lie algebra with a grading subalgebra of type $A_2$ is, modulo its center, the Tits--Kantor--Koecher Lie algebra of a structurable algebra equipped with an idempotent~$t$ satisfying $t+\overline{t}=1$.
	
	Conversely, every structurable algebra with such an idempotent gives rise to an $(AG_2,A_2)$-graded Lie algebra. We describe the corresponding Peirce decomposition and characterize the additional algebraic structure required for an $(AG_2,G_2)$-grading.
	Using the classification of central simple structurable algebras, we determine the central simple examples when the characteristic is also different from $5$.
	Finally, we compare our constructions with the theory of Tits hexagons. For the central simple examples in which all root spaces are non-zero, we determine which ones are Tits-hexagonalizable; in certain cases, this imposes an additional balancing condition on the idempotent.
\end{abstract}

\makeatletter
\@namedef{subjclassname@2020}{\textup{2020} Mathematics Subject Classification}
\makeatother

\maketitle

\section{Introduction}

The original motivation for this paper comes from the recent work of Bernhard M\"uhlherr and Richard Weiss on Tits hexagons~\cite{Muhlherr26}. In that paper, the root system $AG_2$ plays a key role in describing the commutator relations of the root groups.
This root system $AG_2$ is obtained from the usual root system $G_2$ by adjoining the doubles of its six short roots. Although this is not a root system in the usual sense, it is a generalized root system in the sense of \cite{Dimitrov24,Cuntz24}.
\[ \begin{tikzpicture}[scale=1, line cap=round, line join=round]

\def\R{1}                 
\def\L{1.732050807568877} 
\def\D{2}                 

\tikzset{
  short/.style={-Latex, very thick},
  long/.style ={-Latex, very thick, black!20!red},
  dbl/.style  ={dash pattern=on 2pt off 5pt, -Latex, very thick},
}

\foreach \a in {0,60,120,180,240,300}{
  \draw[short] (0,0) -- (\a:\R);
  \draw[short] (0,0) -- (\a+180:\R);

  \draw[long]  (0,0) -- (\a+30:\L);
  \draw[long]  (0,0) -- (\a+210:\L);

  \draw[dbl]   (0,0) -- (\a:\D);
  \draw[dbl]   (0,0) -- (\a+180:\D);
}

\end{tikzpicture} \]
In view of our earlier work on groups and Lie algebras graded by root systems, where the root groups and the root-graded Lie algebras arise side by side from the same algebraic data (see, for instance, \cite{Boelaert2019,DeMedtsMeulewaeter2025,DeMedts2026}), this raises the natural question of what Lie algebras graded by $AG_2$ look like, and how they are related to the Tits hexagons of \cite{Muhlherr26}.
Our aim is to answer these questions, relying on the theory of structurable algebras.

The study of Lie algebras graded by finite root systems was initiated by Berman and Moody~\cite{Berman1992}, who treated the simply-laced root systems of rank at least~$2$, and continued by Benkart and Zelmanov~\cite{Benkart1996} for the non-simply-laced cases.
(The original motivation came from intersection matrix algebras and their relation to toroidal Lie algebras.)
In this theory, a Lie algebra contains a split simple grading subalgebra, and its root-space decomposition is required to be the weight-space decomposition for a split Cartan subalgebra of this grading algebra.
The Lie algebra can then be described, up to central extensions, using a coordinate algebra whose structure depends on the root system.
For instance, in type $A_2$, the coordinate algebra is alternative.
The theory has also been extended to the non-reduced root systems $BC_r$; see, in particular, Allison, Benkart and Gao~\cite{ABG02} for $r\geq 2$.

In our setting, there are two natural choices for the grading algebra.
The long roots of $AG_2$ form a root subsystem of type $A_2$, whereas the short and long roots together form a root subsystem of type $G_2$.
We therefore introduce $(AG_2,A_2)$-graded and $(AG_2,G_2)$-graded Lie algebras.
The first notion is more general, and most of the paper is devoted to this case.
Notice that the roots of $AG_2$ belong to the weight lattice of the long-root $A_2$, but not all of them belong to its root lattice.
Our definition therefore extends the usual notion of an $A_2$-graded Lie algebra (see \cref{rem:omega}).
A related situation was studied by Benkart and Elduque~\cite{BE12}, who considered Lie algebras decomposing, as modules for a subalgebra $\mathfrak{sl}_3$, into copies of the adjoint module, the natural $3$-dimensional module, its dual, and the trivial module.
They also obtained connections with structurable algebras.
Their setting does not include the doubled short-root spaces that we allow here.

As we learned after our paper was essentially finished,
the doubled short-root spaces do appear in the work of Yaseen~\cite{Yaseen22} (see also \cite{Yaseen18}), who studied, in characteristic~$0$, Lie algebras graded by the weight system $(\Theta_3, \mathfrak{sl}_3)$, where
\[ \Theta_3 = \{ 0, \ \pm\varepsilon_i \pm \varepsilon_j, \ \pm\varepsilon_i, \ \pm 2\varepsilon_i \mid 1 \leq i \neq j \leq 3 \} . \]
Since $\varepsilon_1 + \varepsilon_2 + \varepsilon_3 = 0$ for the weights of $\mathfrak{sl}_3$, we have $\varepsilon_i + \varepsilon_j = -\varepsilon_k$ whenever $\{ i,j,k \} = \{ 1,2,3 \}$, and $\Theta_3 \setminus \{ 0 \}$ is precisely the generalized root system $AG_2$, with short roots $\pm\varepsilon_i$, long roots $\varepsilon_i - \varepsilon_j$ and very long roots $\pm 2\varepsilon_i$.
(In other words, $AG_2$ is the image of the root system $BC_3$ under the restriction of weights from $\mathfrak{gl}_3$ to~$\mathfrak{sl}_3$.)
In characteristic~$0$, the $(\Theta_3, \mathfrak{sl}_3)$-graded Lie algebras of \cite[Definition~2.1]{Yaseen22} are exactly the $(AG_2, A_2)$-graded Lie algebras of \cref{def:AG2-Delta-graded} below.
Yaseen decomposes such a Lie algebra into isotypic components for the grading subalgebra~$\mathfrak{sl}_3$, describes the Lie bracket in terms of the corresponding multiplicity spaces, and assembles these into a unital coordinate algebra with involution; see also \cite{BaranovYaseen21} for the analogous question for $\mathfrak{sl}_n$ with larger~$n$.
This is done under an additional assumption which, as we explain in \cref{rem:Yaseen}, excludes precisely those central simple examples that correspond to the large exceptional Tits hexagons from~\cite{Muhlherr26}.

Our approach is different: we identify the coordinate algebra as a structurable algebra equipped with a conjugate-supplementary idempotent, which also yields a converse construction and allows us to determine the central simple examples.

The main ingredient in our description is a structurable algebra equipped with a \emph{conjugate-supplementary idempotent}, i.e., an idempotent $t$ satisfying $t+\overline{t}=1$.
Structurable algebras were introduced by Allison~\cite{Allison1978}; they generalize both Jordan algebras and associative algebras with involution.
They admit a Tits--Kantor--Koecher construction producing a $5$-graded Lie algebra $K(\A)$ from a structurable algebra $\A$.
We show that every $(AG_2,A_2)$-graded Lie algebra $L$ gives rise to a structurable algebra $\A$ with a conjugate-supplementary idempotent, such that
\[ L/Z(L)\cong K(\A). \]
Conversely, every structurable algebra with a conjugate-supplementary idempotent gives rise to an $(AG_2,A_2)$-grading on its TKK Lie algebra.
These are \cref{th:AG2-A2,th:AG2-A2-conv}.

The idempotent $t$ induces a \emph{Peirce decomposition} of $\A$, which we write in matrix form as
\[ \A = \begin{pmatrix} B & J \\ J' & B^\op \end{pmatrix} . \]
We show that $B = t\A t$ is an alternative algebra.
The multiplication and involution of $\A$ give rise to bimodule structures on $J$ and $J'$ and to bilinear maps between these spaces.
We collect these data in the notion of an $AG_2$-algebra.
(We provide a self-contained list of axioms for $AG_2$-algebras in \cref{se:appendix}.)
This description makes the different root spaces of $K(\A)$ explicit.
It also includes the usual $A_2$\dash graded case from \cite{Berman1992}: when $J=J'=0$, the structurable algebra is $B\oplus B^{\op}$ with the exchange involution.

We next determine the additional structure corresponding to an $(AG_2,G_2)$-grading.
It turns out that this is equivalent to the existence of a ``$G_2$-frame'' in~$(\A,t)$, consisting of two elements in the off-diagonal Peirce spaces satisfying a short list of identities.
These elements, together with $t$ and $\overline{t}$, span a $4$-dimensional structurable subalgebra whose TKK Lie algebra is split of type $G_2$.
This gives the characterization in \cref{th:G2-frame}.

Next, we examine the examples arising from the known classification of central simple structurable algebras.
We determine which algebras admit conjugate-supplementary idempotents and describe the corresponding $AG_2$-algebras.
The results are summarized in \cref{th:classification}, and the existence of $G_2$-frames is determined in \cref{co:G2-frame}.
When $\Char(k) \neq 5$, the Allison--Smirnov classification ensures that this list is complete.

Finally, we return to the Tits hexagons that motivated our work.
An $AG_2$\dash grading gives rise to root groups with commutator relations reflecting the root system, but this does not by itself ensure the existence of a corresponding Tits hexagon.
The interaction between opposite root groups, in the form of Weyl elements, imposes additional conditions.
We describe the root groups explicitly, in terms of the $\epsilon$-pseudo-quadratic modules occurring in~\cite{Muhlherr26} for our examples of hermitian type, and in terms of the large exceptional Tits hexagons from \cite[\S 13]{Muhlherr26} for the exceptional examples.
For the central simple examples in which all root spaces are non-zero, we then determine which ones are Tits-hexagonalizable (\cref{co:hexagonalizable}).
In certain cases, this forces a ``balancing condition'' on the idempotent, and there is one further family of exceptions, for which $K(\A)$ is a split orthogonal Lie algebra of type~$D_{3n}$ with $n$ odd.
In particular, we obtain examples of $AG_2$-graded Lie algebras that do not give rise to Tits hexagons in this way.
Our description also includes the large exceptional Tits hexagons and the classical Tits hexagons of ``type~$+1$'' \cite[Observation~12.9]{Muhlherr26}, both of which are excluded by the hypotheses of the main result \cite[Theorem~11.1]{Muhlherr26} of M\"uhlherr and Weiss.
Finally, we interpret these results in terms of an auxiliary $BC_3$-grading, in the spirit of \cite[\S 11]{Muhlherr26}.

The paper is organized as follows.
\Cref{se:Lie AG2} introduces the gradings under consideration.
\Cref{se:idempotents,se:TKK} recall the required facts about structurable algebras and their TKK Lie algebras.
In \cref{se:AG2 to stralg}, we pass from an $(AG_2,A_2)$-graded Lie algebra to a structurable algebra with a conjugate-supplementary idempotent.
\Cref{se:idempotents2,se:AG2 TKK} describe the Peirce decomposition and establish the converse construction.
\Cref{se:AG2 G2} deals with $(AG_2,G_2)$-gradings, \cref{se:examples,se:classification} contain the examples and their classification, and \cref{se:tits hex} discusses the connection with Tits hexagons.

\subsection*{Acknowledgments}
We thank Benny Willen for drawing our attention to the work of Yaseen~\cite{Yaseen22}.
We are grateful to Bernhard M\"uhlherr and Richard Weiss for illuminating discussions about the root system $AG_2$ and Tits hexagons.
Their preprint \cite{Muhlherr26} directly inspired this paper.

In a recent discussion, Petra Schwer pointed out that AI declarations are often too vague; we have tried to follow her advice to make ours more detailed (see the AI Declaration at the end of the paper).

\bigskip

Throughout this paper, $k$ will be a commutative field with $\Char(k) \not\in \{ 2, 3 \}$.

\section{Lie algebras graded by the root system \texorpdfstring{$AG_2$}{AG2}}\label{se:Lie AG2}

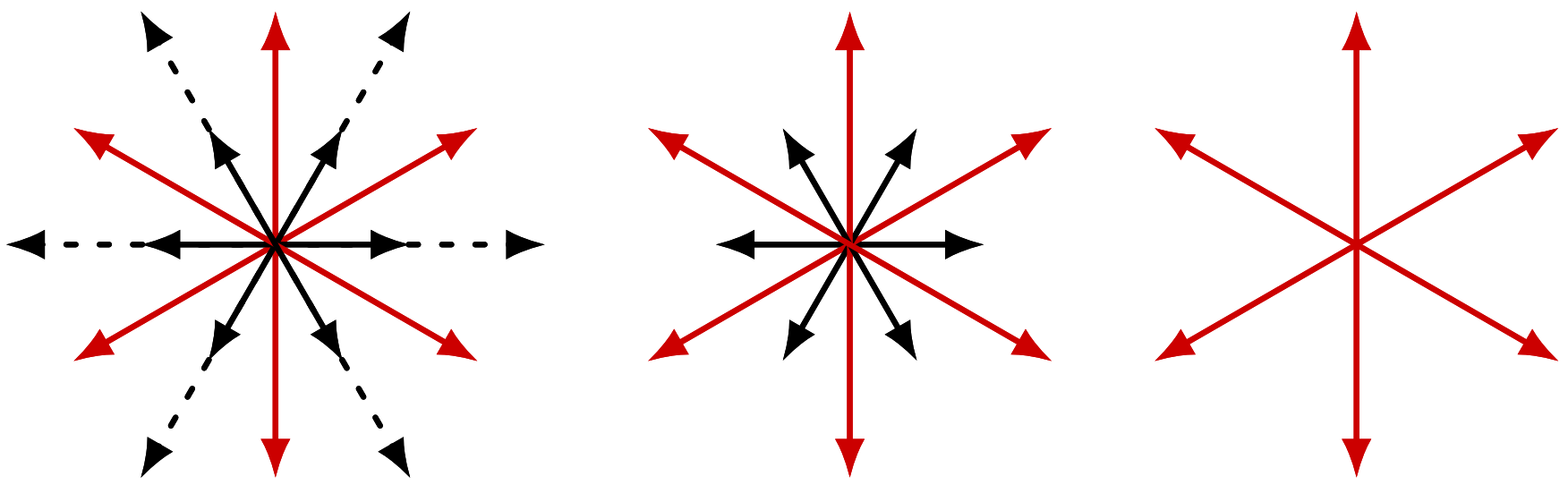
\begin{figure}[ht!]
\[
\begin{tikzpicture}[scale=1, line cap=round, line join=round]

\def\R{1}                 
\def\L{1.732050807568877} 
\def\D{2}                 

\tikzset{
  short/.style={-Latex, very thick},
  long/.style ={-Latex, very thick, black!20!red},
  dbl/.style  ={dash pattern=on 2pt off 5pt, -Latex, very thick},
}

\foreach \a in {0,60,120,180,240,300}{
  \draw[short] (0,0) -- (\a:\R);
  \draw[short] (0,0) -- (\a+180:\R);

  \draw[long]  (0,0) -- (\a+30:\L);
  \draw[long]  (0,0) -- (\a+210:\L);

  \draw[dbl]   (0,0) -- (\a:\D);
  \draw[dbl]   (0,0) -- (\a+180:\D);
}

\end{tikzpicture}
\qquad
\begin{tikzpicture}[scale=1, line cap=round, line join=round]

\def\R{1}                 
\def\L{1.732050807568877} 
\def\D{2}                 

\tikzset{
  short/.style={-Latex, very thick},
  long/.style ={-Latex, very thick, black!20!red},
  dbl/.style  ={dash pattern=on 2pt off 5pt, -Latex, very thick},
}

\foreach \a in {0,60,120,180,240,300}{
  \draw[short] (0,0) -- (\a:\R);
  \draw[short] (0,0) -- (\a+180:\R);

  \draw[long]  (0,0) -- (\a+30:\L);
  \draw[long]  (0,0) -- (\a+210:\L);
}

\end{tikzpicture}
\qquad
\begin{tikzpicture}[scale=1, line cap=round, line join=round]

\def\R{1}                 
\def\L{1.732050807568877} 
\def\D{2}                 

\tikzset{
  short/.style={-Latex, very thick},
  long/.style ={-Latex, very thick, black!20!red},
  dbl/.style  ={dash pattern=on 2pt off 5pt, -Latex, very thick},
}

\foreach \a in {0,60,120,180,240,300}{
  \draw[long]  (0,0) -- (\a+30:\L);
  \draw[long]  (0,0) -- (\a+210:\L);
}

\end{tikzpicture}
\]
\caption{\small The $AG_2$ generalized root system and its distinguished $G_2$ and $A_2$ root subsystems}\label{fig:AG2 GRS}
\end{figure}

To define the ``root system'' $AG_2$, we start with the usual root system $G_2$ and for each short root, we add the double of this root. As a result, we get $18$ roots of three different lengths ($1$, $\sqrt{3}$ and $2$).
This is, of course, not a root system in the usual sense, but it is a ``generalized root system'' in the sense of \cite{Dimitrov24,Cuntz24}.
(It appears as item (2)(vi) in the list of \cite[\S 5.2]{Dimitrov24}.)
We will refer to the roots of each of the lengths as \emph{short}, \emph{long} and \emph{very long}. In particular, the $G_2$ root subsystem consisting of the solid lines in \cref{fig:AG2 GRS} has short and long roots (as usual), and the $A_2$ root subsystem consisting of the red lines has long roots only.
Of course, there are also two other $A_2$ root subsystems consisting of short roots and very long roots, respectively, but these will not be explicitly needed now.

The weakest notion of a Lie algebra reflecting the structure of this root system is the following.
We use the notation%
\footnote{We use the letter $o$ for the origin of the Euclidean space in which the root system lives. It will be convenient to reserve the number $0$ for a $\Z$-gradation.}
$\Phi^o := \Phi \cup \{ o \}$ for any (generalized) root system~$\Phi$.
\begin{definition}\label{def:supported}
	Let $L$ be a Lie algebra over $k$.
	We say that $L$ is \emph{supported by the root system $AG_2$} if
	\begin{enumerate}
		\item we can decompose $L$ as
			\[ L = \bigoplus_{\alpha \in AG_2^o} L_\alpha \]
			in such a way that $[L_\alpha, L_\beta] \subseteq L_{\alpha+\beta}$ for all $\alpha,\beta \in AG_2^o$, where we set
			$L_\gamma = 0$ if $\gamma$ is not contained in $AG_2^o$.
			We call these $L_\alpha$ \emph{root spaces}.
		\item\label{def:supported:2} $L$ is generated by the root spaces $L_\alpha$ with $\alpha \in AG_2$ (excluding $\alpha = o$).
			In particular, by the Jacobi identity, $L_o = \sum_{\alpha \in AG_2} [L_\alpha, L_{-\alpha}]$.
	\end{enumerate}
	We say that $L$ is \emph{strongly supported by the root system $AG_2$} if, in addition, all $L_\alpha$ (for $\alpha \in AG_2^o$) are non-zero.
\end{definition}

This definition in itself is still too weak: it is easy to construct rather degenerate examples of Lie algebras strongly supported by $AG_2$.
For example, start with one-dimensional vector spaces $L_\alpha := k e_\alpha$ for each $\alpha \in AG_2^o$, and declare $[e_\alpha, e_{-\alpha}] = \pm e_o$ for all $\alpha \in AG_2$ and all other Lie brackets $[e_\alpha, e_\beta] = 0$; this results in a two-step nilpotent Lie algebra.

In the spirit of the existing theory of Lie algebras graded by root systems \cite{Berman1992,Benkart1996}, it is natural to impose additional symmetry, giving rise to the following more restrictive definition.

\begin{definition}\label{def:AG2-Delta-graded}
	Let $L$ be a Lie algebra over $k$ and assume that $L = \bigoplus_{\alpha \in AG_2^o} L_\alpha$ is supported by $AG_2$.
	Let $\Delta$ be either the ``red root subsystem'' $A_2$ of $AG_2$ or the ``solid root subsystem'' $G_2$ of $AG_2$ (see~\cref{fig:AG2 GRS}).
	
	We say that $L$ is \emph{$(AG_2, \Delta)$-graded} if
	\begin{enumerate}
		\item $L$ contains a split simple Lie algebra $\fG$ over $k$ of type $\Delta$ with a split Cartan subalgebra $\fH$ such that for each root $\gamma \in \Delta^o$, the root space $\fG_\gamma$ is contained in~$L_\gamma$. In particular, $\fH \subseteq L_o$;
		\item for each $h \in \fH$ and each $x \in L_\alpha$ (with $\alpha \in AG_2^o$), we have $[h, x] = \alpha(h) x$.
	\end{enumerate}
	We call $\fG$ the \emph{grading algebra} of $L$.
\end{definition}
\begin{remark}\label{rem:omega}
	In the original setting of \cite{Berman1992,Benkart1996}, the weights of $L$ with respect to $\fH$ are required to lie in $\Delta^o$, i.e., $L$ is graded by the root system of its grading algebra.
	In our setting, the set $AG_2^o$ of possible weights is strictly larger than $\Delta^o$ for both choices of $\Delta$, and when $\Delta = A_2$, it is not even contained in the root lattice of~$\Delta$.
	Situations of this kind have been considered before: by Allison, Benkart and Gao~\cite{ABG02} for the root systems $BC_r$ (with a grading algebra of type $B_r$, $C_r$ or~$D_r$), by Benkart and Elduque~\cite{BE12}, and, more generally, for finite weight systems in \cite{BaranovYaseen21,Yaseen22}; as explained in the introduction, the $(\Theta_3, \ssl_3)$-graded Lie algebras of \cite{Yaseen22} are exactly our $(AG_2, A_2)$-graded Lie algebras (in characteristic~$0$).
	
	Our two possibilities for $\Delta$ ($A_2$ or $G_2$) do have the same \emph{weight lattice}, which, by~(i), we view as a lattice in $\fH^*$. Since all $\alpha \in AG_2^o$ are contained in this weight lattice, the expression $\alpha(h)$ in condition~(ii) makes sense.
	In \cref{fig:weights}, we show a choice of fundamental roots $\alpha_1, \alpha_2$ for $A_2$ and the corresponding fundamental weights $\omega_1, \omega_2$, such that $\omega_i(h_j) = \delta_{ij}$ for all $i,j \in \{ 1,2 \}$.
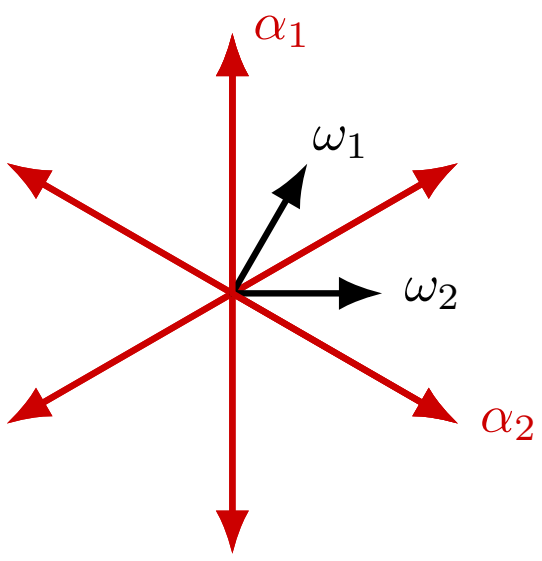
\begin{figure}[ht!]
	\[
	\begin{tikzpicture}[scale=1, line cap=round, line join=round]
		
	\def\R{1}                 
	\def\L{1.732050807568877} 
	\def\D{2}                 
	
	\tikzset{
	  short/.style={-Latex, very thick},
	  long/.style ={-Latex, very thick, black!20!red},
	  dbl/.style  ={dash pattern=on 2pt off 5pt, -Latex, very thick},
	}
	
	\foreach \a in {0,60,120,180,240,300}{
	  \draw[long]  (0,0) -- (\a+30:\L);
	  \draw[long]  (0,0) -- (\a+210:\L);
	}
	\draw[short] (0,0) -- (0:\R) node[right] {$\omega_2$};
	\draw[short] (0,0) -- (60:\R) node[above right=-1ex] {$\omega_1$};

	\draw[long]  (0,0) -- (-30:\L) node[right] {$\alpha_2$};
	\draw[long]  (0,0) -- (90:\L) node[right] {$\alpha_1$};
		
	\end{tikzpicture}
	\]
	\vspace*{-5ex}
	\caption{\small Fundamental roots and fundamental weights for the $A_2$ root system. Observe that $AG_2$ is contained in the weight lattice.}%
	\label{fig:weights}
\end{figure}
\end{remark}

Condition (ii) says that the given $AG_2$-grading is the weight-space decomposition of $L$ for the adjoint action of $\fH$.
The Weyl group of $\fG$ preserves the set $AG_2$ and hence permutes the possible weights.
Passing from an $(AG_2, A_2)$-grading to an $(AG_2,G_2)$-grading introduces the additional reflections corresponding to the short roots of $G_2$.

\bigskip

We will focus mostly on examining the more general $(AG_2, A_2)$-graded Lie algebras;
only afterwards, we will look at the implications of imposing a more restrictive $(AG_2, G_2)$-grading.
Before we proceed, we recall the definition of \emph{structurable algebras} and we introduce conjugate-supplementary idempotents in structurable algebras. These will play a key role in the whole paper.

\section{Conjugate-supplementary idempotents in structurable algebras}\label{se:idempotents}

Structurable algebras have been introduced by Bruce Allison in \cite{Allison1978}.
Recall that $k$ is an arbitrary field with $\Char(k) \neq 2,3$.
Our algebras are not necessarily commutative or associative.
\begin{definition}\label{def:struct algebra}
	Let $\A$ be a unital $k$-algebra equipped with an involution $\sigma \colon x \mapsto \overline{x}$, i.e., a $k$-linear map satisfying $\overline{\overline{x}} = x$ and  $\overline{xy}=\overline{y} .\overline{x}$ for all $x,y \in \A$.
	Let
	\[ V_{x,y}(z) := (x\overline y)z + (z\overline y)x - (z\overline x)y \]
	for all $x,y,z \in \A$.
	If
		\[ [V_{x,y},V_{z,w}]=V_{V_{x,y}(z),w}-V_{z,V_{y,x}(w)} \]
	for all $x,y,z,w\in\A$ (where the left-hand side denotes the Lie bracket of the two operators),
	then we call $\A$ a \textit{structurable algebra}.
	We will also use the notation
	\[ T_x := V_{x,1} \]
	for all $x \in \A$. Explicitly, $T_x(z) = xz + zx - z\overline{x}$ for all $x,z \in \A$.
\end{definition}

\begin{definition}
	Let $\A$ be a structurable algebra; then $\A=\mathcal H\oplus \SS$, with 
		\[ \mathcal H=\{h\in\A\mid\overline h=h\} \text{ and } \SS=\{s\in\A\mid \overline s=-s\}.\]
	The elements of $\mathcal H$ are called \textit{hermitian elements}, the elements of $\SS$ are called \textit{skew elements}. 
	The dimension of $\SS$ is called the \textit{skew-dimension} of $\A$.
\end{definition}
\begin{definition}
	Let $\A$ be a structurable algebra.
	The \emph{skewer map} is the map $\psi \colon \A \times \A \to \SS$ defined by
	\[ \psi(x, y) := x \overline{y} - y \overline{x} \]
	for all $x,y \in \A$.
\end{definition}
We will write the \emph{associator} of three elements $a,b,c \in \A$ as
\[ [a, b, c] := (ab)c - a(bc) . \]
By \cite[Proposition 1]{Allison1978}, structurable algebras are \emph{skew alternative}:
\begin{proposition}\label{pr:skew alternative}
	Let $a,b \in \A$ and $s \in \SS$. Then
	\[ [s, a, b] = - [a, s, b] = [a, b, s] . \]
\end{proposition}

We will be specifically interested in a certain kind of idempotents in $\A$.

\begin{definition}
	Let $t \in \A$ be an idempotent, i.e., $t^2 = t$; notice that also $\overline{t}$ is an idempotent.
	We call $t$ a \emph{conjugate-supplementary idempotent} if $t + \overline{t} = 1$.
\end{definition}
As we will see, many (but of course not all) structurable algebras admit such conjugate-supplementary idempotents.
\begin{lemma}\label{le:csi to skew}
	An element $t \in \A$ is a conjugate-supplementary idempotent if and only if the element $s := 2t - 1$ is a skew element with $s^2 = 1$.
\end{lemma}
\begin{proof}
	We have
	\[ s^2 = 1 \iff 4 t^2 - 4t + 1 = 1 \iff t^2 = t \]
	and
	\[ \overline{s} = -s \iff 2 \overline{t} - 1 = - 2t + 1 \iff t + \overline{t} = 1 . \qedhere  \]
\end{proof}
\begin{corollary}\label{co:assoc}
	Let $t \in \A$ be a conjugate-supplementary idempotent. Then for all $a,b \in \A$, we have
	\[ [t, a, b] = - [a, t, b] = [a, b, t] . \]
	In particular,
	\[ [t, a, t] = [t, a, \overline{t}] = [\overline{t}, a, t] = [\overline{t}, a, \overline{t}] = 0 \]
	and
	\[ t(ta) = ta, \quad (at)t = at, \quad t (\overline{t}a) = 0, \quad (a \overline{t})t = 0 . \]
\end{corollary}
\begin{proof}
	The first claim follows from \cref{pr:skew alternative,le:csi to skew}. Setting $a = t$ then gives $[t, t, a] = 0$, which we rewrite as $t(ta) = t^2 a = ta$ since $t$ is idempotent. Similarly, setting $a = \overline{t}$ gives $[t, \overline{t}, a] = 0$ and hence $t(\overline{t}a) = (t \overline{t})a = 0$ since $t \overline{t} = t(1-t) = 0$. The remaining claims are now clear.
\end{proof}
Notice that, in particular, the expressions $tat$, $ta\overline{t}$, $\overline{t}at$ and $\overline{t}a\overline{t}$ are unambiguous, and we will use them without explicitly referring to \cref{co:assoc}.

\medskip

We will need the following easy computation later.
\begin{lemma}\label{le:V in kt}
	Let $t \in \A$ be a conjugate-supplementary idempotent.
	Then
	\[ V_{t,t} = V_{\overline{t}, \overline{t}} = 0, \quad V_{t, \overline{t}} = T_t, \quad V_{\overline{t}, t} = T_{\overline{t}} . \]
	In particular, for all $x,y \in \B := k[t]$, we have $V_{x,y} = T_z$ for some $z \in \B$.
\end{lemma}
\begin{proof}
	By definition, we have $V_{t,t}(a) = (t\overline t)a + (a\overline t)t - (a\overline t)t = 0$ for all $a \in \A$.
	Similarly, $V_{\overline{t}, \overline{t}} = 0$.
	It follows that $T_t = V_{t,1} = V_{t, t+\overline{t}} = V_{t,t} + V_{t, \overline{t}} = V_{t,\overline{t}}$.
	Similarly, $T_{\overline{t}} = V_{\overline{t}, t}$.
\end{proof}

We will continue our investigation of structurable algebras equipped with a conjugate-supplementary idempotent in \cref{se:idempotents2}.

\section{From structurable algebras to \texorpdfstring{$5$}{5}-graded Lie algebras}\label{se:TKK}

The Tits--Kantor--Koecher construction, producing $3$-graded Lie algebras starting from Jordan algebras, has been generalized by Bruce Allison to a construction producing $5$-graded Lie algebras starting from structurable algebras.
In this short section, we recall this construction.
We have decided to explicitly include the definition of the Lie bracket, because we will extensively use this in \cref{se:AG2 TKK}.

Rather than referring to Allison's original papers \cite{Allison1978,Allison1979}, we refer to the more recent \cite[Chapter 2]{Boelaert2019}, which contains an extended overview of the main results and examples in the theory of structurable algebras.
(The description in \cite[\S 3]{Allison1979} is equivalent, but organized in a slightly different fashion.)

\begin{definition}[{\cite[Definition 2.4.1]{Boelaert2019}}]\label{def:TKK}
	Let $\A$ be a structurable algebra and let $\SS$ be its subspace of skew elements.
	\begin{enumerate}
		\item
			Let $\instrl(\A)$ be the subspace of $\End(\A)$ spanned by the operators $V_{x,y}$ for all $x,y \in \A$.
		\item
			For each $A \in \instrl(\A)$, we define $A^\epsilon \in \instrl(\A)$ and $A^\delta \in \End(\SS)$ by setting
			\begin{align*}
				A^\epsilon &:= A - L_{A(1)+\overline{A(1)}} , \\
				A^\delta &:= A + R_{\overline{A(1)}} ,
			\end{align*}
			where $L_x$ and $R_x$ denote left and right multiplication by $x \in \A$, respectively.
			By \cite[Eq.\@~(4)]{Allison1979},
			\begin{align*}
			    V_{x,y}^\epsilon &= -V_{y,x}, \\
			    V_{x,y}^\delta(s) &= -\psi(x,sy) ,
			\end{align*}
			for all $x,y\in \A$ and $s\in \SS$.
		\item
	Define
	\[ \TKK := \SS_- \oplus \A_- \oplus \instrl(\A) \oplus \A_+ \oplus \SS_+ \]
	equipped with the following Lie bracket, where $x,x',y,y'\in \A$, $s,s',r,r'\in \SS$ and $V_{a,b},V_{a',b'}\in \instrl(\A)$:
	\begin{alignat*}{2}
		\intertext{$\bullet\ [\instrl,K(\A)]$}
			[V_{a,b},V_{a',b'}] &= V_{ V_{a,b} a', b'}-V_{a', V_{b,a} b' } &\ \in \instrl(\A),&\\
			[V_{a,b},x_+] &:= (V_{a,b}x)_+ \in \A_+ , \quad & [V_{a,b},y_-] &:= (V_{a,b}^\epsilon y)_- \\*
			&&&\ =(-V_{b,a}y)_-\in \A_-,\\
			[V_{a,b},s_+] &:= (V_{a,b}^\delta s)_+ & [V_{a,b},r_-] &:= (V_{a,b}^{\epsilon\delta} r)_- \\*
			&\ = -\psi(a,sb)_+\in \SS_+,&&\ =\psi(b,ra)_-\in \SS_-,\\
		\intertext{$\bullet\ [\SS_\pm,\A_\pm]$}
			[s_+,x_+] &:= 0, & [r_-,y_-] &:= 0,\\
			[s_+,y_-] &:= (sy)_+\in \A_+,& [r_-,x_+] &:= (rx)_-\in \A_-,\\
		\intertext{$\bullet\  [\A_\pm,\A_\pm]$}
			[x_+,y_-] &:= V_{x,y}\in \instrl(\A),\\
			[x_+,x_+'] &:= \psi(x,x')_+\in \SS_+,&[y_-,y_-'] &:= \psi(y,y')_-\in \SS_-,\\
		\intertext{$\bullet \ [\SS_\pm,\SS_\pm]$}
			[s_+,s_+']&:=0,&[r_-,r_-']&:=0,\\
			[s_+,r_-]&:=L_{s}L_{r}\in \instrl(\A).
	\end{alignat*}
	This is a $5$-graded Lie algebra, with
	\begin{multline*}
		\TKK_{-2} = \SS_-, \quad \TKK_{-1} = \A_-, \quad \TKK_{0} = \instrl(\A), \\ \TKK_{1} = \A_+, \quad \TKK_{2} = \SS_+
	\end{multline*}
	and $\TKK_{i} = 0$ for all other $i$.
	We call it the \emph{Tits--Kantor--Koecher (TKK) Lie algebra of $\A$}.
	\end{enumerate}
\end{definition}
This TKK construction has been further generalized to \emph{Kantor pairs}. We will only need this notion briefly in the proof of \cref{th:AG2-A2} below, and it suffices to know that Kantor pairs generalize structurable algebras and admit a TKK construction generalizing the TKK construction that we have just presented.
(This is then sometimes referred to as the \emph{Kantor construction}.)
The importance of Kantor pairs is highlighted by the fact that for \emph{every} $5$-graded Lie algebra $L = \bigoplus L_i$, the pair $P = (L_1, L_{-1})$ has the structure of a Kantor pair; this is called the Kantor pair \emph{enveloped by $L$}.
 
We refer to \cite{Allison17} for an excellent introduction to the subject; see, in particular, its \S 4.3 for the Kantor construction.

\section{From \texorpdfstring{$(AG_2, A_2)$}{(AG2, A2)}-graded Lie algebras to structurable algebras}\label{se:AG2 to stralg}

We are now ready to make the connection between $(AG_2, A_2)$-graded Lie algebras and structurable algebras equipped with a conjugate-supplementary idempotent.

So let $L$ be an $(AG_2, A_2)$-graded Lie algebra with grading algebra $\fG \cong \ssl_3(k)$.
\begin{definition}\label{def:efh}
	Inside the grading algebra $\fG \cong \ssl_3(k)$, choose Chevalley generators $e_i, f_i, h_i$ for $i = 1,2$, and set
	\[ e := e_1 + e_2, \quad f := f_1 + f_2, \quad h := h_1 + h_2 . \]
	Then $\{ e, 2h, 2f \}$ forms the $\ssl_2$-triple of an $\ssl_2$-subalgebra that we call \emph{principal}; the elements $e$ and $f$ are regular nilpotent elements.
	See, e.g., \cite{McNinch05}, although we only use the terminology to distinguish this $\ssl_2$ subalgebra from the ``root $\ssl_2$ subalgebras'' spanned by the $\ssl_2$-triples $\{ e_i, h_i, f_i \}$.
\end{definition}
\begin{lemma}\label{le:grading}
	The element $h \in \fH$ defines a $5$-grading $L = L_{-2} \oplus L_{-1} \oplus L_0 \oplus L_1 \oplus L_2$, where $L_i = \{ x \in L \mid [h,x] = ix \}$ for each $i \in \{ -2,\dots,2\}$.
\end{lemma}
\begin{proof}
	Using the fact that $\omega_i(h_j) = \delta_{ij}$ for all $i,j \in \{ 1,2 \}$ (see \cref{rem:omega}), we see that $\omega_1(h) = \omega_2(h) = 1$.
	It follows that for any $\alpha \in AG_2^o$, we have $\alpha(h) = 1$ if and only if $\alpha$ lies on the line through $\omega_1$ and $\omega_2$ (see \cref{fig:weights}), which happens precisely for $\alpha_1$, $\omega_1$, $\omega_2$ and $\alpha_2$. We now see that the roots $\alpha$ with $\alpha(h) = i$ are given by lines parallel to this line. In particular, only the values $i \in \{ -2, -1, 0, 1, 2 \}$ can occur.
	Explicitly, $L_i = \oplus_{\alpha \in \Phi_i} L_\alpha$ where
	\begin{align*}
		\Phi_{-2} &= \{ -(\alpha_1+\alpha_2), -2\omega_1, -2\omega_2 \} , \\
		\Phi_{-1} &= \{ -\alpha_1, -\alpha_2, -\omega_1, -\omega_2 \} , \\
		\Phi_{0} &= \{ o, \pm(\omega_1 - \omega_2), \pm 2 (\omega_1 - \omega_2) \} , \\
		\Phi_{1} &= \{ \alpha_1, \alpha_2, \omega_1, \omega_2 \} , \\
		\Phi_{2} &= \{ \alpha_1+\alpha_2, 2\omega_1, 2\omega_2 \} .
	\end{align*}
	Notice that each $L_i$ is indeed equal to $\{ x \in L \mid [h,x] = ix \}$ because $\Char(k) \neq 2,3$, and that the relations $[L_i, L_j] \subseteq L_{i+j}$ then follow automatically from the Jacobi identity.
\end{proof}

\begin{lemma}\label{le:Tgen}
	Let $L = \bigoplus L_i$ be the grading from \cref{le:grading} and let $T := L_{-1} \oplus L_1$.
	Then $L$ is generated by $T$ as a Lie algebra.
\end{lemma}
\begin{proof}
	Let $M$ be the subalgebra of $L$ generated by $T$.
	In view of \cref{def:supported}\cref{def:supported:2}, it suffices to show that every non-zero $AG_2$-root subspace is contained in $M$.
	
	Certainly $L_{\pm 1} \subseteq M$. We now show that $L_2 \subseteq M$.
	So let $x \in L_\lambda$ with $\lambda(h) = 2$.
	Observe that for each of the three possible weights $\lambda = \alpha_1 + \alpha_2$, $\lambda = 2\omega_1$ and $\lambda = 2\omega_2$, we can choose $i \in \{ 1,2 \}$ such that $\alpha_i + \lambda \not\in AG_2^o$ and $\lambda(h_i) \in \{ 1, 2 \}$.
	Then $[e_i, x] \in L_{\alpha_i + \lambda} = 0$. It follows that
	\[ [e_i, [f_i, x]] = [[e_i, f_i], x] = [h_i, x] = \lambda(h_i) x . \]
	Since $\lambda(h_i) \neq 0$ and $[f_i, x] \in L_1$, we conclude that $x \in [L_1, L_1] \subseteq M$.
	This shows that $L_2 \subseteq M$; the argument for $L_{-2}$ is, of course, similar.

	Next, consider $x \in L_\lambda$ with
	\[ \lambda = \omega_1-\omega_2 \quad \text{or} \quad \lambda = 2(\omega_1-\omega_2). \]
	We now have $\lambda - \alpha_2 \notin AG_2^o$ and $\lambda(h_2) \in \{-1,-2\}$.
	Thus $[f_2,x]=0$, and
	\[ [f_2,[e_2,x]] = [-h_2,x] = -\lambda(h_2)x. \]
	Since $f_2 \in L_{-1}$ and $[e_2, x] \in L_1$, we now conclude that $x \in [L_{-1}, L_1] \subseteq M$.
	Similarly, we have $L_\lambda \subseteq M$ for the two remaining weights $\lambda = -(\omega_1 - \omega_2)$ and $\lambda = -2(\omega_1 - \omega_2)$.
\end{proof}

Our next goal is to show that the Kantor pair $(L_1, L_{-1})$ enveloped by $L$ is isomorphic to the Kantor pair of a structurable algebra; see \cref{th:AG2-A2} below.
A general criterion for this is given in \cite[Lemma~4.16]{Stavrova22}, but its proof relies on \cite[Corollary~15]{Allison1999}, which ultimately relies on \cite[Theorem~8]{Allison1999}.
Although only $2$ and $3$ are assumed to be invertible in \cite{Allison1999}, the proof of \cite[Theorem~8]{Allison1999} is not valid in characteristic~$5$; see also \cite[Remark~3.2.6]{Boelaert2019}.
We therefore follow a different route, which also works in characteristic~$5$.
In the spirit of the proof of \cite[Corollary~15]{Allison1999}, we use an automorphism of $L$ interchanging $L_1$ and $L_{-1}$ to turn $L_1$ into a Kantor triple system with a left unit, which then gives rise to a structurable algebra by \cite[Theorem~2.9]{Elduque14}.
The required automorphism is constructed in the following lemma, using the Weyl element of the grading algebra~$\fG$ corresponding to its highest root.

\begin{lemma}\label{le:Phi}
	There is an automorphism $\Phi$ of $L$ such that $\Phi^2 = \id$, $\Phi(L_i) = L_{-i}$ for all $i$, and $\Phi(e) = f$
	(with $e$ and $f$ as in \cref{def:efh}).
\end{lemma}
\begin{proof}
	Let $\gamma := \alpha_1 + \alpha_2$ be the highest root of $A_2$, and let $e_\gamma := [e_1, e_2] \in L_2$ and $f_\gamma := [f_2, f_1] \in L_{-2}$.
	Since $2\alpha_1 + \alpha_2$ and $\alpha_1 + 2\alpha_2$ are not roots of $A_2$, we have $[e_\gamma, e_1] = [e_\gamma, e_2] = [f_\gamma, f_1] = [f_\gamma, f_2] = 0$, and an easy computation using the Jacobi identity gives
	\[ [f_\gamma, e_1] = f_2, \quad [e_\gamma, f_2] = e_1, \quad [f_\gamma, e_2] = -f_1, \quad [e_\gamma, f_1] = -e_2 . \]
	In particular, $[e_\gamma, f_\gamma] = [[e_\gamma, f_2], f_1] + [f_2, [e_\gamma, f_1]] = h_1 + h_2 = h$, so $(e_\gamma, h, f_\gamma)$ is an $\ssl_2$-triple.

	Since $e_\gamma \in L_2$ and $f_\gamma \in L_{-2}$, the $5$-grading immediately gives $(\ad e_\gamma)^3 = (\ad f_\gamma)^3 = 0$.
	Hence, if $D = \ad e_\gamma$ or $D = \ad f_\gamma$, then the Leibniz rule gives $0 = D^3[x,y] = 3\bigl( [Dx, D^2y] + [D^2x, Dy] \bigr)$ and $0 = D^4[x,y] = 6 [D^2x, D^2y]$ for all $x,y \in L$, and since $\Char(k) \neq 2,3$, this implies that $\exp(D) := \id + D + \tfrac{1}{2}D^2$ is an automorphism of $L$, with inverse $\exp(-D)$.
	Therefore,
	\[ n := \exp(\ad e_\gamma) \exp(-\ad f_\gamma) \exp(\ad e_\gamma) \]
	is an automorphism of $L$.
	A short computation in the $\ssl_2$-triple $(e_\gamma, h, f_\gamma)$ shows that $n(h) = -h$; since $L_i = \{ x \in L \mid [h, x] = ix \}$ by \cref{le:grading}, it follows that $n(L_i) = L_{-i}$ for all~$i$.
	Moreover, the relations above give $n(e_1) = \exp(\ad e_\gamma)(e_1 - f_2) = -f_2$ and $n(e_2) = \exp(\ad e_\gamma)(e_2 + f_1) = f_1$.

	Now let $\tau$ be the linear map acting on each $L_\alpha$ as multiplication by $(-1)^{\alpha(h_1)}$.
	Since $\alpha(h_1) \in \Z$ for all $\alpha \in AG_2^o$, this is an automorphism of $L$ with $\tau^2 = \id$, and clearly $\tau(L_i) = L_i$ for all~$i$.
	Since $\gamma(h_1) = 1$, we have $\tau(e_\gamma) = -e_\gamma$ and $\tau(f_\gamma) = -f_\gamma$, so $\tau \exp(\pm \ad e_\gamma) \tau = \exp(\mp \ad e_\gamma)$ and $\tau \exp(\pm \ad f_\gamma) \tau = \exp(\mp \ad f_\gamma)$; it follows that $\tau n \tau = n^{-1}$.
	Set $\Phi := \tau n$.
	Then $\Phi$ is an automorphism of $L$ with $\Phi^2 = (\tau n \tau) n = \id$ and $\Phi(L_i) = L_{-i}$ for all~$i$.
	Finally, since $\alpha_1(h_1) = 2$ and $\alpha_2(h_1) = -1$, we get $\Phi(e) = \tau(-f_2 + f_1) = f_2 + f_1 = f$.
\end{proof}

\begin{theorem}\label{th:AG2-A2}
	Let $L$ be an $(AG_2,A_2)$-graded Lie algebra. Then:
	\begin{enumerate}
		\item 
			There exists a structurable algebra $\A$ 
			such that
			\[ L/Z(L) \cong K(\A) \]
			as $5$-graded Lie algebras.
			The isomorphism can be chosen in such a way that the composition
			\[ \pi \colon L \to L/Z(L) \to K(\A) \]
			maps the elements $e \in L_1$ and $f \in L_{-1}$ from \cref{def:efh} to $1_+ \in K(\A)_1$ and $1_- \in K(\A)_{-1}$, respectively.
		\item 
			If we write $\pi(e_1) = t_+$ with $t \in \A$, then $t$ is a conjugate-supplementary idempotent.
		\item\label{th:AG2-A2:images}
			The homomorphism $\pi$ maps
			\begin{alignat*}{3}
				e_1 &\mapsto t_+ , & \quad f_1 &\mapsto \overline{t}_- , & \quad h_1 &\mapsto T_t , \\
				e_2 &\mapsto \overline{t}_+ , & \quad f_2 &\mapsto t_- , & \quad h_2 &\mapsto T_{\overline{t}} , \\
				e &\mapsto 1_+ , & \quad f &\mapsto 1_- , & \quad h & \mapsto \id_\A .
			\end{alignat*}
		\item\label{th:AG2-A2:root spaces}
			For each $\alpha \in AG_2^o$, the homomorphism $\pi$ maps the root space $L_\alpha$ onto the weight space
			\[ K(\A)_\alpha := \{ x \in K(\A) \mid [T_t, x] = \alpha(h_1) x \text{ and } [T_{\overline{t}}, x] = \alpha(h_2) x \} . \]
			(These are precisely the root spaces of the $(AG_2, A_2)$-grading of $K(\A)$ from \cref{th:AG2-A2-conv} below.)
	\end{enumerate}
\end{theorem}
\begin{proof}
	By \cref{le:grading}, the Lie algebra $L$ is $5$-graded and $h \in L$ is a grading element for this grading (i.e., $\ad(h)$ is its grading derivation).
	As we have discussed at the very end of \cref{se:TKK}, the pair $P = (L_1, L_{-1})$ equipped with the triple product $\{ x, y, z \} := [[x, y], z]$ (for $x,z \in L_\sigma$ and $y \in L_{-\sigma}$, $\sigma = \pm 1$) is a \emph{Kantor pair}; see, e.g., \cite[p.~445]{Garcia2011} or \cite[\S 4.2]{Allison17}.

	Let $T := L_{-1} \oplus L_1$.
	By \cref{le:Tgen}, $T$ generates $L$ as a Lie algebra.
	By \cite[Lemma 4.4]{Allison17}, $L = [T, T] \oplus T$ with
	\[ [T, T] = \underbrace{[L_{-1}, L_{-1}]}_{L_{-2}} \oplus \underbrace{[L_{-1}, L_1]}_{L_0} \oplus \underbrace{[L_1, L_1]}_{L_2} . \]
	Notice that the center $Z(L)$ is contained in $L_0$ because $h \in L$ is a grading element.
	By \cite[Remark 4.6]{Allison17}, the algebra $\overline{L} := L/Z(L)$ \emph{tightly envelops} the Kantor pair~$P$.
	It follows from the uniqueness of the tight envelope (see \cite[Corollary~4.17]{Allison17}) that $\overline{L}$ is graded-isomorphic to the \emph{standard $5$-graded embedding}~$K(P)$.
	
	Notice that the canonical map $L \to \overline{L}$ restricts to isomorphisms $L_i \to \overline{L}_i$ for each $i \neq 0$, and that the Kantor pair enveloped by $\overline{L}$ is again $P = (L_1, L_{-1})$.
	We now show that $P$ is isomorphic to the Kantor pair of a structurable algebra.
	Let $\Phi$ be the automorphism of $L$ from \cref{le:Phi}, and define a triple product on $L_1$ by
	\[ \langle x, y, z \rangle := [[x, \Phi(y)], z] \quad \text{for all } x, y, z \in L_1 . \]
	Since $\Phi$ is an automorphism with $\Phi^2 = \id$, the pair of maps $(\id_{L_1}, \Phi|_{L_1})$ is an isomorphism from the double $(L_1, L_1)$ of this triple system onto the Kantor pair $P$, so $(L_1, \langle \cdot, \cdot, \cdot \rangle)$ is a Kantor triple system (see \cite[Special Cases~4.1(ii)]{Allison17}).
	Since $\langle e, e, x \rangle = [[e, f], x] = [h, x] = x$ for all $x \in L_1$, the element $e$ is a left unit for this Kantor triple system.
	By \cite[Theorem~2.9]{Elduque14}, there is a structurable algebra $\A$, with underlying vector space $L_1$ and with unit $e$, and an involutive automorphism $\sigma$ of $\A$, such that
	\[ \langle x, y, z \rangle = V_{x, \sigma(y)}(z) \quad \text{for all } x, y, z \in \A . \]
	Now consider the linear bijections $\chi_+ := \id \colon L_1 \to \A$ and $\chi_- := \sigma \circ \Phi \colon L_{-1} \to \A$.
	For $x, z \in L_1$ and $y \in L_{-1}$, we have
	\[ [[x, y], z] = \langle x, \Phi(y), z \rangle = V_{\chi_+(x), \chi_-(y)}\bigl(\chi_+(z)\bigr) , \]
	and for $y, w \in L_{-1}$ and $x \in L_1$, applying the automorphism $\Phi$ and then $\sigma$ gives
	\[ \chi_-\bigl([[y, x], w]\bigr) = \sigma\bigl( \langle \Phi(y), x, \Phi(w) \rangle \bigr) = V_{\chi_-(y), \chi_+(x)}\bigl(\chi_-(w)\bigr) . \]
	Since $[[x_+, y_-], z_+] = (V_{x,y} z)_+$ and $[[y_-, x_+], w_-] = (V_{y,x} w)_-$ in $K(\A)$ by \cref{def:TKK}, this shows that $(\chi_+, \chi_-)$ is an isomorphism from $P$ onto the Kantor pair $(\A_+, \A_-)$ enveloped by $K(\A)$.
	Both $\overline{L}$ and $K(\A)$ tightly envelop their Kantor pairs, so by \cite[Corollaries~4.17 and~4.18]{Allison17}, the isomorphism $(\chi_+, \chi_-)$ extends to an isomorphism $\overline{L} \to K(\A)$ of $5$-graded Lie algebras.
	Finally, $\chi_+(e) = e = 1$ and $\chi_-(f) = \sigma(\Phi(f)) = \sigma(e) = 1$, so this isomorphism maps $\overline{e}$ to $1_+$ and $\overline{f}$ to $1_-$.
	
	Next, consider $\pi \colon L \to \overline{L} \to K(\A)$ and write $\pi(e_1) = t_+$ for some $t \in \A$.
	Since $\pi(e) = 1_+$ by construction and $e = e_1 + e_2$, we have $\pi(e_2) = (1-t)_+$.
	Because $h_1 = [e_1, f_1] = [e_1, f]$, we have $\pi(h_1) = [t_+, 1_-]$, so by \cref{def:TKK}, we get
	\[ \pi(h_1) = V_{t,1} = T_t . \]
	Now
	\[ [h_1, e] = [h_1, e_1 + e_2] = 2e_1 - e_2 = 3 e_1 - e , \]
	and applying $\pi$ gives $[T_t, 1_+] = (3t - 1)_+$.
	By \cref{def:TKK} again, $[T_t, 1_+] = (T_t (1))_+ = (2t - \overline{t})_+$,
	and we see that $3t - 1 = 2t - \overline{t}$, i.e.,
	\[ \overline{t} = 1 - t . \]
	Similarly, $[h_1, e_1] = 2e_1$, and again applying $\pi$ and \cref{def:TKK} gives $T_t (t) = 2t$.
	Since $T_t(t) = 2t^2 - t \overline{t} = 3t^2 - t$, it follows that $3t^2 = 3t$. Since $\Char(k) \neq 3$, we conclude that
	\[ t^2 = t \]
	and thus $t$ is a conjugate-supplementary idempotent.
	
	Next, since $h_2 = [e_2, f]$, we get
	\[ \pi(h_2) = [\overline{t}_+, 1_-] = T_{\overline{t}} . \]
	Now $[h_1, f] = [h_1, f_1 + f_2] = -2 f_1 + f_2 = f - 3 f_1$, hence $[T_t, 1_-] = 1_- - 3 \pi(f_1)$.
	By \cref{def:TKK} again, we have
	\[ [T_t, 1_-] = [V_{t,1}, 1_-] = (-V_{1,t}(1))_- = (t - 2\overline{t})_- = (1 - 3\overline{t})_- , \]
	so we deduce that $\pi(f_1) = \overline{t}_-$.
	Since $f_2 = f - f_1$, we then also get $\pi(f_2) = t_-$.
	
	Finally, let $\alpha \in AG_2^o$ and $x \in L_\alpha$.
	Since $\pi(h_1) = T_t$ and $\pi(h_2) = T_{\overline{t}}$, we get $[T_t, \pi(x)] = \pi([h_1, x]) = \alpha(h_1) \pi(x)$ and, similarly, $[T_{\overline{t}}, \pi(x)] = \alpha(h_2) \pi(x)$, so $\pi(L_\alpha) \subseteq K(\A)_\alpha$.
	Since $\Char(k) \neq 2,3$, distinct $\alpha \in AG_2^o$ give distinct pairs $(\alpha(h_1), \alpha(h_2))$ in $k^2$ (both entries lie in $\{ -2, \dots, 2 \}$), so the sum of the subspaces $K(\A)_\alpha$ is direct.
	As $\pi$ is surjective and $L = \bigoplus_{\alpha \in AG_2^o} L_\alpha$, it follows that $\pi(L_\alpha) = K(\A)_\alpha$ for all $\alpha \in AG_2^o$.
\end{proof}

\section{From structurable algebras to \texorpdfstring{$AG_2$}{AG2}-algebras}\label{se:idempotents2}

We now take a fresh start, continuing where we left off after \cref{se:idempotents,se:TKK}.
So we begin with a structurable algebra $\A$ equipped with a conjugate-supplementary idempotent $t \in \A$ and we let $L = K(\A)$ be the TKK Lie algebra constructed from~$\A$.
We aim to construct an $(AG_2, A_2)$-grading on $L$ and understand the algebraic structure of the different components that arise in this grading in terms of the data $\A$ and $t$.

\medskip

Our first step is to show that conjugate-supplementary idempotents give rise to a certain Peirce decomposition of the algebra $\A$.
\begin{proposition}\label{pr:peirce}
	Let $t \in \A$ be a conjugate-supplementary idempotent. Then the left and right multiplication operators $L_t \colon a \mapsto ta$ and $R_t \colon a \mapsto at$ on $\A$ are diagonalizable with only eigenvalues $0$ and $1$, and they commute. They induce a Peirce decomposition
	\[ \A = \A_{11} \oplus \A_{10} \oplus \A_{01} \oplus \A_{00} \]
	where
	\[ \A_{ij} = \{ a \in \A \mid L_t(a) = ia, \ R_t(a) = ja \} . \]
	Explicitly, we have
	\[ \A_{11} = t \A t, \quad \A_{10} = t \A \overline{t} , \quad \A_{01} = \overline{t} \A t , \quad \A_{00} = \overline{t} \A \overline{t} . \]
\end{proposition}
\begin{proof}
	By \cref{co:assoc}, we have $L_t^2 - L_t = 0$, so $L_t$ has only eigenvalues $0$ and~$1$; the same holds for $R_t$.
	Now simply observe that because $t + \overline{t} = 1$, any $a \in \A$ can be decomposed as
	\[ a = tat + ta\overline{t} + \overline{t}at + \overline{t}a\overline{t} . \]
	The rest is now clear.
\end{proof}
It will be convenient to write this decomposition in a ``matrix form''
\begin{equation}\label{eq:matrix form}
	\A = \begin{pmatrix} \A_{11} & \A_{10} \\ \A_{01} & \A_{00} \end{pmatrix}
		 = \begin{pmatrix} t \A t & t \A \overline{t} \\ \overline{t} \A t & \overline{t} \A \overline{t} \end{pmatrix} .
\end{equation}
\begin{example}\label{ex:SMA}
	It might be instructive to see what happens in the case of the ``structurable matrix algebras'' $M(J,\eta)$; see \cite[Section~8.v]{Allison1978} and \cite[Section~4]{Allison1984}.
	We follow the description from \cite[Example 1.9]{Allison1990}.
	Let $J$ be a (possibly degenerate) cubic norm structure, i.e., $J$ is a Jordan algebra over $k$, $T \colon J \times J \to k$ is a symmetric bilinear form, and $N \colon J \to k$ is a cubic form on $J$ such that either
	\begin{enumerate}[(a)]
		\item\label{ex:SMA:1} $J$ is a non-degenerate cubic norm structure with norm $N$, trace $T$ and base point $1$ (so $N(1)=1$);
		\item\label{ex:SMA:2} $J$ is the Jordan algebra of a non-degenerate quadratic form $q$ with base point~$1$ (so $q(1)=1$), with trace $T$ and $N=0$;
		\item\label{ex:SMA:3} $J=0$ (with $T=0$ and $N=0$).
	\end{enumerate}
	Such a cubic norm structure comes equipped with a quadratic ``adjoint map''
	\[ \sharp \colon J \to J \colon v \mapsto v^\sharp \]
	satisfying $N(v+w) = N(v) + T(w, v^\sharp) + T(v, w^\sharp) + N(w)$ for all $v,w \in J$, and a ``Freudenthal cross product'' $\times \colon J \times J \to J$ given by
	\[ v \times w := (v+w)^\sharp - v^\sharp - w^\sharp . \]
	(Both $\sharp$ and $\times$ are the zero maps in cases \cref{ex:SMA:2,ex:SMA:3} listed above.)
	Now let $\eta \in k^\times$ and let
	\[ 
	\A=\left\{ 
	\begin{pmatrix}
		a & v \\ 
		w & b
	\end{pmatrix} 
	\Bigm| a,b\in k,\ v,w\in J 
	\right\},
	\]
	and define the multiplication and the involution by
	\begin{align*}
	\begin{pmatrix}
		a_1 & v_1 \\ 
		w_1 & b_1
	\end{pmatrix}
	\begin{pmatrix}
		a_2 & v_2 \\ 
		w_2 & b_2
	\end{pmatrix}
	&=
	\begin{pmatrix}
		a_1a_2+\eta T(v_1,w_2)
		&
		a_1v_2+b_2v_1+\eta(w_1\times w_2)
		\\
		a_2w_1+b_1w_2+v_1\times v_2
		&
		b_1b_2+\eta T(w_1,v_2)
	\end{pmatrix},
	\\
	\overline{
	\begin{pmatrix}
		a & v \\ 
		w & b
	\end{pmatrix}}
	&=
	\begin{pmatrix}
		b & v \\ 
		w & a
	\end{pmatrix},
	\end{align*}
	for all $a_1,a_2,b_1,b_2\in k$ and $v_1,v_2,w_1,w_2\in J$.
	Observe that
	\[ \SS = \left\{ \begin{pmatrix}
			a & 0 \\ 0 & -a
		\end{pmatrix} \Bigm| a \in k \right\} , \]
	so in particular, $\A$ has skew-dimension $1$.
	
	Notice that, by \cite[Theorem 1.13]{Allison1990} together with \cref{le:csi to skew}, any central simple structurable algebra of skew-dimension $1$ admitting a conjugate-supplementary idempotent is of this form.
	
	Now let $t$ be the idempotent $t = \idempotentmatrix \in \A$ and observe that $\overline{t} = \begin{psmallmatrix} 0 & 0 \\ 0 & 1 \end{psmallmatrix}$, so $t$ is a conjugate-supplementary idempotent.
	For any $x=\begin{psmallmatrix}
			a & v \\ w & b
		\end{psmallmatrix}\in\A$, we now compute
	\[
		txt=\begin{pmatrix} a&0\\0&0\end{pmatrix},\quad
		tx\overline{t}=\begin{pmatrix}0&v\\0&0\end{pmatrix},\quad
		\overline{t}xt=\begin{pmatrix}0&0\\w&0\end{pmatrix},\quad
		\overline{t}x\overline{t}=\begin{pmatrix}0&0\\0&b\end{pmatrix},
	\]
	so we see that the matrix form from~\eqref{eq:matrix form} is precisely the given matrix form of $M(J, \eta)$.
\end{example}
\begin{example}\label{ex:split comp}
	Another interesting class of examples are the split composition algebras of dimension $\geq 2$.
	There are three cases, which we describe as in \cite[pp.\@~157--158, \S 2.4]{McCrimmon2004}.
	\begin{enumerate}[(a)]
		\item $\A = k \oplus k$, equipped with the exchange involution. In this case, let $J := 0$, let $T \colon J \times J \to k$ be the zero map, and let $\times \colon J \times J \to J$ be the zero map.
		\item $\A = M_2(k)$, equipped with the standard (symplectic) involution. In this case, let $J := k$, let $T \colon J \times J \to k \colon (a, b) \mapsto ab$, and let $\times \colon J \times J \to J$ be the zero map.
		\item $\A = \Zorn(k)$, equipped with the standard involution. In this case, let $J := k^3$, let $T \colon J \times J \to k$ be the usual dot product
			\[ T \bigl( (a_1, a_2, a_3), \ (b_1, b_2, b_3)\bigr) = a_1 b_1 + a_2 b_2 + a_3 b_3 \]
			and let $\times \colon J \times J \to J$ be the usual cross product on $k^3$, i.e.,
			\[ (a_1, a_2, a_3) \times (b_1, b_2, b_3) = (a_2 b_3 - a_3 b_2, \ a_3 b_1 - a_1 b_3, \ a_1 b_2 - a_2 b_1) \]
			for all $a_i, b_i \in k$.
	\end{enumerate}
	We can now describe all three cases simultaneously in the ``Zorn matrix form''
	\[ \A = \begin{pmatrix} k & J \\ J & k \end{pmatrix} , \]
	equipped with multiplication and involution
	\begin{align*}
		\begin{pmatrix} a_1 & v_1 \\ w_1 & b_1 \end{pmatrix} \begin{pmatrix} a_2 & v_2 \\ w_2 & b_2 \end{pmatrix}
		&= \begin{pmatrix} a_1 a_2 + T(v_1, w_2) & a_1 v_2 + b_2 v_1 - w_1 \times w_2 \\ a_2 w_1 + b_1 w_2 + v_1 \times v_2 & b_1 b_2 + T(v_2, w_1) \end{pmatrix} , \\
		\overline{\begin{pmatrix} a & v \\ w & b \end{pmatrix}}
		&= \begin{pmatrix} b & -v \\ -w & a \end{pmatrix},
	\end{align*}
	exactly as in \cite[p.\@~158, Example 2.4.2]{McCrimmon2004}.
	Let $t = \idempotentmatrix$; then $t$ is indeed a conjugate-supplementary idempotent. For any element $x = \begin{psmallmatrix}
				a & v \\ w & b
			\end{psmallmatrix} \in \A$, we now get
	\[ txt = \begin{psmallmatrix} a & 0 \\ 0 & 0 \end{psmallmatrix}, \quad
		tx\overline{t} = \begin{psmallmatrix} 0 & v \\ 0 & 0 \end{psmallmatrix}, \quad
		\overline{t}xt = \begin{psmallmatrix} 0 & 0 \\ w & 0 \end{psmallmatrix}, \quad
		\overline{t}x\overline{t} = \begin{psmallmatrix} 0 & 0 \\ 0 & b \end{psmallmatrix},
	\]
	so we see that the matrix form from~\eqref{eq:matrix form} precisely matches our Zorn matrix description.
	
	Notice that when $\A = k \oplus k$ (i.e., $J=0$), we recover \cref{ex:SMA} with $J=0$, but in the other cases, the examples are different; this can be seen most easily by comparing the skew-dimension of $\A$, which is always $1$ in \cref{ex:SMA}, but equal to $\dim \A - 1 = 1 + 2 \dim J$ in our current example.
\end{example}
These examples already give some motivation to also write our Peirce decomposition in the matrix form~\eqref{eq:matrix form}, but the next result further strengthens this.
Notice that when we write an element of $\A$ in a matrix form, this implicitly always refers to the description from~\eqref{eq:matrix form}.
\begin{proposition}\label{pr:Peirce multiplication}
	The multiplication in $\A$ with respect to the Peirce decomposition~\eqref{eq:matrix form} is given by
	\[	\begin{pmatrix}
		a & b \\ c & d
	\end{pmatrix} \begin{pmatrix}
		e & f \\ g & h
	\end{pmatrix} = \begin{pmatrix}
		ae + bg & af + bh + cg \\ ce + dg + bf & cf + dh
	\end{pmatrix} . \]
	In particular, this gives precise multiplication rules (also called a ``fusion law'') between the Peirce subspaces $\A_{ij}$.
\end{proposition}
\begin{proof}
	Assume that $x \in \A_{ij}$ and $y \in \A_{\ell m}$, with $i,j,\ell,m \in \{ 0,1 \}$.
	By \cref{co:assoc}, we have $[t, x, y] = - [x, t, y] = [x, y, t]$. Using $tx = ix$, $xt = jx$, $ty = \ell y$ and $yt = my$, this gives
	\[ i.xy - t(xy) = \ell.xy - j.xy = (xy)t - m.xy , \]
	hence $t(xy) = (i+j-\ell)xy$ and $(xy)t = (\ell+m-j)xy$.
	Therefore, $xy \in \A_{i+j-\ell, \ell+m-j}$.
	In particular, if $i+j-\ell$ or $\ell+m-j$ is different from $0$ or $1$, then $xy = 0$;
	this implies that $ag = ah = be = ch = de = df = 0$. All other products $xy$ belong to the indicated Peirce subspaces.
\end{proof}
We now come to our first structural result.
Recall that an algebra $B$ is called \emph{alternative} if $[x,x,y] = 0$ and $[y,x,x] = 0$ for all $x,y \in B$.
\begin{proposition}\label{pr:B+Bop}
	The subalgebra $B = t\A t$ of $\A$ is an alternative algebra (with unit $t$). Moreover, $\overline{B} = \overline{t} \A \overline{t} \cong B^\op$, the opposite algebra of $B$, and $t \A t \oplus \overline{t} \A \overline{t}$ is a subalgebra of $\A$ isomorphic to $B \oplus B^\op$.
\end{proposition}
\begin{proof}
	Observe that by \cref{pr:Peirce multiplication}, $B = t\A t = \A_{11}$, $B' = \overline{t} \A \overline{t} = \A_{00}$ and $B \oplus B'$ are indeed subalgebras of~$\A$.
	Moreover, the involution of $\A$ maps $t \A t$ to $\overline{t} \A \overline{t}$ and hence induces an isomorphism $B' \cong B^\op$.
	
	Next, notice that any element $b \in B$ corresponds to a skew element
	\begin{equation}\label{eq:s(b)}
		s(b) = \begin{pmatrix} b \\ & -\overline{b} \end{pmatrix} \in \A .
	\end{equation}
	(The map $s \colon B \to \A \colon b \mapsto s(b)$ is \emph{not} an algebra homomorphism, however.)
	By~\cref{pr:skew alternative}, the elements $s(b)$ satisfy alternativity, i.e., $[s(b), s(b), s(c)] = 0$ and $[s(c), s(b), s(b)] = 0$ for all $b,c \in B$.
	Comparing the $\A_{11}$-components of both sides then gives $[b, b, c] = 0$ and $[c, b, b] = 0$ for all $b,c \in B$.
\end{proof}
\begin{remark}
	The appearance of an alternative algebra is not surprising, and is essentially due to the grading subalgebra of type $A_2$.
	Compare, for instance, with \cite[Proposition (3.17)]{Berman1992}, dealing with Lie algebras graded by $A_2$ in a stronger sense than our \cref{def:AG2-Delta-graded}; see \cref{rem:omega}. The situation described in \cite{Berman1992} corresponds to a structurable algebra of the form $\A = B \oplus B^\op$ with the exchange involution---the ``feeble case'', see \cref{def:AG2-algebra}, \cref{rem:feeble etc} and \cref{th:classification}\cref{th:classification:feeble}---which we have encountered as a subalgebra in \cref{pr:B+Bop}.
\end{remark}

It will be convenient to always write the elements of $B^\op$ as $\overline{b}$ for $b \in B$.
The multiplication in $\A$ now gives $J := \A_{10}$ and $J' := \A_{01}$ the structure of a kind of module over $B$, in the following sense.
\begin{definition}\label{def:bimodules}
	\begin{enumerate}
		\item 
			Let $J$ be a $k$-vector space equipped with an involutory $k$-linear map $v \mapsto \overline{v}$.
			Assume, moreover, that there are $k$-bilinear maps
			\begin{align*}
				B \times J \to J &\colon (b, v) \mapsto bv \quad \text{and} \\
				J \times B^\op \to J &\colon (v, \overline{b}) \mapsto v\overline{b} .
			\end{align*}
			We call $J$ an \emph{involutory $B$-$B^\op$-bimodule} if $v\overline{b} = \overline{b \overline{v}}$ for all $b \in B$ and $v \in J$.
		\item 
			Applying (i) to $B^\op$ instead of $B$, and identifying $(B^\op)^\op$ with $B$ via $\overline{\overline{b}} = b$, we obtain the notion of an \emph{involutory $B^\op$-$B$-bimodule}.
			(Explicitly, this is a $k$-vector space $J'$ with an involutory $k$-linear map $v' \mapsto \overline{v'}$ and $k$-bilinear maps $(\overline{b}, v') \mapsto \overline{b}v'$ and $(v', b) \mapsto v'b$ such that $v'b = \overline{\overline{b}\,\overline{v'}}$ for all $b \in B$ and $v' \in J'$.)
	\end{enumerate}
\end{definition}

\begin{proposition}\label{pr:bimodules}
	The subspaces $J := \A_{10}$ and $J' := \A_{01}$ form an involutory $B$\dash $B^\op$\dash bimodule and $B^\op$-$B$-bimodule, respectively, where the module maps are simply given by the corresponding multiplication inside the algebra $\A$.
\end{proposition}

\begin{proof}
	This is immediately clear from \cref{pr:Peirce multiplication}.
\end{proof}

\begin{definition}\label{def:maps}
	\begin{enumerate}
		\item 
			Define a map $\times \colon J \times J \to J'$ as the restriction of the multiplication in $\A$ to $J \times J$, and observe that by \cref{pr:Peirce multiplication}, the product of two elements in $J$ indeed lands in $J'$.
			Similarly, we define a map $\times' \colon J' \times J' \to J$ by restricting the multiplication in $\A$ to $J' \times J'$.
		\item 
			Define a map $T \colon J \times J' \to B$ as the restriction of the multiplication in $\A$ to $J \times J'$, and observe that by \cref{pr:Peirce multiplication}, the product of an element in $J$ with an element in $J'$ indeed lands in $B$.
			Similarly, we define a map $T' \colon J' \times J \to B^\op$ by restricting the multiplication in $\A$ to $J' \times J$.
	\end{enumerate}
\end{definition}

\begin{lemma}\label{le:T'}
	We have ${T'(v', v)} = \overline{T(\overline{v}, \overline{v'})}$ for all $v \in J$ and $v' \in J'$.
\end{lemma}

\begin{proof}
	Since both $T$ and $T'$ are simply given by multiplication in $\A$, this follows immediately from the fact that $x \mapsto \overline{x}$ is an involution of $\A$.
\end{proof}

Using these notations, the multiplication from \cref{pr:Peirce multiplication} takes the following form.

\begin{proposition}\label{pr:Peirce multiplication new}
	The involution on $\A$ is given by
	\[
	\overline{\begin{pmatrix}
		c & v \\ v' & \overline{d}
	\end{pmatrix}} = \begin{pmatrix}
		d & \overline{v} \\ \overline{v'} & \overline{c}
	\end{pmatrix}
	\]
	and the multiplication in $\A$ is given by
	\[
	\begin{pmatrix}
		c & v \\ v' & \overline{d}
	\end{pmatrix}
	\begin{pmatrix}
		p & w \\ w' & \overline{q}
	\end{pmatrix}
	=
	\begin{pmatrix}
		cp + T(v, w') & cw + v\overline{q} + v' \times' w' \\
		v'p + \overline{d}w' + v \times w & \overline{T(\overline{w}, \overline{v'}) + qd}
	\end{pmatrix} .
	\]
\end{proposition}

\begin{definition}\label{def:AG2-algebra}
	\begin{enumerate}
		\item 
			Let $B$ be an alternative algebra over $k$ and let $b \mapsto \overline{b}$ be a fixed isomorphism from $B$ to $B^\op$.
			Let $J$ be an involutory $B$-$B^\op$-bimodule and let $J'$ be an involutory $B^\op$-$B$-bimodule.
			Let $\times \colon J \times J \to J'$, $\times' \colon J' \times J' \to J$ and $T \colon J \times J' \to B$ be $k$-bilinear maps.
			
			We call the tuple $(B, J, J', \times, \times', T)$ an \emph{$AG_2$-algebra} if the algebra
			\[
			\A := \begin{pmatrix}
				B & J \\ J' & \overline{B}
			\end{pmatrix} ,
			\]
			with involution and multiplication given by \cref{pr:Peirce multiplication new}, is a structurable algebra.
		\item 
			We call an $AG_2$-algebra \emph{strong} if both $J$ and $J'$ contain non-zero skew elements;
			we call it \emph{weak} if $J \cap \SS = J' \cap \SS = 0$.
			If $J = J' = 0$, we call the $AG_2$-algebra \emph{feeble}.
			
			An $AG_2$-algebra is called \emph{mighty} if it is strong and, in addition, the maps $\times$ and $\times'$ are not identically zero.
		\item\label{it:AG2-algebra:associated}
			If $\A$ is a structurable algebra and $t \in \A$ is a conjugate-supplementary idempotent, then the alternative algebra $B$ defined in \cref{pr:B+Bop}, the bimodules $J$ and $J'$ defined in \cref{pr:bimodules}, and the maps $\times$, $\times'$ and $T$ introduced in \cref{def:maps} define an $AG_2$-algebra, which we refer to as the $AG_2$-algebra \emph{associated with} $\A$ and $t$.
	\end{enumerate}
\end{definition}

\begin{remark}
	It is possible for an $AG_2$-algebra to be neither weak nor strong. For example, let $\A = \langle t, \overline{t}, j, j' \rangle$ with $B = k$, $J = kj$ with $\overline{j} = j$, $J' = kj'$ with $\overline{j'} = -j'$, with scalar module actions, and with $\times$, $\times'$ and $T$ all zero. Then this algebra $\A$ is structurable, but $J'$ contains non-zero skew elements while $J$ does not.
	
	However, we will see that when $\A$ is a finite-dimensional \emph{central simple} structurable algebra (with $\Char(k) \neq 5$), this intermediate case cannot occur; see \cref{th:classification}.
\end{remark}

\begin{remark}
%
	The defining identity of a structurable algebra is $5$-linear, so a direct translation into identities for the maps $\times$, $\times'$ and $T$ would produce $4^5=1024$ component identities. Using Allison's characterization of structurable algebras by skew-alternativity together with two conditions on hermitian elements \cite[Theorem~13]{Allison1978}, this can be reduced to a rather short list, which we have computed with AI assistance, and which we have reproduced in \cref{se:appendix}.
\end{remark}

\section{From \texorpdfstring{$AG_2$}{AG2}-algebras to \texorpdfstring{$AG_2$}{AG2}-graded Lie algebras}\label{se:AG2 TKK}

Now that we have analyzed the internal structure of the structurable algebra~$\A$ induced by the existence of the conjugate-supplementary idempotent $t \in \A$, we are ready to investigate the corresponding TKK Lie algebra $K(\A)$.
\begin{lemma}\label{le:split A2}
	Let $\A$ be a structurable algebra with conjugate-supplementary idempotent $t \in \A$ and consider the subalgebra $\B := k[t]$. Then:
	\begin{enumerate}
		\item $\B \cong k \oplus k$ equipped with the exchange involution $(x,y) \mapsto (y,x)$;
		\item $\B$ has TKK Lie algebra $K(\B) \cong \ssl_3(k)$;
		\item $K(\B)$ embeds into $\TKK$.
	\end{enumerate}
\end{lemma}
\begin{proof}
	\begin{enumerate}[leftmargin=1.5em]
		\item 
			Since $t^2 = t$ and $\overline{t} = 1-t$, we have $\B = k[t] \cong k[x] / (x^2 - x) \cong k \oplus k$.
			(Note that $t \not\in k$.)
			The involution of $\A$ stabilizes $\B$ and maps $t = (1,0)$ to $1-t = (0,1)$, so its restriction to $\B$ corresponds to the exchange involution on $k \oplus k$.
		\item
			Notice that the skew elements of $\B \cong k \oplus k$ are the multiples of $t - \overline{t}$.
			Observe that $\B$ has no non-zero derivations, so $\instrl(\B) = T_\B = \langle T_t, T_{\overline{t}} \rangle$, where $T_x := V_{x,1}$. (See, for example, \cite[\S 1]{Allison1979}.)

			Now let $e_i, f_i, h_i$ ($i = 1,2$) be Chevalley generators of $\ssl_3(k)$ as in \cref{def:efh}.
			Since $e_1, e_2, [e_1, e_2], f_1, f_2, [f_2, f_1], h_1, h_2$ form a basis of $\ssl_3(k)$ and $t_\pm, \overline{t}_\pm, (t - \overline{t})_\pm, T_t, T_{\overline{t}}$ form a basis of $K(\B)$, we can define a linear isomorphism $\ssl_3(k) \to K(\B)$ by
			\begin{alignat*}{3}
				e_1 &\mapsto t_+ , & \quad f_1 &\mapsto \overline{t}_- , & \quad h_1 &\mapsto T_t , \\
				e_2 &\mapsto \overline{t}_+ , & \quad f_2 &\mapsto t_- , & \quad h_2 &\mapsto T_{\overline{t}} , \\
				[e_1, e_2] &\mapsto (t - \overline{t})_+ , & \quad [f_2, f_1] &\mapsto (t - \overline{t})_- .
			\end{alignat*}
			Using \cref{def:TKK}, it is now easy to check that this is an isomorphism of Lie algebras.
			(Observe that it is given by the same formulas as the map $\pi$ in \cref{th:AG2-A2}\cref{th:AG2-A2:images}.)
		\item
			Notice that this is not automatic, because the TKK construction $\A \mapsto \TKK$ is not functorial.
			However, in our situation, $\instrl(\B) = T_\B$, and the map $\B \to T_\B \colon x \mapsto T_x$ is an isomorphism of vector spaces.
			(Indeed, if $T_x = 0$, then $\overline{x} = 2x$ implies $x = 4x$ so $x = 0$.)
			Similarly, the map $\A \to T_\A \colon x \mapsto T_x$ is an isomorphism of vector spaces, so the inclusion $\B \into \A$ induces an injective map
			\[ \iota \colon \instrl(\B) = T_\B \into T_\A \into \instrl(\A) \]
			and therefore also an inclusion $K(\B) \into K(\A)$ \emph{as $5$-graded vector spaces}.
			
			To verify that this inclusion preserves the Lie bracket, it remains to check that for all $x,y \in \B$, the Lie bracket $[x_+, y_-]$ computed in $K(\B)$ is mapped by $\iota$ to the Lie bracket $[x_+, y_-]$ computed in $\TKK$. In other words, we have to show that
			\[ \iota(V_{x,y}^\B) = V_{x,y}^\A \]
			for all $x,y \in \B$. By \cref{le:V in kt}, we can rewrite $V_{x,y} = T_z$ for some $z \in \B$, and this equality is valid not only in $\B$ \emph{but also in $\A$}, so we get the required result.
		\qedhere
	\end{enumerate}
\end{proof}

Our next step is to use elements of this $\ssl_3$-subalgebra to construct new gradings of the Lie algebra $L = \TKK$.
Notice that $T_1 = V_{1,1}$ is the identity operator on~$\A$, so by definition of the Lie bracket of $L$, $\ad(T_1)$ is the grading derivation of the given grading. Since $1 = t + \overline{t}$, the element $T_1 = T_t + T_{\overline{t}}$ corresponds to $h_1 + h_2$ under the isomorphism in the proof of \cref{le:split A2}, so the natural choice for finding a new grading is to use the element $T_t = V_{t,1}$ corresponding to $h_1$. (We could also use the element $T_{\overline{t}}$ corresponding to $h_2$, but since $T_{\overline{t}} = T_1 - T_t$, this would give nothing new.)
\begin{proposition}\label{pr:5by5}
	Let $\A$ be a structurable algebra with conjugate-supplementary idempotent $t \in \A$ and let $L = \TKK$.
	Let $T_t = V_{t,1} \in L_0$.
	Then $\ad(T_t)$ is diagonalizable, stabilizes each $L_i$, and has eigenvalues $-2, -1, 0, 1, 2$ only. If we write $L'_i$ for the $i$-eigenspaces of $\ad(T_t)$, then we get a $(5 \times 5)$-grading
	\[ L = \bigoplus_{i,j = -2}^2 (L_i \cap L'_j) . \]
\end{proposition}
\begin{proof}
	Observe that by definition, $T_t(z) = tz + zt - z\overline{t}$ for all $z \in \A$. Together with the fact that $\overline{t} = 1-t$, we get
	\[ T_t = L_t + 2 R_t - \id . \]
	By \cref{pr:peirce}, therefore, $T_t$ acts on each subspace $\A_{ij}$ as multiplication by $i + 2j - 1$.
	In fact, since $\Char(k) \neq 2,3$, each $\A_{ij}$ is precisely the eigenspace of $T_t$ with eigenvalue $i + 2j - 1$.
	(Notice that $i$ and $j$ are $0$ or $1$, so $i + 2j - 1$ takes the values $2, 1, 0, -1$ only.)
	
	In order to determine the eigenspaces of $\ad(T_t)$ in its action on $L = \TKK$, we have to use \cref{def:TKK}, so in particular, we will need
	\begin{align*}
		T_t^\epsilon &= T_t - L_{T_t(1) + \overline{T_t(1)}} = T_t - \id, \\
		T_t^\delta &= T_t + R_{\overline{T_t(1)}} = (L_t + 2 R_{t} - \id) + (2 \id - 3 R_t) = L_t - R_t + \id , \\
		T_t^{\epsilon\delta} &= (T_t - \id)^\delta = T_t^\delta - 2 \id = L_t - R_t - \id .
	\end{align*} 
	(Observe that $L_t - R_t$ indeed stabilizes the subspace $\SS$ of skew elements.)
	The action of $\ad(T_t)$ on $L_{-2}$, $L_{-1}$, $L_1$ and $L_2$ is now given by applying $T_t^{\epsilon\delta}$ on $\SS$, $T_t^{\epsilon}$ on $\A$, $T_t$ on $\A$ and $T_t^\delta$ on $\SS$, respectively.
	We again use \cref{pr:peirce} to see that $L_t - R_t$ acts on each $\A_{ij}$ as multiplication by $i - j$.
	To summarize, we already get the following decompositions for the action of $\ad(T_t)$ on each $L_i$ for $i \neq 0$, where $\ad(T_t)$ acts on each $L'_j$ as multiplication by~$j$.
	We use the notation $\SS_{00 \oplus 11} := \SS \cap (\A_{00} \oplus \A_{11})$, etc.
	(Note that $\SS_{00 \oplus 11} \neq \SS_{00} \oplus \SS_{11}$!)
	\begingroup
	\renewcommand{\arraystretch}{1.3}
	\[
	\begin{array}{l|ccccc}
	& L_{-2} & L_{-1} & L_{0} & L_{1} & L_{2} \\
	\hline
	L'_{2} & & & & (\A_{11})_+ & (\SS_{10})_+ \\
	L'_{1} & & (\A_{11})_- & & (\A_{01})_+ & (\SS_{00 \oplus 11})_+ \\
	L'_{0} & (\SS_{10})_- & (\A_{01})_- & & (\A_{10})_+ & (\SS_{01})_+ \\
	L'_{-1} & (\SS_{00 \oplus 11})_- & (\A_{10})_- & & (\A_{00})_+ & \\
	L'_{-2} & (\SS_{01})_- & (\A_{00})_- & & &
	\end{array}
	\]
	\endgroup
	To fill in the missing information for the action of $\ad(T_t)$ on $L_0$, we recall that $L_0$ is spanned by the $V_{x,y} = [x_+, y_-]$ for all $x,y \in \A$. Using the action of $\ad(T_t)$ on $L_{-1}$ and $L_1$, we can already conclude that $L_0$ is contained in $\oplus_{j=-3}^3 L'_j$, so it remains to show that
	\begin{align*}
		&V_{x,y} = 0 \quad \text{for all } x,y \in \A_{11} , \\
		&V_{x,y} = 0 \quad \text{for all } x,y \in \A_{00} .
	\end{align*}
	So let $x,y \in \A_{11} = t \A t$ and let $z \in \A$ be arbitrary. Write
	\[ V_{x,y}(z) = (x \overline{y}) z + (z \overline{y})x - (z \overline{x})y . \]
	By \cref{pr:Peirce multiplication}, we have $x \overline{y} = 0$. Moreover, if we write $z = \begin{psmallmatrix} a & b \\ c & d \end{psmallmatrix}$, then we see that $(zh)e = 0$ for all $h \in \A_{00}$ and all $e \in \A_{11}$, so also the two remaining terms $(z \overline{y})x$ and $(z \overline{x})y$ are~$0$.
	Similarly, we have $V_{x,y} = 0$ for all $x,y \in \A_{00}$.
	
	We conclude that for all $x,y \in \A$, the element $V_{x,y} = [x_+, y_-] \in L_0$ is contained in $\oplus_{j=-2}^2 L'_j$.
	This shows that we get a new $5$-grading $L = L'_{-2} \oplus L'_{-1} \oplus L'_{0} \oplus L'_{1} \oplus L'_{2}$, which is compatible with the original grading in the sense that
	\[ L = \bigoplus_{i,j = -2}^2 (L_i \cap L'_j) . \qedhere \]
\end{proof}

\medskip

Our next step is to understand the different components of the $L_0$-part in terms of $\A$ and $t$ and its associated $AG_2$-algebra (see \cref{def:AG2-algebra}\cref{it:AG2-algebra:associated}).
We use the notation $J = \A_{10} = t \A \overline{t}$ and $J' = \A_{01} = \overline{t} \A t$ from \cref{se:idempotents2}.

We begin with a lemma, showing how to rewrite all $V$\dash operators in $L_0 \cap L'_i$ with $i \neq 0$ in the form $V_{*, t}$ or $V_{*, \overline{t}}$.
Recall the notation $\varphi^\epsilon$ for $\varphi \in \instrl(\A)$ from \cref{def:TKK}.
\begin{lemma}\label{le:rewrite V}
	Let $\A$, $t \in \A$ and $L = \TKK$ with its $(5 \times 5)$-grading as in \cref{pr:5by5}.
	Then:
	\begin{enumerate}[itemsep=.5ex]
		\item\label{it:V_tbart} $T_t = V_{t, \overline{t}}$ and $V_{t, \overline{t}} + V_{\overline{t}, t} = \id$.
		\item\label{it:iphi} For each $\varphi \in L_0 \cap L'_i$, we have $i \varphi = - V_{\varphi(t), \overline{t}} - V_{t, \varphi^\epsilon(\overline{t})} = V_{\varphi(\overline{t}), t} + V_{\overline{t}, \varphi^\epsilon(t)}$.
		\item $V_{v, a} = V_{v\overline{a}, t}$ for all $v \in J, a \in B$.
		\item $V_{v', \overline{a}} = V_{v' a, \overline{t}}$ for all $v' \in J', \overline{a} \in B^\op$.
		\item $V_{v', a} = V_{\overline{a}v', t} = - V_{\overline{v'}a, t} = \tfrac{1}{2} V_{\overline{a} v' - \overline{v'} a, t}$ for all $v' \in J', a \in B$.
		\item\label{it:V v bara} $V_{v, \overline{a}} = V_{av, \overline{t}} = - V_{\overline{av}, \overline{t}} = \tfrac{1}{2} V_{av - \overline{av}, \overline{t}}$ for all $v \in J, \overline{a} \in B^\op$.
		\item\label{it:V JJ} $V_{v, w} = V_{v\overline{w} + \overline{w}v - \overline{v}w, \overline{t}}$ for all $v,w \in J$.
		\item\label{it:V J'J'} $V_{v', w'} = V_{v'\overline{w'} + \overline{w'}v' - \overline{v'}w', t}$ for all $v',w' \in J'$.
	\end{enumerate}
\end{lemma}
\begin{proof}
	\begin{enumerate}[leftmargin=1.5em]
		\item
			This was shown already in \cref{le:V in kt}.
		\item
			Let $\varphi \in L_0 \cap L'_i$. By \cref{pr:5by5} together with \cref{it:V_tbart}, we have $[\varphi, V_{t, \overline{t}}] = - [T_t, \varphi] = -i \varphi$. On the other hand, since $\varphi \in \instrl(\A)$, we have
			\[ [\varphi, V_{t, \overline{t}}] = V_{\varphi(t), \overline{t}} + V_{t, \varphi^\epsilon(\overline{t})} . \]
			By \cref{it:V_tbart} again, we also have
			\[ [\varphi, V_{t, \overline{t}}] = - [\varphi, V_{\overline{t}, t}] = - V_{\varphi(\overline{t}), t} - V_{\overline{t}, \varphi^\epsilon(t)} . \]
		\item
			Let $v \in J, a \in B$ and let $\varphi = V_{v,a}$; then $\varphi \in L_0 \cap L'_1$.
			Since $t_+ \in L_1 \cap L'_2$, we have $[\varphi, t_+] \in L_1 \cap L'_3 = 0$, hence $\varphi(t) = 0$.
			Similarly, since $t_- \in L_{-1} \cap L'_1$, we have $[\varphi, t_-] \in L_{-1} \cap L'_2 = 0$, hence $\varphi^\epsilon(t) = 0$.
			By~\cref{it:iphi}, therefore, $\varphi = V_{\varphi(\overline{t}), t}$.
			Now $\varphi(\overline{t}) = \varphi(1) = V_{v,a}(1) = v \overline{a}$, and we conclude that $V_{v,a} = V_{v\overline{a}, t}$.
		\item
			This is completely similar to the previous item.
		\item
			Let $v' \in J', a \in B$ and let $\varphi = V_{v',a}$; then $\varphi \in L_0 \cap L'_2$.
			Again, $\varphi(t) = 0$ and $\varphi^\epsilon(t) = 0$.
			By~\cref{it:iphi}, therefore, $2 \varphi = V_{\varphi(\overline{t}), t}$.
			Now $\varphi(\overline{t}) = \varphi(1) = V_{v',a}(1) = \overline{a} v' - \overline{v'} a$, hence
			\begin{equation}\label{eq:V_v'a}
				V_{v',a} = \tfrac{1}{2} V_{\overline{a} v' - \overline{v'} a, t} .
			\end{equation}
			Setting $a = t$ in \cref{eq:V_v'a} now gives $V_{v' + \overline{v'}, t} = 0$.
			Replacing $v'$ by $\overline{a} v'$ in this equation gives $V_{\overline{a} v' + \overline{v'} a, t} = 0$, so combined with \cref{eq:V_v'a} again, we get $V_{v',a} = V_{\overline{a}v', t} = - V_{\overline{v'}a, t}$.
		\item
			This is completely similar to the previous item.
		\item
			Let $v,w \in J$ and let $\varphi = V_{v,w}$; then $\varphi \in L_0 \cap L'_{-1}$.
			We have $\varphi(\overline{t}) = 0$ and $\varphi^\epsilon(\overline{t}) = 0$, so by \cref{it:iphi}, $\varphi = V_{\varphi(t), \overline{t}}$.
			Now $\varphi(t) = \varphi(1) = V_{v,w}(1) = v\overline{w} + \overline{w}v - \overline{v}w$, so the result follows.
		\item
			This is completely similar to the previous item.
		\qedhere
	\end{enumerate}
\end{proof}

\begin{proposition}
	Let $\A$, $t \in \A$ and $L = \TKK$ with its $(5 \times 5)$-grading as in \cref{pr:5by5}.
	Then:
	\begin{align*}
		&L_0 \cap L'_{-2} = \{ V_{v, \overline{t}} \mid v \in J \cap \SS \} , \\
		&L_0 \cap L'_{-1} = \{ V_{v', \overline{t}} \mid v' \in J' \} , \\
		&L_0 \cap L'_{1} = \{ V_{v, t} \mid v \in J \} , \\
		&L_0 \cap L'_{2} = \{ V_{v', t} \mid v' \in J' \cap \SS \} .
	\end{align*}
\end{proposition}
\begin{proof}
	By construction, $L_0 = \instrl(\A) = \langle V_{x,y} \mid x,y \in \A \rangle$, so by the grading from \cref{pr:5by5}, we have
	\begin{align*}
		&L_0 \cap L'_{-2} = \langle V_{v, \overline{a}} \mid v \in J, \overline{a} \in B^\op \rangle + \langle V_{\overline{a}, v} \mid v \in J, \overline{a} \in B^\op \rangle, \\
		&L_0 \cap L'_{-1} = \begin{multlined}[t]
			\langle V_{v', \overline{a}} \mid v' \in J', \overline{a} \in B^\op \rangle + \langle V_{\overline{a}, v'} \mid v' \in J', \overline{a} \in B^\op \rangle \\
			+ \langle V_{v,w} \mid v,w \in J \rangle
			\end{multlined} \\
		&L_0 \cap L'_{1} =
			\langle V_{v, a} \mid v \in J, a \in B \rangle + \langle V_{a, v} \mid v \in J, a \in B \rangle
			+ \langle V_{v',w'} \mid v',w' \in J' \rangle \\
		&L_0 \cap L'_{2} = \langle V_{v', a} \mid v' \in J', a \in B \rangle + \langle V_{a, v'} \mid v' \in J', a \in B \rangle.
	\end{align*}
	We now use the property
	\begin{equation}\label{eq:Vxy+Vyx}
		V_{x,y} + V_{y,x} = T_{x \overline{y} + y \overline{x}}
	\end{equation}
	for all $x,y \in \A$.
	In particular, $V_{y,x} = - V_{x,y}$ whenever $x \overline{y} = 0$, so
	\[ V_{\overline{a}, v} = - V_{v, \overline{a}} \quad \text{and} \quad V_{a, v'} = - V_{v', a}  \]
	for all $a \in B, v \in J, v' \in J'$.
	Next, notice that by \cref{le:rewrite V}\cref{it:V v bara}, we have $V_{v + \overline{v}, \overline{t}} = 0$ for all $v \in J$, hence
	$T_{v + \overline{v}} = V_{v + \overline{v}, t}$ for all $v \in J$.
	By~\cref{eq:Vxy+Vyx}, we then get $V_{v,a} + V_{a,v} = V_{v \overline{a} + a \overline{v}, t}$ for all $v \in J, a \in B$.
	Since $v \overline{a} + a \overline{v} \in J$ and $t \in B$, therefore,
	\[ \langle V_{a, v} \mid v \in J, a \in B \rangle \subseteq \langle V_{v, a} \mid v \in J, a \in B \rangle . \]
	Similarly,
	\[ \langle V_{\overline{a}, v'} \mid v' \in J', \overline{a} \in B^\op \rangle \subseteq \langle V_{v', \overline{a}} \mid v' \in J', \overline{a} \in B^\op \rangle . \]
	The result now follows from \cref{le:rewrite V}.
\end{proof}
Notice that this result shows that the components $L_0 \cap L'_i$ with $i \neq 0$ can be parametrized by $J \cap \SS$, $J'$, $J$ and $J' \cap \SS$, respectively. Indeed, we have $V_{v,t}(1) = v$ for all $v \in J$ and $V_{v', t}(1) = 2v'$ for all $v' \in J' \cap \SS$, and similarly $V_{v',\overline{t}}(1) = v'$ for all $v' \in J'$ and $V_{v, \overline{t}}(1) = 2v$ for all $v \in J \cap \SS$.
We can summarize this in \cref{ta:Lij}.
(We do not further specify the ``middle part'' $Z = L_0 \cap L'_0$.)
Recall our notation $s(b)$ for $b \in B$ from \cref{eq:s(b)}.

\begin{table}[ht!]
\begingroup
\renewcommand{\arraystretch}{1.3}
\[
\begin{array}{l|ccccc}
& L_{-2} & L_{-1} & L_{0} & L_{1} & L_{2} \\
\hline
L'_{2} & & & V_{J' \cap \SS, t} & B_+ & (J \cap \SS)_+ \\
L'_{1} & & B_- & V_{J, t} & J'_+ & (s(B))_+ \\
L'_{0} & (J \cap \SS)_- & J'_- & Z & J_+ & (J' \cap \SS)_+ \\
L'_{-1} & (s(B))_- & J_- & V_{J', \overline{t}} & B^\op_+ & \\
L'_{-2} & (J' \cap \SS)_- & B^\op_- & V_{J \cap \SS, \overline{t}} & &
\end{array}
\]
\endgroup
\smallskip
\caption{\small The components $L_{ij} = L_i \cap L'_j$}\label{ta:Lij}
\end{table}

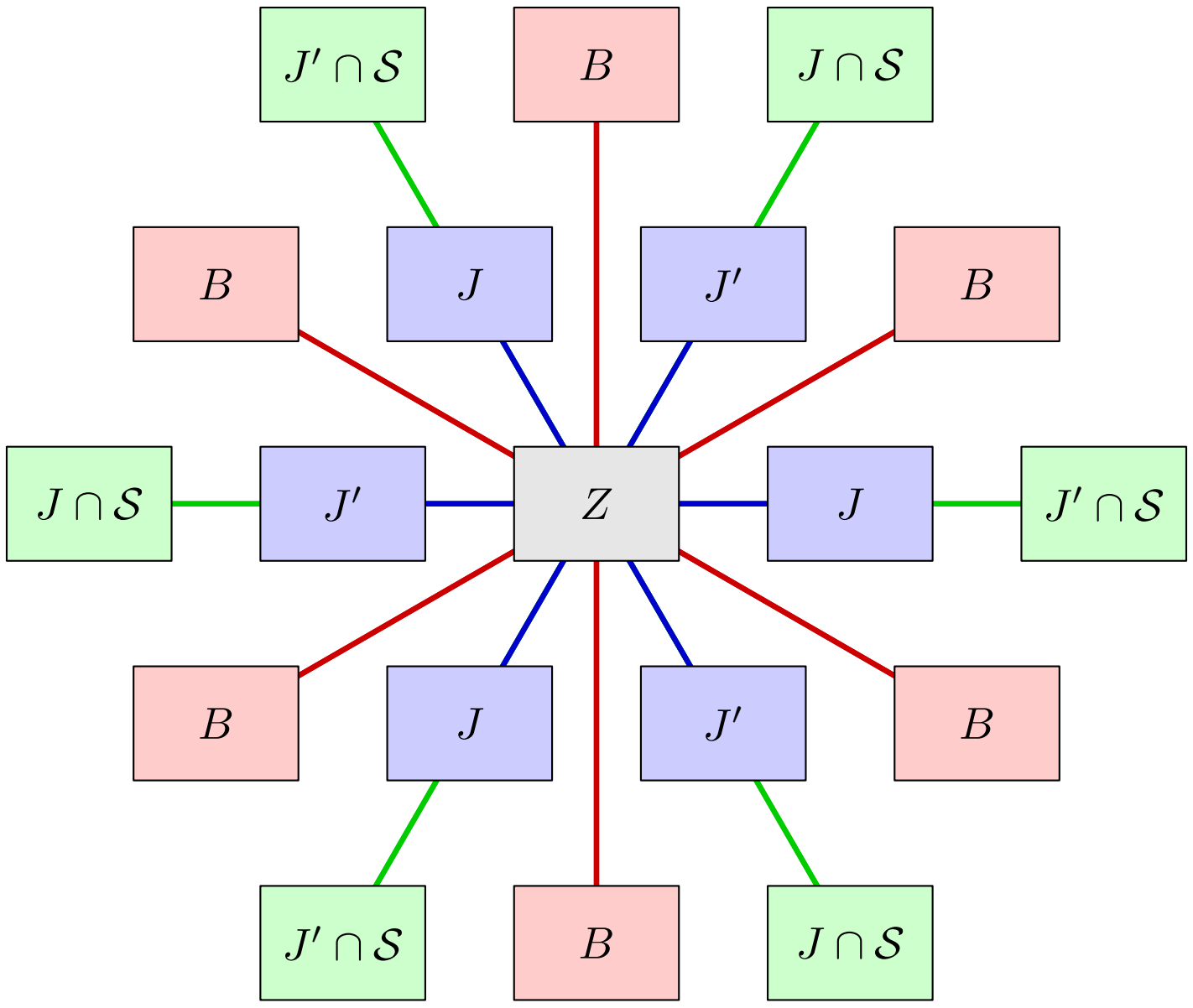
\begin{figure}
\begin{tikzpicture}[scale=2, line cap=round, line join=round]

\def\R{1}
\def\L{1.732050807568877}
\def\D{2}

\tikzset{
  short/.style={-Latex, very thick, black!20!blue},
  long/.style ={-Latex, very thick, black!20!red},
  dbl/.style  ={-Latex, very thick, black!20!green},
  box/.style ={draw, rectangle, minimum width=1.3cm, minimum height=0.9cm, inner sep=0pt, align=center},
  boxS/.style={box, fill=blue!20},
  boxL/.style={box, fill=red!20},
  boxD/.style={box, fill=green!20},
  boxC/.style={box, fill=gray!20},
}

\foreach \a in {0,60,120,180,240,300}{
  \draw[long]  (0,0) -- (\a+30:\L);
  \draw[long]  (0,0) -- (\a+210:\L);

  \draw[dbl]   (0,0) -- (\a:\D);
  \draw[dbl]   (0,0) -- (\a+180:\D);

  \draw[short] (0,0) -- (\a:\R);
  \draw[short] (0,0) -- (\a+180:\R);
}


\node[boxC] at (0,0) {$Z$};

\node[boxS] at (0:\R)   {$J$};
\node[boxS] at (60:\R)  {$J'$};
\node[boxS] at (120:\R) {$J$};
\node[boxS] at (180:\R) {$J'$};
\node[boxS] at (240:\R) {$J$};
\node[boxS] at (300:\R) {$J'$};

\node[boxL] at (30:\L)  {$B$};
\node[boxL] at (90:\L)  {$B$};
\node[boxL] at (150:\L) {$B$};
\node[boxL] at (210:\L) {$B$};
\node[boxL] at (270:\L) {$B$};
\node[boxL] at (330:\L) {$B$};

\node[boxD] at (0:\D)   {$J'\cap \SS$};
\node[boxD] at (60:\D)  {$J\cap \SS$};
\node[boxD] at (120:\D) {$J'\cap \SS$};
\node[boxD] at (180:\D) {$J\cap \SS$};
\node[boxD] at (240:\D) {$J'\cap \SS$};
\node[boxD] at (300:\D) {$J\cap \SS$};

\end{tikzpicture}
\caption{The $AG_2$-grading of the TKK Lie algebra}\label{fig:AG2}
\end{figure}

At this point, we cannot resist redrawing this parametrization in a slightly different way; see \cref{fig:AG2}.
In particular, we have shown that $L$ is $AG_2$-supported. In fact, combining everything, we get the following converse to \cref{th:AG2-A2}.

\begin{theorem}\label{th:AG2-A2-conv}
	Let $\A$ be a structurable algebra with conjugate-supplementary idempotent $t \in \A$ and let $L = \TKK$.
	Then $L$ is $(AG_2, A_2)$-graded.
\end{theorem}
\begin{proof}
	Let $\omega_1, \omega_2$ be the fundamental weights as in \cref{fig:weights}.
	Then all roots of $AG_2^o$ are of the form $j \omega_1 + (i-j) \omega_2$ with $i, j \in \{ -2, -1, 0, 1, 2 \}$.
	More precisely, we see that $j \omega_1 + (i-j) \omega_2 \in AG_2^o$ if and only if the item in column $L_i$ and row~$L'_j$ in \cref{ta:Lij} is non-empty.
	This shows that the first condition of \cref{def:supported} is satisfied.
	The second condition is also clear by \cref{def:TKK}, because $L_0 = \instrl(\A)$ is equal to $[L_1, L_{-1}]$, so certainly $Z = L_o$ is contained in $[L_1, L_{-1}]$.
	
	Next, by \cref{le:split A2}, $L$ contains a split simple Lie algebra $\fG$ of type $A_2$, with Cartan subalgebra $\fH = \langle \id, T_t \rangle \leq \instrl(\A)$.
	By \cref{def:TKK}, $\ad(\id)$ is precisely the grading derivation of $L$ with respect to the $L_i$, and by \cref{pr:5by5}, $\ad(T_t)$ is the grading derivation of $L$ with respect to the $L'_j$.
	Therefore, for each element $h = p \cdot \id + q \cdot T_t \in \fH$, $\ad(h)$ acts on each $L_{ij} = L_i \cap L'_j$ as multiplication by $pi + qj$.
	Since $\id$ and $T_t$ correspond to $h_1 + h_2$ and $h_1$ under the isomorphism in the proof of \cref{le:split A2}, we have $(j \omega_1 + (i-j) \omega_2)(h) = pi + qj$, and we conclude that $L$ is indeed $(AG_2, A_2)$-graded with grading algebra $\fG$ as in \cref{def:AG2-Delta-graded}.
\end{proof}

\begin{remark}\label{rem:feeble etc}
	If the $AG_2$-algebra $(B, J, J', \times, \times', T)$ is \emph{weak} (as introduced in \cref{def:AG2-algebra}), then the $AG_2$-grading from \cref{fig:AG2} is, in fact, a $G_2$-grading. If the $AG_2$-algebra is \emph{feeble}, then it is simply an $A_2$-grading.
	
	The $AG_2$-algebra is \emph{strong} if and only if all 19 pieces in the grading are non-trivial, i.e., if and only if $K(\A)$ is strongly supported by $AG_2$ as in \cref{def:supported}. A~strong $AG_2$-algebra is \emph{mighty} precisely when the Lie brackets between two distinct $J$-pieces and between two distinct $J'$-pieces are non-trivial. 
\end{remark}

\section{\texorpdfstring{$(AG_2, G_2)$}{(AG2, G2)}-graded Lie algebras}\label{se:AG2 G2}

We now investigate the effect of imposing a more restrictive grading algebra of type $G_2$ rather than $A_2$.
The following definition will be useful.
\begin{definition}\label{def:G2-frame}
	Let $\A$ be a structurable algebra and let $t \in \A$ be a conjugate-supplementary idempotent.
	Let $J = t \A \overline{t}$ and $J' = \overline{t} \A t$ as in \cref{se:idempotents2}.
	A pair of elements $u \in J$, $u' \in J'$ is called a \emph{$G_2$-frame} for $(\A, t)$ if
	\[ \overline{u} = u, \quad \overline{u'} = u', \quad u^2 = 2u', \quad u'^2 = 2u, \quad uu' = 3t, \quad u'u = 3\overline{t} . \]
	In particular, the $4$-dimensional subspace with basis $t, \overline{t}, u, u'$ is a subalgebra of~$\A$, which we denote by $M(u,u')$.
\end{definition}
First, we need to generalize \cref{le:V in kt} to this larger subalgebra $M(u,u')$.
\begin{lemma}\label{le:G2-frame}
	Let $\A$ be a structurable algebra, let $t \in \A$ be a conjugate-supplementary idempotent, and let $(u,u')$ be a $G_2$-frame for $(\A, t)$. Consider the subalgebra $\B := M(u,u')$. Then for all $x,y \in \B$, there exists some $z \in \B$ such that $V_{x,y} = T_z$ (as operators on $\A$).
\end{lemma}
\begin{proof}
	We claim that the following table expresses $V_{x,y}$ as $T_z$ for each choice of $x$ and $y$ in the basis $\{ t, \overline{t}, u, u' \}$.
	The result then follows, of course, by bilinearity.
	\[
	\begin{array}{c|cccc}
		V_{x,y}&y=t&y=\overline{t}&y=u&y=u'\\[.3ex]
		\hline
		x=t&0&T_t&T_u&0\\
		x=\overline{t}&T_{\overline{t}}&0&0&T_{u'}\\
		x=u&T_u&0&2T_{u'}&T_t+2T_{\overline{t}}\\
		x=u'&0&T_{u'}&2T_t+T_{\overline{t}}&2T_u
	\end{array}
	\]
	\Cref{le:V in kt} already shows our claim when both $x$ and $y$ belong to $\{ t, \overline{t} \}$.
	Moreover, \cref{le:rewrite V}, together with \cref{eq:Vxy+Vyx}, gives
	\begin{alignat*}{2}
		V_{u,t}=V_{t,u}&=T_u, \qquad
		&
		V_{u,\overline{t}}=V_{\overline{t},u}&=0,\\
		V_{u',\overline{t}}=V_{\overline{t},u'}&=T_{u'},
		&
		V_{u',t}=V_{t,u'}&=0.
	\end{alignat*}
	Since $\overline{u}=u$ and $\overline{u'}=u'$, \cref{le:rewrite V}\cref{it:V JJ,it:V J'J'} also give
	\[ V_{u,u} = V_{u^2,\overline{t}} = 2V_{u',\overline{t}} = 2T_{u'} \]
	and
	\[ V_{u',u'} = V_{u'^2,t} = 2V_{u,t} = 2T_u. \]
	Thus only $V_{u,u'}$ and $V_{u',u}$ remain.
	We compute the Lie bracket $[T_u, T_{u'}]$ in two different ways. First, we simply invoke the structurable identity (see \cref{def:struct algebra}) to get
	\begin{align}
		[T_u,T_{u'}]
		&=[V_{u,1},V_{u',1}] \notag \\
		&=V_{T_u(u'),1}-V_{u',V_{1,u}(1)} \notag \\
		&=3T_t-V_{u',u}. \label{eq:Tuv 1}
	\end{align}
	On the other hand, we use the structurable identity with $(x,y,z,w) = (u,t,u,u)$. Since we already know that $V_{u,t} = T_u$ and $V_{u,u} = 2 T_{u'}$, this becomes
	\[ [T_u, 2 T_{u'}] = V_{T_u(u), u} - V_{u, V_{t,u}(u)} . \]
	Now $T_u(u) = u^2 = 2u'$ and $V_{t,u}(u) = T_u(u) = 2u'$, so we get, after dividing by $2$,
	\[ [T_u, T_{u'}] = V_{u', u} - V_{u, u'} . \]
	Comparing this with \cref{eq:Tuv 1}, we obtain
	\[
	V_{u,u'}-2V_{u',u}=-3T_t.
	\]
	By~\cref{eq:Vxy+Vyx}, on the other hand, we have
	\[ V_{u,u'}+V_{u',u} = T_{u\overline{u'}+u'\overline{u}} = T_{uu'+u'u} = 3T_t+3T_{\overline{t}} . \]
	Combining the previous two equalities, we find that indeed
	\[ V_{u,u'}=T_t+2T_{\overline{t}} \qquad \text{and} \qquad V_{u',u}=2T_t+T_{\overline{t}}. \qedhere \]
\end{proof}

We can now prove a result similar to \cref{le:split A2}.
Recall the algebras $M(J, \eta)$ from \cref{ex:SMA}. Observe that for $J = k$, with the only possible non-degenerate cubic norm structure $N(a) = a^3$ and $T(a,b) = 3ab$, the algebra is $4$-dimensional.

\begin{lemma}\label{le:split G2}
	Let $\A$ be a structurable algebra, let $t \in \A$ be a conjugate-supplementary idempotent, and let $(u,u')$ be a $G_2$-frame for $(\A, t)$. Consider the subalgebra $\B := M(u,u')$. Then:
	\begin{enumerate}
		\item\label{it:B=Mk1} $\B \cong M(k, 1)$ as structurable algebras;
		\item $\B$ has TKK Lie algebra $K(\B) \cong \fG_2(k)$;
		\item $K(\B)$ embeds into $K(\A)$.
	\end{enumerate}
\end{lemma}
\begin{proof}
	\begin{enumerate}[leftmargin=1.5em]
		\item 
		This follows directly from the definitions in \cref{ex:SMA}. (Notice that $a \times b = 2ab$ for all $a,b \in k$.)
		\item
		It is well known that $M(k,1)$ has TKK Lie algebra $\fG_2(k)$; see, for example, \cite[p.\@~1871]{Allison1979} for an argument relying on the classification of the simple Lie algebras.
		Alternatively, set
		\[ e_\alpha := T_u, \qquad f_\alpha := T_{u'}, \qquad e_\beta := \overline{t}_+, \qquad f_\beta := t_- , \]
		and check that these elements generate $K(\B)$ and satisfy the Chevalley--Serre relations for the Cartan matrix
		$\begin{psmallmatrix} 2&-3\\ -1&2 \end{psmallmatrix}$.
		(Notice that $e_\beta$ and $f_\beta$ are the images of $e_2$ and $f_2$ under the isomorphism $\ssl_3(k) \to K(k[t])$ from the proof of \cref{le:split A2}; in particular, $\beta = \alpha_2$.)
		\item
		By \cref{le:G2-frame}, the exact same argument as in the proof of \cref{le:split A2} applies here.
		\qedhere
	\end{enumerate}
\end{proof}
Notice that the embedding of $K(\B)$ into $K(\A)$ is compatible with the distinguished long-root $A_2$, i.e., the restriction of this embedding to $K(k[t])$ coincides with the embedding of $K(k[t])$ into $K(\A)$ from \cref{le:split A2}.

\medskip

The following result might be known, but we could not find it explicitly in the literature.
(Recall that the TKK construction is not functorial!)
\begin{proposition}\label{pr:KB-KA}
	Let $\A$ and $\B$ be structurable algebras. Assume that
	\[ \Phi \colon K(\B) \to K(\A) \]
	is a Lie algebra homomorphism preserving the $5$-grading and satisfying $\Phi(1_+) = 1_+$ and $\Phi(1_-) = 1_-$.
	Then there is a homomorphism $\varphi \colon \B \to \A$ (as $k$-algebras with involution) inducing $\Phi$, i.e., such that $\Phi(x_\pm) = \varphi(x)_\pm$ for all $x \in \B$.
\end{proposition}
\begin{proof}
	The crucial ingredients are Faulkner's \emph{reconstruction formulas} from \cite[Lemma 1.7]{Faulkner1994} (see also \cite[Theorem 2.9]{Elduque14} and its proof):
	If we write $\{x,y,z\} := V_{x,y}(z)$, then
	\begin{align}
		\overline{x} &= 2x - \{ x, 1, 1 \}, \quad \text{and} \label{eq:reconstruction1} \\
		xy &= \tfrac{1}{2} \bigl(	\{ x, 1, y \} + \{ 1, x, y \} + \{ x, y, 1 \} - \{ 1, y, x \} \bigr). \label{eq:reconstruction2}
	\end{align}
	Since $\Phi$ preserves the grading, there exist $k$-linear maps $\varphi_\pm \colon \B \to \A$ such that $\Phi(x_\pm) = \varphi_\pm(x)_\pm$ for all $x \in \B$. We know that $\varphi_\pm(1) = 1$.
	By \cref{def:TKK}, $\{ x,y,z \}_+ = [[x_+, y_-], z_+]$ for all $x,y,z \in \B$.
	Since $\Phi$ is a Lie algebra homomorphism,
	$\Phi(\{ x,y,z \}_+) = [[\Phi(x_+), \Phi(y_-)], \Phi(z_+)]$, hence
	\begin{equation}\label{eq:varphi+}
		\varphi_+( \{ x,y,z \}) = \{ \varphi_+(x), \varphi_-(y), \varphi_+(z) \}
	\end{equation}
	for all $x,y,z \in \B$.
	By \cref{eq:reconstruction1}, this already implies that $\varphi_+(\overline{x}) = \overline{\varphi_+(x)}$ for all $x \in \B$.
	Similarly, $\varphi_-(\overline{x}) = \overline{\varphi_-(x)}$ for all $x \in \B$.
	
	We now apply \cref{eq:varphi+} with $x=z=1$ and $y$ arbitrary; note that $\{ 1,y,1 \} = 2\overline{y} - y$.
	Using the fact that $\varphi_-(\overline{y}) = \overline{\varphi_-(y)}$ and that $\varphi_+(1) = 1$,
	we get
	\[ \varphi_+(2\overline{y} - y) =  \varphi_+( \{ 1,y,1 \}) = \{ 1, \varphi_-(y), 1 \} = \varphi_-(2 \overline{y} - y) . \]
	Because the map $y \mapsto 2 \overline{y} - y$ is invertible (with inverse $x \mapsto \tfrac{1}{3}(2\overline{x} + x)$), this shows that $\varphi_+ = \varphi_-$. Writing $\varphi := \varphi_\pm$ now, it follows from \cref{eq:reconstruction2,eq:varphi+}, together with $\varphi(1)=1$, that $\varphi$ preserves the multiplication.
\end{proof}

We now have the following extension of \cref{th:AG2-A2,th:AG2-A2-conv} to the setting of $(AG_2, G_2)$-graded Lie algebras.
\begin{theorem}\label{th:G2-frame}
	\begin{enumerate}
		\item 
			Let $L$ be an $(AG_2,G_2)$-graded Lie algebra with grading algebra $\fG$.
			Then there exist a structurable algebra $\A$, a conjugate-supplementary idempotent $t\in \A$, a $G_2$-frame $(u,u')$ for $(\A,t)$, and an isomorphism
			\[ L/Z(L) \cong \TKK . \]
			Moreover, the isomorphism can be chosen so that the image of $\fG$ is precisely
			\[ K(M(u,u')) \cong \fG_2(k) . \]
		\item 
			Conversely, let $\A$ be a structurable algebra, let $t\in \A$ be a conjugate\dash supplementary idempotent, and assume that $(\A,t)$ admits a $G_2$-frame $(u,u')$.
			Then $\TKK$ is $(AG_2,G_2)$-graded, with grading algebra
			\[ K(M(u,u')) \cong \fG_2(k) . \]
	\end{enumerate}
\end{theorem}
\begin{proof}
	\begin{enumerate}[leftmargin=1.5em]
		\item 
			Let $L$ be $(AG_2,G_2)$-graded, with grading algebra $\fG \cong \fG_2(k)$.
			Let $\fG_\ell$ be the subalgebra of $\fG$ generated by the long-root spaces; this is a split simple Lie algebra of type $A_2$, and $L$ is $(AG_2,A_2)$-graded with grading algebra $\fG_\ell$.
			
			Now let $\B = M(k,1)$, equipped with its standard conjugate-supplementary idempotent $t_0$ and its standard $G_2$-frame $(u_0,u'_0)$.
			By \cref{le:split G2} and \cref{ta:Lij}, $K(\B)\cong\fG_2(k)$, with split Cartan subalgebra $\fH_0 := \langle \id, T_{t_0} \rangle$.
			By the isomorphism theorem for split simple Lie algebras, we can choose an isomorphism
			\[ \theta\colon K(\B)\to \fG \quad \text{with} \quad \theta(\fH_0) = \fH . \]
			Then $\theta$ maps long-root spaces to long-root spaces, and hence maps $K(k[t_0])$ onto $\fG_\ell$.
			We choose the Chevalley generators of $\fG_\ell$ (see \cref{def:efh}) to be
			\[ e_1 := \theta\bigl((t_0)_+\bigr), \quad e_2 := \theta\bigl((\ov{t_0})_+\bigr), \quad f_1 := \theta\bigl((\ov{t_0})_-\bigr), \quad f_2 := \theta\bigl((t_0)_-\bigr) ; \]
			this is a valid choice, because $\theta(\fH_0) = \fH$ and because $(t_0)_+$, $(\ov{t_0})_+$, $(\ov{t_0})_-$ and $(t_0)_-$ are the images of Chevalley generators of $\ssl_3(k)$ under the isomorphism $\ssl_3(k) \to K(k[t_0])$ from the proof of \cref{le:split A2}.
			With this choice, \cref{th:AG2-A2} gives a structurable algebra $\A$, a conjugate-supplementary idempotent $t\in \A$, and a homomorphism
			\[ \pi\colon L\to \TKK \]
			with kernel $Z(L)$, such that $\pi(e_1) = t_+$, $\pi(e)=1_+$ and $\pi(f)=1_-$.
			Since $\theta(1_+)=e$, $\theta(1_-)=f$ and $\theta(\id) = h$, the map
			\[ \Phi := \pi\circ\theta\colon K(\B)\to K(\A) \]
			is a $5$-graded Lie algebra homomorphism satisfying $\Phi(1_\pm) = 1_\pm$.
			Note that $\Phi$~is injective because $\fG \cap Z(L) = 0$.
			
			We can now invoke \cref{pr:KB-KA} to get a homomorphism of structurable algebras
			\[ \varphi\colon \B\to \A \]
			such that $\Phi(x_\pm)=\varphi(x)_\pm$ for all $x\in \B$. Since $\Phi$ is injective, so is $\varphi$. Since $\Phi((t_0)_+) = \pi(e_1) = t_+$, we have $\varphi(t_0)=t$. Set
			\[ u := \varphi(u_0) \qquad \text{and} \qquad u' := \varphi(u'_0). \]
			Then $(u,u')$ is a $G_2$-frame for $(\A,t)$, and $\varphi(\B) = M(u,u')$. Moreover,
			\[ \pi(\fG) = \Phi\bigl(K(\B)\bigr) = K(M(u,u')) \cong \fG_2(k). \]
		\item 
			Conversely, let $\A$ be a structurable algebra with a conjugate-supplementary idempotent $t$ and a $G_2$-frame $(u,u')$.
			Put $\B := M(u,u')$. By \cref{le:split G2}, the Lie algebra
			\[ \fG := K(\B) \cong \fG_2(k) \]
			embeds into $\TKK$. By the compatibility statement following \cref{le:split G2}, its long-root $A_2$-subalgebra is the copy of $K(k[t])$ used in \cref{th:AG2-A2-conv}, and its Cartan subalgebra is
			\[ \fH = \langle \id_\A, T_t \rangle . \]
			The $AG_2$-grading constructed in \cref{th:AG2-A2-conv} is precisely the weight-space decomposition of $L := \TKK$ for the action of this Cartan subalgebra.
			Since $\fH \subseteq \fG$, it follows that each root space $\fG_\gamma$ of $\fG$ (with respect to $\fH$) is contained in the corresponding weight space $L_\gamma$ of $L$.
			It remains to check that the roots of $\fG$ are precisely the roots of the ``solid'' root subsystem $G_2$ of $AG_2$ (see \cref{fig:AG2 GRS}).
			Now $\B \cap \SS = k(t - \overline{t})$ and $\instrl(\B) = T_\B = \langle T_t, T_{\overline{t}}, T_u, T_{u'} \rangle$ by \cref{le:G2-frame}, so
			\[ \fG = (\B \cap \SS)_- \oplus \B_- \oplus T_\B \oplus \B_+ \oplus (\B \cap \SS)_+ . \]
			By \cref{ta:Lij}, the elements $t_\pm$, $\overline{t}_\pm$ and $(t - \overline{t})_\pm$ lie in the six long-root spaces (one in each), the elements $u_\pm$, $u'_\pm$, $T_u = V_{u,t}$ and $T_{u'} = V_{u',\overline{t}}$ lie in the six short-root spaces (one in each), and $T_t, T_{\overline{t}} \in L_o$.
			Hence the roots of $\fG$ are precisely the long and short roots of $AG_2$.
			Thus both conditions of \cref{def:AG2-Delta-graded} are satisfied, and $\TKK$ is $(AG_2,G_2)$-graded with grading algebra $K(M(u,u')) \cong \fG_2(k)$.
		\qedhere
	\end{enumerate}
\end{proof}

\section{Conjugate-supplementary idempotents in central simple structurable algebras}\label{se:examples}

Our next goal is to classify all finite-dimensional central simple structurable algebras admitting a conjugate-supplementary idempotent.
We recall from \cite[\S 2.3]{Boelaert2019} that, when $\Char(k) \neq 5$, every finite-dimensional central simple structurable algebra belongs to one of the following six (non-disjoint) classes; this is a result due to Allison and Smirnov \cite{Allison1978,Smirnov1990,Smirnov1992}.
(Remember that we have assumed throughout the paper that $\Char(k) \neq 2,3$.)

\begin{enumerate}[(1)]
	\item\label{cssa:1} central simple associative algebras-with-involution%
	\footnote{The point of our somewhat unusual addition of the hyphens in ``algebra-with-involution'' is the fact that we are emphasizing that ``central simple'' refers to the ``algebra with involution'', which is not the same as a ``central simple algebra'' equipped with an involution. See \cref{le:csa assoc}.},
	\item\label{cssa:2} central simple Jordan algebras,
	\item\label{cssa:3} structurable algebras constructed from a non-degenerate hermitian form over a central simple associative algebra with involution,
	\item\label{cssa:4} central simple structurable algebras of skew-dimension 1; these are forms of the ``structurable matrix algebras'' from~\cref{ex:SMA},
	\item\label{cssa:5} forms of the tensor product of two composition algebras,
	\item\label{cssa:6} ``Smirnov algebras'', exceptional 35-dimensional algebras arising as subalgebras of the tensor product of two copies of the same octonion algebra.
\end{enumerate}
We will rely on the explicit description of these examples that can be found in detail in~\cite[\S 2.3]{Boelaert2019}.

The algebras in these six classes also make sense when $\Char(k) = 5$, and our analysis below applies to them without change; what is not known in characteristic~$5$ is whether the list is complete.
In this direction, Stavrova has shown that when $\Char(k) \neq 2,3$, every finite-dimensional central simple Lie algebra with a non-trivial $5$-grading is of Chevalley type \cite[Theorem~3.1]{Stavrova22}; see also \cite[Theorem~1.1 and \S 1]{Stavrova22}.
In particular, if $L$ is an $(AG_2,A_2)$-graded Lie algebra such that $L/Z(L)$ is finite-dimensional and central simple, then $L/Z(L) \cong \TKK$ (see \cref{th:AG2-A2}) is of Chevalley type.

\subsection{Central simple associative algebras-with-involution}\label{ss:ex1}

In the associative case, conjugate-supplementary idempotents have been studied before.
It turns out that they are directly related to the notion of \emph{hyperbolic involutions}.
We follow the approach from \cite{BayerFluckiger1993}; see also \cite[Chapter~6]{Knus1998}.
Notice that, as explained in \cite[pp.\@~20--21]{Knus1998}, we have to translate ``central simple as an algebra-with-involution'' (which coincides with central simple as a structurable algebra) to the usual context of central simple algebras.
The argument of the following statement can also be found in \cite[p.\@~208]{Jacobson1968}, but we present it for completeness.
\begin{lemma}\label{le:csa assoc}
	Let $A$ be an associative $k$-algebra with involution.
	Then $A$ is central simple as a structurable $k$-algebra if and only if
	\begin{enumerate}
		\item\label{it:csa assoc:1} either $A$ is a central simple $k$-algebra, equipped with an involution of the first kind;
		\item\label{it:csa assoc:2} or $A$ is a central simple $K$-algebra for some quadratic extension~$K/k$, equipped with an involution of the second kind fixing $k$ but not $K$;
		\item\label{it:csa assoc:3} or $A \cong B \oplus B^\op$ for some central simple $k$-algebra $B$, equipped with the exchange involution.
	\end{enumerate}
\end{lemma}
\begin{proof}
	Observe that the ideals of $A$ as a structurable algebra are precisely the involution-invariant two-sided ideals of~$A$, and that the center of $A$ as a structurable algebra consists precisely of the elements of $Z(A)$ fixed by the involution.
	In particular, the ``if'' direction is clear: in cases \ref{it:csa assoc:1} and \ref{it:csa assoc:2}, $A$ has no proper non-zero two-sided ideals at all, and in case \ref{it:csa assoc:3}, the only involution-invariant two-sided ideals of $B \oplus B^\op$ are $0$ and $A$; moreover, in all three cases, the elements of $Z(A)$ fixed by the involution are precisely the elements of~$k$.
	
	We now prove the ``only if'' direction.
	Suppose first that $A$ is simple as a $k$\dash algebra.
	Then $Z(A)$ is a field, and the involution of $A$ acts on $Z(A)$. Since $A$ is central as a structurable algebra, the subfield of $Z(A)$ fixed by the involution is precisely~$k$.
	Hence either $Z(A) = k$ and we are in case \ref{it:csa assoc:1}, or $Z(A) = K$ for some separable quadratic extension $K/k$ and we are in case \ref{it:csa assoc:2}.
	
	Suppose next that $A$ is not simple as a $k$-algebra, so there exists a proper non-zero ideal $B \unlhd A$. Since $A$ is simple as a structurable algebra, the ideal $B + \overline{B}$ coincides with the $k$-algebra~$A$. Since $B \cap \overline{B}$ is now an ideal invariant under the involution, we must have $B \cap \overline{B} = 0$, so we already have $A = B \oplus \overline{B} \cong B \oplus B^\op$, and the involution of $A$ is the exchange involution.
	It is also clear that $B$ is simple, because a proper non-zero ideal $I \unlhd B$ would give rise to a proper non-zero involution-invariant ideal $I \oplus \overline{I} \unlhd A$.
	Finally, $B$ is a central $k$-algebra, because the hermitian central elements of $A$ are precisely the elements $\{ b + \overline{b} \mid b \in Z(B) \}$, which, by assumption, coincides with $k$.
	We conclude that we are in case \ref{it:csa assoc:3}.
\end{proof}
Observe now that \cite[Definition (6.8)]{Knus1998} says precisely that a central simple algebra-with-involution $(A, \sigma)$ is called \emph{hyperbolic} if and only if $A$ contains an idempotent $e$ such that $\sigma(e) = 1-e$, which is exactly our definition of a conjugate-supplementary idempotent. If $(A, \sigma)$ is as in cases \ref{it:csa assoc:1} or \ref{it:csa assoc:2} of \cref{le:csa assoc}, then \cite[Proposition (6.7)]{Knus1998} gives equivalent formulations for being hyperbolic.
See also \cite[Theorem 2.1]{BayerFluckiger1993}.
Using \cite[Theorem 2.2]{BayerFluckiger1993}, we push this result a bit further to get precise information about the structure of the resulting algebra with involution.
Notice that for \emph{orthogonal} involutions, the following proposition is essentially contained in \cite[p.\@~145, Exercise 2]{Knus1998}.
\begin{proposition}\label{pr:assoc idemp}
	Let $A$ be an associative $k$-algebra with involution, central simple as a structurable $k$-algebra.
	Assume that $t \in A$ is a conjugate-supplementary idempotent. Then
	\begin{enumerate}
		\item\label{it:assoc idemp:1} either $A \cong M_2(B)$ for some central simple $K$-algebra $B$ invariant under the involution of $A$, where either $K=k$ or $K/k$ is a quadratic field extension, and the isomorphism can be chosen so that $t$ is mapped to $\idempotentmatrix$ and the involution on $M_2(B)$ is given by
		\[ \overline{\begin{pmatrix} a & b \\ c & d \end{pmatrix}} = \begin{pmatrix} \overline{d} & \epsilon \overline{b} \\ \epsilon\overline{c} & \overline{a} \end{pmatrix} , \]
		where $b \mapsto \overline{b}$ is the restriction to $B$ of the involution of $A$.
		If the involution of $A$ is of the first kind, then the involution of $B$ is orthogonal, and $\epsilon = 1$ or $\epsilon = -1$ according to whether the involution of $A$ is orthogonal or symplectic.
		If the involution of $A$ is of the second kind, then so is the involution of $B$, and $\epsilon = 1$;
		\item\label{it:assoc idemp:2} or $A \cong B \oplus B^\op$ for some central simple $k$-algebra $B$ and the isomorphism maps $t$ to an element $(u, \overline{1-u})$ for some (arbitrary) idempotent $u \in B$.	
	\end{enumerate}
\end{proposition}
\begin{proof}
	Assume first that we are in case \ref{it:csa assoc:1} or \ref{it:csa assoc:2} of \cref{le:csa assoc}.
	Then we can apply \cite[Theorem 2.2]{BayerFluckiger1993} to find an involution-invariant subalgebra $M$ of $A$ with $M \cong M_2(K)$, where $K=k$ or $K/k$ is a quadratic field extension, depending on whether we are in case \ref{it:csa assoc:1} or \ref{it:csa assoc:2} of \cref{le:csa assoc};
	moreover, the restriction of the involution on $A$ to the subalgebra $M = M_2(K)$ is given by
	\begin{equation}\label{eq:M2L}
		\overline{\begin{pmatrix}
			a & b \\ c & d
		\end{pmatrix}} = \begin{pmatrix}
			\overline{d} & \epsilon \overline{b} \\ \epsilon\overline{c} & \overline{a}
		\end{pmatrix}
	\end{equation}
	for all $a,b,c,d \in K$,
	where $\epsilon = 1$ if the involution is orthogonal, $\epsilon = -1$ if the involution is symplectic, and $\epsilon = 1$ if the involution is of the second kind.
	In addition, the proof of \cite[Theorem 2.2]{BayerFluckiger1993} shows%
	\footnote{In one direction, this is stated explicitly in the proof, but in the other direction, which we need, this requires going through some details, using the proof of implication (d)$\Rightarrow$(a) of \cite[Theorem 2.1]{BayerFluckiger1993}.}
	that the idempotent $t$ corresponds to $\idempotentmatrix \in M_2(K)$.
	
	Now let $B$ be the centralizer of $M$ in $A$ and note that $M$ is a central simple $K$-subalgebra of $A$.
	By the double centralizer theorem (see, e.g., \cite[(1.5)]{Knus1998}), $B$ is again a central simple $K$-subalgebra of $A$ and moreover,
	the multiplication in $A$ defines a canonical isomorphism $A = M \otimes_K B$.
	Since $M$ is involution-invariant, so is its centralizer $B$.
	If the involution is of the first kind, then the involution on $B$ is orthogonal, because the involution \eqref{eq:M2L} on $M \cong M_2(k)$ has the same type as the involution of $A$, and the type of a tensor product of involutions of the first kind is the product of the types \cite[(2.23)]{Knus1998}.
	If the involution is of the second kind, then it is non-trivial on the center $K \subseteq B$, so its restriction to $B$ is again of the second kind.
	This implies that the involution on $A = M \otimes_K B \cong M_2(B)$ is indeed given by the same formula~\eqref{eq:M2L}, where now $a,b,c,d \in B$.
		
	Assume now that we are in case \ref{it:csa assoc:3} of \cref{le:csa assoc}.
	Write $t = (u, \overline{v}) \in B \oplus B^\op$ with $u,v \in B$.
	Then $\overline{t} = (v, \overline{u})$, so $t + \overline{t} = 1$ immediately implies that $u+v = 1$.
	Conversely, if $u \in B$ is an arbitrary idempotent, then $t := (u, \overline{1-u})$ is a conjugate-supplementary idempotent in $A$.
\end{proof}

\subsection{Central simple Jordan algebras}\label{ss:ex2}

Recall from \cite[p.~135, Remark~(ii)]{Allison1978} or \cite[\S 2.3.2]{Boelaert2019} that $\A$ is a Jordan algebra if and only if $\SS = 0$. By \cref{le:csi to skew}, such an algebra cannot admit a conjugate-supplementary idempotent.

\subsection{Structurable algebras from a non-degenerate hermitian form}\label{ss:ex3}

We recall the construction from \cite[\S 7]{Allison1979}; see also \cite[Definition 2.3.1]{Boelaert2019}.
\begin{definition}\label{def:stralg hermitian}
	Let $E$ be an associative $k$-algebra with $k$-linear involution $x \mapsto \overline{x}$.
	Let $M \neq 0$ be a right%
	\footnote{\cite{Allison1978} and \cite{Boelaert2019} use left $E$-modules and \cite{Allison1979} uses right $E$-modules, but as pointed out in \cite{Allison1979}, this choice is irrelevant because $E$ has an involution.}
	$E$-module, and let $h \colon M \times M \to E$ be a hermitian form, i.e., $h$~is $k$-bilinear and satisfies $\overline{h(m,n)} = h(n,m)$ and $h(mr, ns) = \overline{r} h(m,n) s$ for all $m,n \in M$ and all $r,s \in E$.
	Let $\A := E \oplus M$, equipped with multiplication
	\[ (r, m) \cdot (s, n) := \bigl( rs + h(n,m), \  nr + m\overline{s} \bigr) \]
	and involution
	\[ \overline{(r,m)} := (\overline{r}, m) . \]
	Then $\A$ is a structurable algebra.
	Moreover, $\A$ is central simple if and only if $E$ is a central simple algebra-with-involution, $\dim(E) \dim(M) \neq 1$ and $h$ is non-degenerate; see \cite[Proposition 4.1]{Allison1986}.
\end{definition}

Determining the conjugate-supplementary idempotents in such a structurable algebra immediately reduces to the case of associative algebras with involution.
\begin{proposition}\label{pr:csi herm}
	Let $\A = E \oplus M$ be as in \cref{def:stralg hermitian}.
	Then $t \in \A$ is a conjugate-supplementary idempotent if and only if $t = (u, 0)$ for some conjugate-supplementary idempotent $u \in E$.
\end{proposition}
\begin{proof}
	Observe that all skew elements of $\A$ are contained in $E \oplus 0$, because the involution fixes $M$ elementwise.
	Since $\A$ has identity element $1 = (1, 0)$, the result now follows from \cref{le:csi to skew}.
\end{proof}
The following result will be relevant when we look at the corresponding Peirce decomposition; see \cref{ex:herm:B+B,ex:herm:M2B}.

\begin{proposition}\label{pr:Mt=Mbart}
	Let $E$ be an associative $k$-algebra with involution, central simple as a structurable $k$-algebra, and let $t \in E$ be a conjugate-supplementary idempotent.
	Let $M$ be a right $E$-module, and consider the $k$-subspaces $Mt$ and $M\overline{t}$.
	\begin{enumerate}
		\item If $E$ is as in case~\ref{it:assoc idemp:1} of \cref{pr:assoc idemp}, then $Mt \cong M\overline{t}$.
		\item If $E$ is as in case~\ref{it:assoc idemp:2} of \cref{pr:assoc idemp} and there exists a non-degenerate hermitian form $h \colon M \times M \to E$, then $Mt \cong M\overline{t}$.
	\end{enumerate}
\end{proposition}
\begin{proof}
	Observe that $M = Mt \oplus M\overline{t}$ as vector spaces. An element $m \in M$ belongs to~$Mt$ if and only if $mt = m$, and similarly, $m$ belongs to $M\overline{t}$ if and only if $m\overline{t} = m$.
	
	Assume first that $E = M_2(B)$ as in case~\ref{it:assoc idemp:1} of \cref{pr:assoc idemp}, with $t = e_{11}$ and $\overline{t} = e_{22}$, where we write $e_{ij}$ for the matrix with $1$ on position $(i,j)$ as usual. Let
	\begin{align*}
		\alpha &\colon Mt \to M\overline{t} \colon m \mapsto me_{12} , \\
		\beta &\colon M\overline{t} \to Mt \colon m \mapsto me_{21} .
	\end{align*}
	Notice that indeed $me_{12} \in M\overline{t}$ because $me_{12} \overline{t} = me_{12}e_{22} = me_{12}$ for all $m \in Mt$; similarly, $me_{21} \in Mt$ for all $m \in M\overline{t}$.
	Moreover, $(\beta \circ \alpha)(m) = me_{11} = mt = m$ for all $m \in Mt$, and similarly $(\alpha \circ \beta)(m) = m$ for all $m \in M\overline{t}$. We conclude that $Mt \cong M\overline{t}$.
	
	Assume next that $E = B \oplus B^\op$ as in case~\ref{it:assoc idemp:2} of \cref{pr:assoc idemp}, with $t = (u, \overline{1-u})$ for some idempotent $u \in B$.
	On the one hand, we have our decomposition $M = Mt \oplus M\overline{t}$, but on the other hand, we can also decompose $M$ with respect to the \emph{central} idempotents $e := (1,0)$ and $f := \overline{e} = (0,1)$ as $M = Me \oplus Mf$.
	Since $e$ and $f$ are central, the parts $Mt$ and $M\overline{t}$ also decompose accordingly, so we get
	\[ M = Mte \oplus Mtf \oplus M\overline{t}e \oplus M\overline{t}f . \]
	
	Now let $h \colon M \times M \to E$ be a non-degenerate hermitian form, so $h(mb, nc) = \overline{b} h(m,n) c$ for all $b,c \in E$.
	In particular, $h(me, ne) = f h(m,n) e = 0$ and $h(mf, nf) = e h(m,n) f = 0$ for all $m,n \in M$.
	Consequently, the restriction of $h$ to $Mf \times Me$ is still non-degenerate.
	Also observe that $h(mf, ne) = e h(m,n) e = h(m,n) e \in B \oplus 0$ for all $m,n \in M$.
	
	Since $B$ is a central simple $k$-algebra, it comes equipped with a \emph{reduced trace map} $\Trd \colon B \to k$. Consider the $k$-bilinear form
	\[ \varphi \colon Mf \times Me \to k \colon (mf, ne) \mapsto \Trd(h(m,n)e) , \]
	which is again non-degenerate.
	Now observe that
	\[ \varphi(mtf, nte) = \Trd(h(mt,nt) e) = \Trd(\overline{t} h(m,n) t e) = \Trd(h(m,n) te\overline{t}) = 0 \]
	for all $m,n \in M$ because $e$ is central and $t \overline{t} = 0$; similarly, $\varphi(m\overline{t}f, n\overline{t}e) = 0$ for all $m,n \in M$.
	Since $Mf = Mtf \oplus M\overline{t}f$ and $Me = Mte \oplus M\overline{t}e$, the restriction of $\varphi$ to $M\overline{t}f \times Mte$ is still non-degenerate. It follows that $M\overline{t}f \cong Mte$ as vector spaces; similarly, $Mtf \cong M\overline{t}e$.
	We conclude that
	\[ Mt = Mte \oplus Mtf \cong M\overline{t}f \oplus M\overline{t}e = M\overline{t} . \qedhere \]
\end{proof}

\subsection{Structurable algebras of skew-dimension one}\label{ss:ex4}

Structurable algebras of skew-dimension one have been studied in detail in \cite{Allison1990}.

\begin{proposition}\label{pr:skew-dim one}
	Let $\A$ be a central simple structurable algebra with $\dim \SS = 1$.
	Then $\A$ admits a conjugate-supplementary idempotent $t$ if and only if $\A \cong M(J, \eta)$ as in \cref{ex:SMA}.
	Moreover, in this case, we can choose the isomorphism in such a way that $t$ is mapped to $\idempotentmatrix$.
\end{proposition}
\begin{proof}
	Write $\SS = k s_0$ for some fixed $s_0 \in \SS$.
	By \cite[Lemma 2.1(b)]{Allison1984}, the element $\mu := s_0^2$ lies in $k \setminus \{ 0 \}$.
	Moreover, by \cite[Theorem 1.13]{Allison1990}, $\mu$ is a square if and only if $\A$ is isomorphic to a structurable matrix algebra as defined in \cref{ex:SMA}.
	Now observe that by \cref{le:csi to skew}, $t \in \A$ is a conjugate-supplementary idempotent if and only if $s := 2t - 1$ lies in $\SS$ and satisfies $s^2 = 1$. In particular, such a $t$ exists if and only if there exists an $s \in \SS$ with $s^2 = 1$, if and only if $\mu$ is a square.
	
	Finally, notice that the proof of \cite[Theorem 1.13]{Allison1990} relies on \cite[Theorem~4.5]{Allison1986}, which, in turn, relies on the arguments of \cite[pp.\@~153--154]{Allison1978}. In particular, the isomorphism obtained explicitly in the paragraph directly following equation (65) on p.\@~154 of \emph{loc.\@~cit.\@} maps $e_1 := \tfrac{1}{2}(1 + s_0)$ (with $s_0^2 = 1$) to $\idempotentmatrix$.
\end{proof}

\subsection{Tensor products of composition algebras}\label{ss:ex5}

\begin{proposition}\label{pr:csi C1xC2}
	\begin{enumerate}
		\item\label{pr:csi C1xC2:1}
			Let $\A \cong C_1 \otimes C_2$ be a tensor product of two composition algebras.
			Then $t \in \A$ is a conjugate-supplementary idempotent if and only if, up to interchanging $C_1$ and $C_2$,
			$C_1$ is split of dimension $\geq 2$ and $t = t_1 \otimes 1$ for some conjugate-supplementary idempotent $t_1 \in C_1$.
		\item 
			Let $\A$ be a form of a tensor product of two composition algebras.
			If $\A$ contains a conjugate-supplementary idempotent, then $\A$ is already itself a tensor product of two composition algebras.
	\end{enumerate}
\end{proposition}
\begin{proof}
	\begin{enumerate}
		\item 
			By \cref{le:csi to skew}, $t \in \A$ is a conjugate-supplementary idempotent if and only if $s := 2t - 1$ lies in $\SS$ and satisfies $s^2 = 1$.
			Write $C_1 = \SS_1 + \HH_1$ and $C_2 = \SS_2 + \HH_2$ with $\SS_i$ and $\HH_i$ the skew elements and hermitian elements of $C_i$, respectively. Since each $C_i$ is a composition algebra equipped with the standard involution, we have $\HH_1 = \HH_2 = k$. In particular, $\SS = \SS_1 \otimes 1 \oplus 1 \otimes \SS_2$, so we can write $s = s_1 \otimes 1 + 1 \otimes s_2$ with $s_1 \in \SS_1$ and $s_2 \in \SS_2$. Since $s^2 = 1$, we have
			\[ 1 \otimes 1 = \underbrace{s_1^2 \otimes 1 + 1 \otimes s_2^2}_{\in \HH_1 \otimes \HH_2} + \underbrace{2 s_1 \otimes s_2}_{\in \SS_1 \otimes \SS_2} \, , \]
			hence $s_1 \otimes s_2 = 0$. Up to interchanging $C_1$ and $C_2$, we may assume that $s_2 = 0$ and thus $s = s_1 \otimes 1$, with $s_1^2 = 1$.
			We conclude that $t = t_1 \otimes 1$, where $t_1 := \tfrac{1}{2}(s_1 + 1)$, so $t_1 \in C_1$ is a conjugate-supplementary idempotent.
			In particular, $C_1$ must be split of dimension $\geq 2$.
		\item 
			If $\A$ is a form of a tensor product of two composition algebras that is not itself such a tensor product,
			then by \cite[Theorem~2.1]{Allison1988} and \cite[Proposition~7.9]{Allison1992} (see also \cite[Proposition~2.3.14]{Boelaert2019}), $\A$ is a so-called \emph{twisted $(m,m)$-product algebra} for some $m \in \{2,4,8\}$.
			This means that there are a quadratic field extension $K/k$ with non-trivial automorphism~$\sigma$ and a composition algebra $C$ over $K$ of dimension~$m$ such that $\A$ is the $k$-subalgebra of $C \otimes_K {}^\sigma C$ of elements fixed by the $\sigma$-semilinear automorphism $\tau \colon x \otimes y \mapsto y \otimes x$; in particular, $\A \otimes_k K \cong C \otimes_K {}^\sigma C$.
			Suppose that $t \in \A$ were a conjugate-supplementary idempotent.
			Then $t$ is also a conjugate-supplementary idempotent of $C \otimes_K {}^\sigma C$, so by \cref{pr:csi C1xC2:1} (applied over~$K$), we have $t = t_1 \otimes 1$ or $t = 1 \otimes t_2$.
			In the first case, $\tau(t) = 1 \otimes t_1 \neq t$, because $t_1 \not\in K.1$ (indeed, $\ov{t_1} = 1 - t_1 \neq t_1$); the second case is similar.
			This contradicts $\tau(t) = t$.
		\qedhere
	\end{enumerate}
\end{proof}
\begin{remark}\label{rem:comp-idemp}
	If $C$ is a composition algebra (with standard involution), then any non-trivial idempotent $t \in C$ is automatically conjugate-supplementary.
	Indeed, any element $x \in C$ satisfies the quadratic equation
	\[ x^2 - T(x)x + N(x) = 0 , \]
	where $T(x) = x + \overline{x} \in k$ and $N(x) = x \overline{x} \in k$. In particular, since $t^2 = t$, we have $(1 - T(t))t + N(t) = 0$,
	and since $t \neq 0,1$, the elements $t$ and $1$ are linearly independent over $k$; it follows that $T(t) = 1$ and $N(t) = 0$, so $t + \overline{t} = 1$.
\end{remark}
\begin{remark}\label{rem:splitC1}
	There are three possibilities for the composition algebra $C_1$ and the element $t_1 \in C_1$, up to isomorphism:
	\begin{enumerate}
		\item\label{it:splitC1:1} $C_1 = k \oplus k$ and $t_1 = (1,0)$;
		\item\label{it:splitC1:2} $C_1 = M_2(k)$ and $t_1 = \idempotentmatrix$;
		\item\label{it:splitC1:3} $C_1 = \Zorn(k)$ and $t_1 = \idempotentmatrix$.
	\end{enumerate}
	(See also \cref{ex:split comp}.)
	Correspondingly, the algebra $\A \cong C_1 \otimes C_2$ can be written as $C_2 \oplus C_2^\op$ in case \ref{it:splitC1:1}, as $M_2(C_2)$ in case \ref{it:splitC1:2}, and as $\Zorn(k) \otimes C_2$ in case \ref{it:splitC1:3}.
\end{remark}

\subsection{The Smirnov algebra}\label{ss:ex6}

\begin{proposition}
	The $35$-dimensional Smirnov algebras over $k$ do not admit conjugate\dash supplementary idempotents.
\end{proposition}
\begin{proof}
	Let $\A$ be the Smirnov algebra $\A = T(C)$ associated to the octonion algebra $C$ over $k$.
	By \cite[(1.4) and Proposition 1.9]{Allison1993}, $\A$ is isomorphic to the $35$-dimensional subalgebra of $\widehat\A = C \otimes_k C$ given by
	\[ \A :=  \bigl\langle a \otimes a - N(a) 1 \otimes 1 \mid a \in C \bigr\rangle \]
	and has
	\[ \SS = \{ s \otimes 1 + 1 \otimes s \mid s \in C, \ \overline{s} = -s \} \subset \A . \]
	Note, however, that $1 \otimes 1 \not\in \A$; in fact, by \cite[p.\@~48]{Allison1993}, if $(e_1,\dots,e_8)$ and $(e_1^*,\dots,e_8^*)$ are two bases of $C$ dual with respect to the polar form $N(x,y) := N(x+y) - N(x) - N(y)$, then
	\[ 1_{\A} = 1 \otimes 1 - \tfrac14 \textstyle\sum_{i=1}^8 e_i \otimes e_i^* . \]
	In particular, the rank (as a tensor) of $1_\A$ is $7$ or $8$. (It is $7$ only when $\Char(k) = 7$, but this is not relevant.)
	
	Suppose now that $t \in \A$ is a conjugate-supplementary idempotent, and let $s := 2t - 1_\A$ be the corresponding skew element with $s^2 = 1_\A$ as in \cref{le:csi to skew}.
	Write $s = s_0 \otimes 1 + 1 \otimes s_0$; then $s^2 = 2 \bigl( s_0 \otimes s_0 - N(s_0) 1 \otimes 1 \bigr)$.
	This is a contradiction with the previous paragraph, because this element has rank at most $2$.
\end{proof}

\section{Classification of the central simple examples}\label{se:classification}

Now that we have determined the possible conjugate-supplementary idempotents in all central simple structurable algebras in \cref{se:examples}, we proceed with their explicit description as $AG_2$-algebras as in \cref{def:AG2-algebra}.
We organize them in a different order compared to \cref{se:examples}.
We also investigate whether the modules $J$ and $J'$ contain non-zero skew elements, and whether the maps $\times$ and $\times'$ are non-zero.

\begin{example}\label{ex:B+B}
	Let $B$ be a central simple alternative algebra.
	By Kleinfeld's classification theorem \cite{Kleinfeld1953}, $B$ is either a central simple associative $k$-algebra or an octonion algebra over $k$; see  also \cite[Theorem 3.17]{Schafer1966}.
	
	Let $\A$ be the structurable algebra $\A = B \oplus B^\op$ equipped with the exchange involution.
	\begin{enumerate}
		\item\label{it:B+B:central}
			The conjugate-supplementary idempotent $t = (1, 0)$ gives rise to an $AG_2$-algebra $(B, J, J', \times, \times', T)$, where
			\[ \A = \begin{pmatrix}	
				B & 0 \\ 0 & B^\op
		
		\end{pmatrix} , \]
			with $J = J' = 0$, so also the maps $\times$, $\times'$ and $T$ are the zero maps.
		\item\label{it:B+B:non-central}
			On the other hand, let $t = (u, \overline{1-u})$, where $u \in B$ is any non-trivial idempotent. Write $v := 1-u$.
			Then this gives rise to an $AG_2$-algebra $(C, J, J', \times, \times', T)$, where
			\[ \A = \begin{pmatrix}	
				uBu \oplus \overline{vBv} & uBv \oplus \overline{uBv} \\ vBu \oplus \overline{vBu} & vBv \oplus \overline{uBu}
				\end{pmatrix} , \]
			with $C = uBu \oplus \overline{vBv} \subseteq B \oplus B^\op$.
			In fact, we can write $\A = B \oplus B^\op$ with
			\[ B = \begin{pmatrix}	
				uBu & uBv \\ vBu & vBv
				\end{pmatrix} , \]
			which is the usual Peirce decomposition with respect to the arbitrary non-trivial idempotent $u \in B$.
			Recall that the Peirce subspace $uBv$ is always non-trivial.
			(Indeed, if $uBv = 0$, then $Bu$ would be a two-sided ideal in $B$, contradicting the simplicity of $B$.)
			In particular, the modules $J = uBv \oplus \overline{uBv}$ and $J' = vBu \oplus \overline{vBu}$ are always non-zero.
			Moreover, $J$ and $J'$ always contain non-zero skew elements, since any element of the form $(x, \overline{-x})$ is skew.
			
			When $B$ is associative, the maps $\times$ and $\times'$ are the zero maps.
			In contrast, when $B$ is (necessarily split) octonion, there is only one choice for $u$, up to isomorphism, namely $u = \idempotentmatrix$ (see also \cref{rem:comp-idemp});
			the resulting modules $J$ and $J'$ are isomorphic to $k^3 \oplus k^3$. The map $\times$ is given by the usual cross product on the first of the two $k^3$-components and by the negative of the usual cross product on the second component; the map $\times'$ is given by the negative of the map $\times$.
	\end{enumerate}
\end{example}
\begin{remark}
	\begin{enumerate}
		\item 
			If $B$ admits an involution $x \mapsto \overline{x}$, then $B^\op \cong B$ and the algebra $\A = B \oplus B^\op$ with the exchange involution is isomorphic to the algebra $\A' := B \oplus B$ equipped with the involution $(a,b) \mapsto (\overline{b}, \overline{a})$. This algebra $\A'$ with this involution can also be viewed as the tensor product $(k \oplus k) \otimes_k B$ of the two algebras with involution $k \oplus k$ (with the exchange involution) and $B$ (with the involution $x \mapsto \overline{x}$).
		\item
			Assume now that $B$ is an octonion algebra.
			If $t$ is a central idempotent $t = (1, 0)$ as in \cref{ex:B+B}\cref{it:B+B:central}, then $\A' \cong (k \oplus k) \otimes B$ and $t$ are as in \cref{rem:splitC1}\cref{it:splitC1:1}.
			On the other hand, if $t$ is non-central, say $t = (u, \overline{1-u})$ with $u = \idempotentmatrix$ as in \cref{ex:B+B}\cref{it:B+B:non-central}, then $\A' \cong \Zorn(k) \otimes (k \oplus k)$ and $t$ are as in \cref{rem:splitC1}\cref{it:splitC1:3}.
	\end{enumerate}
\end{remark}

\begin{example}\label{ex:M2B}
	Let $B$ be a central simple alternative $K$-algebra, equipped with an involution $\sigma \colon b \mapsto \overline{b}$, where either $K=k$ and the involution is of the first kind, or $K/k$ is a quadratic extension and the involution is of the second kind (fixing $k$ but not $K$).
	Let $\epsilon = \pm 1$.
	If $B$ is an octonion algebra, we require $K=k$, $\epsilon = -1$, and $\sigma$ equal to the standard involution.
	
	Now let $\A$ be the structurable algebra $M_2(B)$ equipped with the involution
	\[ \overline{\begin{pmatrix}
			a & b \\ c & d
		\end{pmatrix}} = \begin{pmatrix}
			\overline{d} & \epsilon \overline{b} \\ \epsilon\overline{c} & \overline{a}
		\end{pmatrix} . \]
	Then $\A$ with corresponding conjugate-supplementary idempotent $t = \idempotentmatrix$ gives rise to an $AG_2$-algebra $(B, J, J', \times, \times', T)$, where
	\[ \A = \begin{pmatrix}	
		B & B \\ B & B
	\end{pmatrix} , \]
	with $J = J' = B$.
	(Compared to the general form of $\A$ in \cref{def:AG2-algebra}, the bottom-right entry is $B$ rather than $\overline{B}$; the two are identified via the involution of~$B$.)
	
	If $B = k$, then the case $\epsilon=1$ results in no non-zero skew elements in $J$ and $J'$, whereas the case $\epsilon = -1$ results in the situation where all elements in $J$ and $J'$ are skew.
	If $B \neq k$, then there are always non-zero skew elements in $J$ and $J'$.
	Observe that the case without non-zero skew elements in $J$ is precisely the case where $\A$ has skew-dimension $1$.

	In this example, the maps $\times$ and $\times'$ are always the zero maps, even when $B$ is octonion, simply because the multiplication in $\A$ is given by the usual matrix multiplication of $2 \times 2$-matrices (with entries in $B$).
\end{example}
\begin{remark}
	\begin{enumerate}
		\item If $B$ is an octonion algebra over $k$ and we use $\epsilon = 1$ or $\sigma$~a non-standard involution in \cref{ex:M2B}, then the resulting algebra $\A = M_2(B)$ with this involution is no longer a structurable algebra. (This follows from the classification, but it can also be checked directly that the defining identity from \cref{def:struct algebra} fails in this case.)
			Also observe that if $B$ is octonion over $k$ (with its standard involution) and $\epsilon = -1$, then the structurable algebra $\A = M_2(B)$ is isomorphic to $M_2(k) \otimes B$, where both $M_2(k)$ and $B$ are equipped with their standard involution, so we are in the situation of \cref{rem:splitC1}\cref{it:splitC1:2}.
		\item If $K/k$ is a quadratic extension (and thus the involution is of the second kind), then the choice of $\epsilon$ is irrelevant. Indeed, an explicit isomorphism from the algebra with $\epsilon = 1$ to the algebra with $\epsilon = -1$ is given by
		\[ \begin{pmatrix}
			a & b \\ c & d
		\end{pmatrix} \mapsto \begin{pmatrix}
			a & ub \\ u^{-1}c & d
		\end{pmatrix} \]
		where $u$ is an invertible skew element of $K$.

		On the other hand, when $K=k$ (and thus the involution is of the first kind), then the choice $\epsilon=1$ gives the involution on $\A$ the same type as that on $B$, whereas the choice $\epsilon=-1$ gives the opposite type.
	\end{enumerate}
\end{remark}

In our next two examples, we deal with structurable algebras from a non-degenerate hermitian form defined over an associative algebra $E$. We first consider the case where $E$ is as in \cref{ex:B+B}, and then the case where $E$ is as in \cref{ex:M2B}.

\begin{example}\label{ex:herm:B+B}
	Let $E = B \oplus B^\op$ for some associative central simple $k$-algebra~$B$.
	Let $M \neq 0$ be a right $E$-module, let $h \colon M \times M \to E$ be a non-degenerate hermitian form, and let $\A = E \oplus M$ be the corresponding structurable algebra as in \cref{def:stralg hermitian}.
	Then $\A = B \oplus B^\op \oplus M$.
	By \cref{pr:csi herm}, any conjugate-supplementary idempotent $t \in \A$ belongs to $E$,
	and by \cref{pr:assoc idemp}\ref{it:assoc idemp:2}, $t = (u, \overline{1-u}) \in E = B \oplus B^\op$ for some idempotent $u \in B$.
	Recall also that by \cref{pr:Mt=Mbart}, $M = Mt \oplus M\overline{t}$ with $Mt \cong M\overline{t}$.
	\begin{enumerate}
		\item
			The case $u=1$ gives rise to an $AG_2$-algebra $(B, J, J', \times, \times', T)$, where
			\[ \A = \begin{pmatrix}	
				B & Mt \\ M\overline{t} & B^\op			
			\end{pmatrix} , \]
			with $J = Mt$ and $J' = M\overline{t}$. Since the involution of $\A$ fixes $M$, the modules $J$ and $J'$ do not contain non-zero skew elements.
			Moreover, the map $\times \colon J \times J \to J'$ is given by the product map in $\A$, so by definition, we have
			\[ mt \times nt = h(nt, mt) = \overline{t} h(n,m) t = 0 \]
			because $t$ is central in $E$ and $\overline{t} t = 0$. Similarly, the map $\times'$ is the zero map.
		\item
			Now let $u \in B$ be a non-trivial idempotent and write $v := 1-u$.
			Using \cref{ex:B+B}, we get an $AG_2$-algebra $(C, J, J', \times, \times', T)$ with
			\[ \A = \begin{pmatrix}	
				uBu \oplus \overline{vBv} \oplus 0 & uBv \oplus \overline{uBv} \oplus Mt \\ vBu \oplus \overline{vBu} \oplus M\overline{t} & vBv \oplus \overline{uBu} \oplus 0
				\end{pmatrix} , \]
			with $C = uBu \oplus \overline{vBv} \subseteq B \oplus B^\op$.
			We can also write this as $\A = B \oplus B^\op \oplus M$ with
			\[ B = \begin{pmatrix}	
				uBu & uBv \\ vBu & vBv
				\end{pmatrix} \quad \text{and} \quad M = \begin{pmatrix}	
				0 & Mt \\ M\overline{t} & 0
				\end{pmatrix} .
			\]
			As in \cref{ex:B+B}, the modules
			\[ J = uBv \oplus \overline{uBv} \oplus Mt \quad \text{and} \quad J' = vBu \oplus \overline{vBu} \oplus M\overline{t} \]
			always contain non-zero skew elements, namely those of the form $(x, \overline{-x}, 0)$.
			Moreover, the maps $\times$ and $\times'$ are non-zero.
			Explicitly, we have
			\[ (x_1, \overline{y_1}, m_1) \times (x_2, \overline{y_2}, m_2) =
				\bigl( h(m_2, m_1), \ m_2(x_1, \overline{y_1}) + m_1(y_2, \overline{x_2}) \bigr) \]
			for all $x_1, x_2, y_1, y_2 \in uBv$ and all $m_1, m_2 \in Mt$.
			Observe that indeed
			\begin{align*}
				h(m_2, m_1) &\in \overline{t} E t = vBu \oplus \overline{vBu} \quad \text{and} \\
				m_2(x_1, \overline{y_1}) + m_1(y_2, \overline{x_2}) &\in Mt \cdot tE\overline{t} \subseteq M\overline{t} .
			\end{align*}
	\end{enumerate}
\end{example}

\begin{example}\label{ex:herm:M2B}
	Let $E = M_2(B)$ as in \cref{pr:assoc idemp}\cref{it:assoc idemp:1}.
	Let $M \neq 0$ be a right $E$\dash module, let $h \colon M \times M \to E$ be a non-degenerate hermitian form, and let $\A = E \oplus M$ be the corresponding structurable algebra as in \cref{def:stralg hermitian}.
	Then $\A = M_2(B) \oplus M$.
	By \cref{pr:csi herm}, any conjugate-supplementary idempotent $t \in \A$ belongs to~$E$,
	and by \cref{pr:assoc idemp}\ref{it:assoc idemp:1}, we may assume that $t = \idempotentmatrix$.
	By \cref{pr:Mt=Mbart} again, $M = Mt \oplus M\overline{t}$ with $Mt \cong M\overline{t}$.
	Then $\A$ with the conjugate-supplementary idempotent $t$ gives rise to an $AG_2$\dash algebra $(B, J, J', \times, \times', T)$, where
	\[ \A = \begin{pmatrix}	
		B & B \oplus Mt \\ B \oplus M\overline{t} & B
	\end{pmatrix} , \]
	with $J = B \oplus Mt$ and $J' = B \oplus M\overline{t}$.
	Since the involution of $\A$ fixes $M$, the skew elements of $J$ and $J'$ are exactly as in \cref{ex:M2B}.
	In particular, $J$ and $J'$ contain non-zero skew elements unless $B=k$ and $\epsilon=1$ (in the notation of \cref{ex:M2B}).
	Again, observe that the case without non-zero skew elements in $J$ is precisely the case where $\A$ has skew-dimension $1$.
	
	Moreover, the maps $\times$ and $\times'$ are non-zero.
	Explicitly, we have
	\[ (b, m) \times (c, n) = \Bigl( h(n,m), \ n.\begin{psmallmatrix} 0 & b \\ 0 & 0 \end{psmallmatrix} + m.\begin{psmallmatrix} 0 & \epsilon \overline{c} \\ 0 & 0 \end{psmallmatrix} \Bigr) \]
	for all $b,c \in B$ and all $m,n \in Mt$. Observe that indeed
	\begin{align*}
		h(n, m) &\in \overline{t} E t = \begin{psmallmatrix} 0 & 0 \\ B & 0 \end{psmallmatrix} \quad \text{and} \\
		n.\begin{psmallmatrix} 0 & b \\ 0 & 0 \end{psmallmatrix} + m.\begin{psmallmatrix} 0 & \epsilon \overline{c} \\ 0 & 0 \end{psmallmatrix} &\in M\overline{t} .
	\end{align*}
\end{example}

\begin{example}\label{ex:skewdim1}
	Let $\A = M(J, \eta)$ be the structurable algebra of skew-dimension $1$ from \cref{ex:SMA} and let $t = \idempotentmatrix \in \A$.
	Then $\A$ with the conjugate-supplementary idempotent $t$ gives rise to an $AG_2$\dash algebra $(k, J, J', \times, \times', \eta T)$, where
	\[ \A = \begin{pmatrix}	
		k & J \\ J & k
	\end{pmatrix} , \]
	with $J = J'$.
	(Notice that the ``trace map'' of the $AG_2$-algebra is $\eta T$, where $T$ is as in \cref{ex:SMA}.)
	Since the involution of $\A$ fixes $J$ and $J'$, there are no non-zero skew elements in $J$ and $J'$.
	The map $\times$ is precisely the map called $\times$ in \cref{ex:SMA}, while the map $\times'$ is given by
	\[ v' \times' w' := \eta (v' \times w') \]
	for all $v', w' \in J'$; see the multiplication formula in \cref{ex:SMA}.
	In particular, the maps $\times$ and $\times'$ are non-zero in case \cref{ex:SMA:1} of \cref{ex:SMA}, but they are the zero maps in cases \cref{ex:SMA:2,ex:SMA:3}.
\end{example}

\begin{example}\label{ex:C1xC2}
	Let $C_1 = \Zorn(k)$ and let $C_2$ be an arbitrary composition algebra over $k$ (of dimension $1$, $2$, $4$ or $8$).
	Let $\A = C_1 \otimes C_2$, equipped with the tensor product of the standard involutions on $C_1$ and $C_2$,
	and let $t = t_1 \otimes 1$, where $t_1 = \idempotentmatrix \in C_1$.
	
	Then $\A$ with the conjugate-supplementary idempotent $t$ gives rise to an $AG_2$\dash algebra $(C_2, J, J', \times, \times', T)$, where
	\[ \A = \begin{pmatrix}	
		C_2 & k^3 \otimes C_2 \\ k^3 \otimes C_2 & C_2^\op
	\end{pmatrix} , \]
	with $J = J' = k^3 \otimes C_2$.
	Since the involution acts as $-1$ on $k^3$ and as the standard involution on $C_2$, the subspace of skew elements of $J$ and $J'$ is equal to $k^3 \otimes 1$.
	Moreover, the maps $\times$ and $\times'$ are non-zero.
	Indeed, since $\A$ is the tensor product $\Zorn(k) \otimes C_2$ and the multiplication on $\Zorn(k)$ restricted to $J = k^3 \otimes 1$ is given by the usual ``cross product'' (see again \cite[p.\@~158, Example 2.4.2]{McCrimmon2004}), the map $\times$ is the scalar extension of this cross product to $C_2$. Explicitly, if we write $J = J' = C_2 \oplus C_2 \oplus C_2$, we get
	\[ (a_1, a_2, a_3) \times (b_1, b_2, b_3) = (a_2 b_3 - a_3 b_2, \ a_3 b_1 - a_1 b_3, \ a_1 b_2 - a_2 b_1) \]
	for all $a_i, b_i \in C_2$. On the other hand, the map $\times' \colon J' \times J' \to J$ is given by the negative of this formula, i.e.,
	\[ (a_1, a_2, a_3) \times' (b_1, b_2, b_3) = (- a_2 b_3 + a_3 b_2, \ - a_3 b_1 + a_1 b_3, \ - a_1 b_2 + a_2 b_1) \]
	for all $a_i, b_i \in C_2$.
	Note that the anti-symmetry law $w \times v = - v \times w$ from $\Zorn(k)$ does \emph{not} hold for $\Zorn(k) \otimes C_2$ when $C_2$ is non-commutative.
	It is also worth recording that the trace map $T \colon J \times J' \to C_2$ of the $AG_2$-algebra is given by
	\[ T \bigl( (a_1, a_2, a_3),\ (b_1, b_2, b_3) \bigr) = a_1 b_1 + a_2 b_2 + a_3 b_3 \]
	for all $a_i, b_i \in C_2$.
\end{example}

\begin{remark}
	The following overview illustrates that we have now covered all possible examples.
	(Notice that there are no examples in \S\ref{ss:ex2} and \S\ref{ss:ex6}, and that the examples in \S\ref{ss:ex5} where both $C_1$ and $C_2$ are associative---i.e., have dimension $1$, $2$ or $4$---are covered by \S\ref{ss:ex1}.)

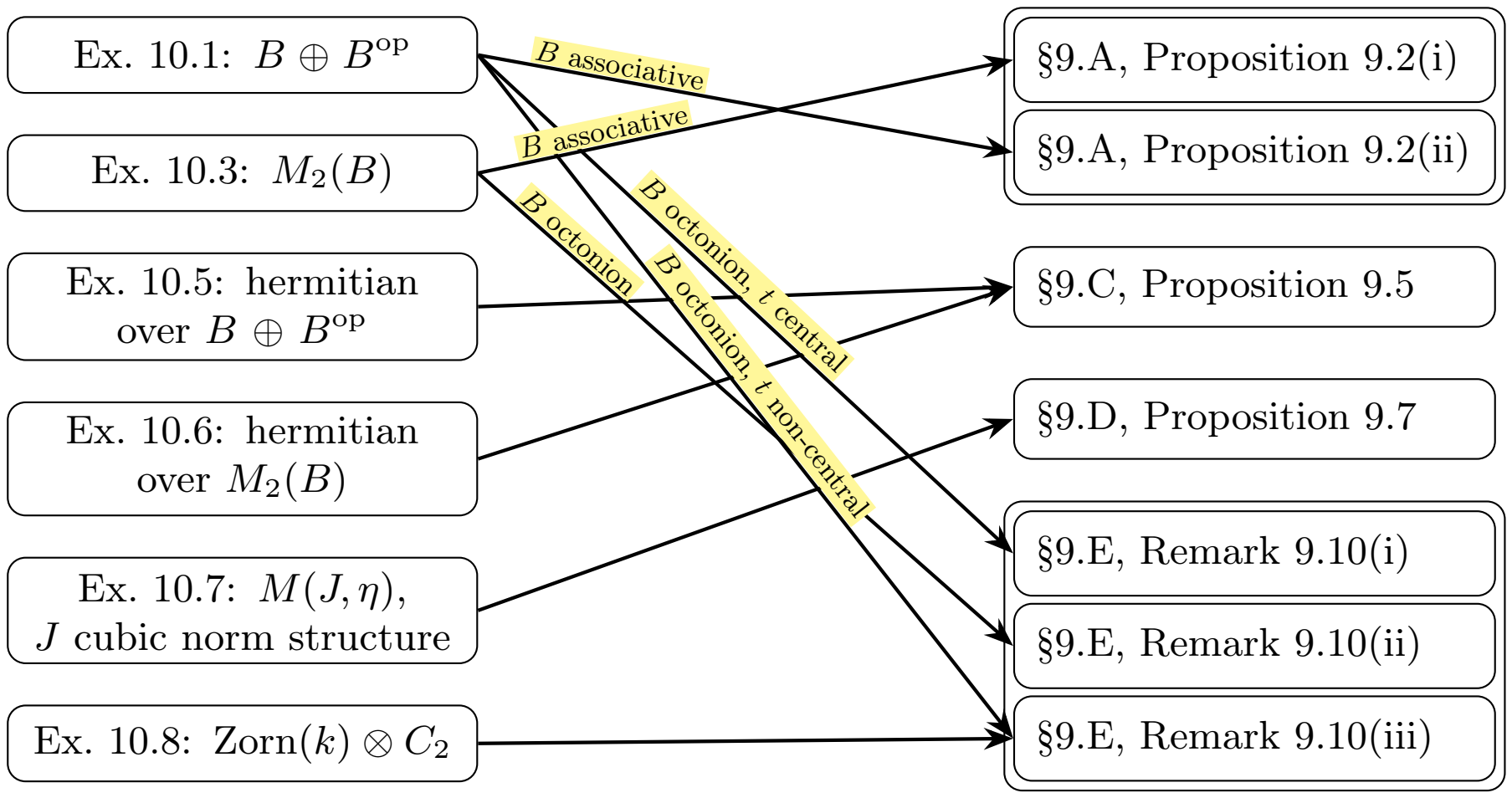
\begin{center}
\begin{tikzpicture}
\small 
\tikzset{
  itemleft/.style={
    draw, rounded corners,
    align=center,
    text width=3.5cm,
    minimum height=6mm,
    inner sep=3pt
  },
  itemright/.style={
    draw, rounded corners,
    align=left,
    text width=3.45cm,
    minimum height=6mm,
    inner sep=5pt
  },
  group/.style={
    draw,
    rounded corners,
    inner sep=2pt
  },
  arrowtext/.style={
  	sloped,
  	above=.01,
  	scale=.72,
  	fill=yellow!50,
  	inner sep=1pt
  },
  arrow/.style={-{Stealth}, thick}
}

\matrix[row sep=3.2mm] (L) {
  \node[itemleft] (L1) {Ex.\@~\ref{ex:B+B}: $B \oplus B^\op$}; \\
  \node[itemleft] (L2) {Ex.\@~\ref{ex:M2B}: $M_2(B)$}; \\
  \node[itemleft] (L3) {Ex.\@~\ref{ex:herm:B+B}: hermitian \\ over $B \oplus B^\op$}; \\
  \node[itemleft] (L4) {Ex.\@~\ref{ex:herm:M2B}: hermitian \\ over $M_2(B)$}; \\
  \node[itemleft] (L5) {Ex.\@~\ref{ex:skewdim1}: $M(J, \eta)$, \\ $J$ cubic norm structure}; \\
  \node[itemleft] (L6) {Ex.\@~\ref{ex:C1xC2}: $\Zorn(k) \otimes C_2$}; \\
};

\matrix[row sep=4mm, right=4cm of L] (R) {
  \node[itemright] (R1a) {\S\ref{ss:ex1}, \cref{pr:assoc idemp}\cref{it:assoc idemp:1}}; \\[-3.5mm]
  \node[itemright] (R1b) {\S\ref{ss:ex1}, \cref{pr:assoc idemp}\cref{it:assoc idemp:2}}; \\
  \node[itemright] (R3) {\S\ref{ss:ex3}, \cref{pr:csi herm}}; \\
  \node[itemright] (R4) {\S\ref{ss:ex4}, \cref{pr:skew-dim one}}; \\
  \node[itemright] (R5a) {\S\ref{ss:ex5}, \cref{rem:splitC1}\cref{it:splitC1:1}}; \\[-3.5mm]
  \node[itemright] (R5b) {\S\ref{ss:ex5}, \cref{rem:splitC1}\cref{it:splitC1:2}}; \\[-3.5mm]
  \node[itemright] (R5c) {\S\ref{ss:ex5}, \cref{rem:splitC1}\cref{it:splitC1:3}}; \\
};

\node[group, fit=(R1a)(R1b)] (G1) {};
\node[group, fit=(R5a)(R5b)(R5c)] (G5) {};

\draw[arrow] (L3.east) -- (R3.west);
\draw[arrow] (L4.east) -- (R3.west);
\draw[arrow] (L5.east) -- (R4.west);
\draw[arrow] (L6.east) -- (R5c.west);
\draw[arrow] (L1.east) -- node[pos=0.26, arrowtext] {$B$ associative} (R1b.west);
\draw[arrow] (L1.east) -- node[pos=0.47, arrowtext] {$B$ octonion, $t$ central} (R5a.west);
\draw[arrow] (L2.east) -- node[pos=0.17, arrowtext] {$B$ octonion} (R5b.west);
\draw[arrow] (L1.east) -- node[pos=0.5, arrowtext] {$B$ octonion, $t$ non-central} (R5c.west);
\draw[arrow] (L2.east) -- node[pos=0.24, arrowtext] {$B$ associative} (R1a.west);

\end{tikzpicture}
\end{center}

\end{remark}

We can summarize our results in the following classification theorem.
\begin{theorem}\label{th:classification}
	Let $\A$ be a finite-dimensional central simple structurable algebra over $k$ and let $t \in \A$ be a conjugate-supplementary idempotent.
	Assume that $\Char(k) \neq 5$, or that $\A$ belongs to one of the classes \ref{cssa:1}--\ref{cssa:6} listed at the beginning of \cref{se:examples}.
	Let $\Xi = (B, J, J', \times, \times', T)$ be the associated $AG_2$-algebra.
	Then $\Xi$ is either weak or strong. Moreover:
	\begin{enumerate}
		\item\label{th:classification:feeble}\label{it:class:1} If $\Xi$ is feeble, then up to isomorphism, $\A = B \oplus B^\op$ (with $B$ either associative or octonion) and $t = (1, 0)$.
		\item\label{it:class:2} If $\Xi$ is weak but not feeble, then up to isomorphism, one of the following holds:
			\begin{enumerate}[\rm (a)]
				\item\label{it:class:2a} $\A$ is of hermitian type over $E = B \oplus B^\op$ and $t = (1, 0) \in E$;
				\item\label{it:class:2b} $\A = M(J, \eta)$ for some (possibly degenerate) non-zero cubic norm structure $J$ and $t = \idempotentmatrix$.
			\end{enumerate}
		\item\label{it:class:3} If $\Xi$ is strong but not mighty, then up to isomorphism, one of the following holds:
			\begin{enumerate}[\rm (a)]
				\item\label{it:class:3a} $\A = B \oplus B^\op$ with $B$ associative, and $t$ is non-central;
				\item\label{it:class:3b} $\A = M_2(B)$ (with $B$ either associative or octonion), where if $B=k$ or $B$ is octonion, the involution is standard; $t = \idempotentmatrix$.
			\end{enumerate}
		\item\label{it:class:4} If $\Xi$ is mighty, then up to isomorphism, one of the following holds:
			\begin{enumerate}[\rm (a)]
				\item\label{it:class:4a} $\A$ is of hermitian type over $E = B \oplus B^\op$ and $t$ is a non-central idempotent in $E$;
				\item\label{it:class:4b} $\A$ is of hermitian type over $E = M_2(B)$ and $t = \idempotentmatrix \in E$, where if $B=k$, the involution is standard;
				\item\label{it:class:4c} $\A = \Zorn(k) \otimes C$ for some composition algebra $C$ (of dimension $1$, $2$, $4$ or $8$) and $t = \idempotentmatrix \otimes 1$.
			\end{enumerate}
	\end{enumerate}
\end{theorem}
\begin{proof}
	This follows from the discussion in \cref{ex:B+B,ex:M2B,ex:herm:B+B,ex:herm:M2B,ex:skewdim1,ex:C1xC2}.
	Note that the mighty $AG_2$-algebra in \cref{ex:B+B} (with $B$ octonion and $t$ non-central) is captured by \cref{ex:C1xC2}, so this is part of item~\ref{it:class:4}\ref{it:class:4c}.
\end{proof}

It is now easy to check for which of these examples the $(AG_2, A_2)$-grading extends to an $(AG_2, G_2)$-grading, or equivalently, by \cref{th:G2-frame}, which of the corresponding pairs $(\A, t)$ admit a $G_2$-frame.
\begin{corollary}\label{co:G2-frame}
	Let $\A$ and $t \in \A$ be as in \cref{th:classification}.
	Let $(B, J, J', \times, \times', T)$ be the associated $AG_2$-algebra.
	Then $(\A, t)$ admits a $G_2$-frame if and only if we are in one of the following cases of \cref{th:classification}:
	\begin{itemize}
		\item[\ref{it:class:2}\ref{it:class:2b}] $\A = M(J, \eta)$ and there exists some $x \in J$ with $N(x) = \eta^{-1}$;
		\item[\ref{it:class:4}\ref{it:class:4a},\ref{it:class:4b}] $\A$ is of hermitian type over $E$ and there exists $m \in Mt$ such that for $s := \tfrac{1}{2} h(m,m) \in \overline{t}Et$, there exists $r \in tE\overline{t}$ with $rs = t$ and $sr = \overline{t}$;
		\item[\ref{it:class:4}\ref{it:class:4c}] $\A = \Zorn(k) \otimes C$ for some composition algebra $C$ of dimension $4$ or $8$.
	\end{itemize}
\end{corollary}
\begin{proof}
	The existence of a $G_2$-frame implies, in particular, that the maps $\times$ and $\times'$ are non-zero.
	This already excludes cases \ref{it:class:1}, \ref{it:class:2}\ref{it:class:2a} and \ref{it:class:3}.
	\begin{itemize}
		\item[\ref{it:class:2}\ref{it:class:2b}]
			Write $u = \begin{psmallmatrix} 0 & x \\ 0 & 0 \end{psmallmatrix}$ and $u' = \begin{psmallmatrix} 0 & 0 \\ y & 0 \end{psmallmatrix}$.
			The equations from \cref{def:G2-frame} become
			\[ x \times x = 2y , \qquad \eta(y \times y) = 2x, \qquad \eta T(x,y) = 3 . \]
			For a cubic norm structure, $x \times x = 2x^\sharp$, so necessarily $y = x^\sharp$. Using
			\[ (x^\sharp)^\sharp = N(x) x , \qquad T(x, x^\sharp) = 3 N(x) , \]
			the remaining equations are equivalent to $\eta N(x) = 1$.
		\item[\ref{it:class:4}\ref{it:class:4a},\ref{it:class:4b}] Assume now that $\A$ is of hermitian type over $E$, so $\A = E \oplus M$ as in \cref{def:stralg hermitian}.
			Write $u := (r, m) \in J$ and $u' := (s, n) \in J'$ with $r,s \in E$ and $m,n \in M$.
			Notice that $r \in tE\overline{t}$ and $s \in \overline{t}Et$; in particular, $r^2 = s^2 = 0$.
			The equations from \cref{def:G2-frame} become
			\begin{alignat*}{2}
				r &= \overline{r}, &\quad s &= \overline{s} , \\
				\bigl( h(m,m),\ 2mr \bigr) &= 2(s, n), &\quad \bigl( h(n,n),\ 2ns \bigr) &= 2(r, m), \\
				\bigl( rs + h(n,m), \  nr + ms \bigr) &= (3t, 0), &\quad \bigl( sr + h(m,n), \  ms + nr \bigr) &= (3\overline{t}, 0). 
			\end{alignat*}
			We immediately see that $n = mr$, $m = ns$, $s = \tfrac{1}{2} h(m,m)$, and $r = \tfrac{1}{2} h(n,n)$.
			Next, $3t = rs + h(n,m) = rs + h(mr,m) = rs + rh(m,m) = 3rs$, so $rs = t$; similarly, $sr = \overline{t}$.
			
			Conversely, if there exists an element $m \in Mt$ such that $s := \tfrac{1}{2} h(m,m) \in \overline{t}Et$ has the property that there exists an element $r \in tE\overline{t}$ with $rs = t$ and $sr = \overline{t}$, then setting $n := mr$ gives a solution to the $G_2$-frame equations.
		\item[\ref{it:class:4}\ref{it:class:4c}] Assume, finally, that $\A = \Zorn(k) \otimes C$ for some composition algebra $C$.
			If $C$ is commutative, then $u \times u = 0$ for all $u \in C^3$, so no $G_2$-frame exists.
			
			Suppose now that $C$ is not commutative, i.e., $\dim C \in \{ 4, 8 \}$.
			Choose a quaternion subalgebra $D$ of $C$; then $D$ is generated by skew elements $p,q \in D$ satisfying
			\[ p^2 = \alpha, \qquad q^2 = \beta, \qquad pq = -qp \]
			with $\alpha,\beta \in k^\times$.
			Write $r := pq$, so $r^2 = -\alpha\beta$.  Now set
			\[ u := \bigl( p,\ q,\ -(\alpha\beta)^{-1} r \bigr) \quad \text{and} \quad u' := \bigl( \alpha^{-1} p, \ \beta^{-1} q, \ r \bigr)  . \]
			Then it follows from the formula for $\times$ and $T$ in \cref{ex:C1xC2} that the pair $(u, u')$ forms a $G_2$-frame.
		\qedhere
	\end{itemize}
\end{proof}

\begin{remark}\label{rem:isom-idemp}
	The condition appearing in the hermitian case that $t = rs$ and $\overline{t} = sr$ for certain $r,s \in E$ is precisely the defining condition for \emph{isomorphic idempotents} in the sense of \cite[Proposition (21.20)]{Lam01}.
	In particular, the left $E$-modules $Et$ and $E\overline{t}$ are isomorphic (as are the right $E$-modules $tE$ and $\overline{t}E$).
	It also follows that $tE\overline{t}$ and $\overline{t}Et$ are isomorphic as $k$-vector spaces, with (non-canonical) isomorphism
	\begin{align*}
		tE\overline{t} \to \overline{t}Et & \colon x \mapsto sxs , \quad \text{with inverse} \\
		\overline{t}Et \to tE\overline{t} & \colon y \mapsto ryr .
	\end{align*}
	Combined with the fact that $Mt \cong M\overline{t}$ (by \cref{pr:Mt=Mbart}), we conclude that $J \cong J'$.
	
	The condition that $t$ and $\overline{t} = 1-t$ be isomorphic idempotents is a genuine condition.
	See, for instance, \cref{ex:M3k} below.
\end{remark}

\begin{remark}\label{rem:Yaseen}
	We make the comparison with~\cite{Yaseen22} from the introduction more precise.
	Let $L$ be an $(AG_2, A_2)$-graded Lie algebra with grading algebra $\fG \cong \ssl_3(k)$, let $V$ be the natural $\fG$-module with highest weight $\omega_1$, and let $V^*$ be its dual.
	Yaseen decomposes $L$ into isotypic components as
	\[ L = (\fG \otimes A) \oplus (S^2V \otimes C) \oplus (S^2V^* \otimes C') \oplus (V^* \otimes E) \oplus (V \otimes E') \oplus D , \]
	where $D$ is the sum of the trivial submodules; see \cite[(2.2)]{Yaseen22}, where $V^*$ and $V$ appear as $\wedge^2 V$ and $\wedge^2 V^*$.
	\begin{enumerate}
		\item
			Let $(B, J, J', \times, \times', T)$ be the $AG_2$-algebra associated with $L$ as in \cref{th:AG2-A2}, so that $L/Z(L) \cong \TKK$.
			Comparing weights (see \cref{fig:weights,fig:AG2}) and using the action of the root vectors $t_+$ and $\overline{t}_-$ of $\fG$ (for instance, $[t_+, x_+] = (\overline{x} - x)_+$ for $x \in J$), one finds that the $S^2V$- and $V^*$-isotypic parts of the root space $J_+$ of weight $\omega_2$ are $(J \cap \SS)_+$ and $(J \cap \HH)_+$, respectively.
			It follows that
			\[ A \cong B, \quad C \cong J \cap \SS, \quad C' \cong J' \cap \SS, \quad E \cong J \cap \HH, \quad E' \cong J' \cap \HH , \]
			so the summands of Yaseen's coordinate algebra $\mathfrak{a} = A^+ \oplus A^- \oplus C \oplus C' \oplus E \oplus E'$ correspond to the hermitian and skew parts of $B \oplus B^\op$, $J$ and $J'$.
		\item
			Yaseen's additional assumption \cite[(3.1)]{Yaseen22} is that $[S^2V \otimes C, S^2V \otimes C] = 0$ and $[S^2V^* \otimes C', S^2V^* \otimes C'] = 0$.
			Since $S^2V^*$ is the only simple constituent of $S^2V \otimes S^2V$ occurring in $L$, and it occurs with multiplicity one, the first bracket is determined by its restriction to $(J \cap \SS)_+ \times (J \cap \SS)_+$, where $[x_+, y_+] = \psi(x,y)_+ = -2 (x \times y)_+$ by skew-alternativity.
			Hence this assumption amounts to the vanishing of $\times$ and $\times'$ on pairs of skew elements of $J$ and $J'$, respectively.
			Among the central simple examples, this excludes precisely the algebras $\Zorn(k) \otimes C$ from \cref{th:classification}\cref{it:class:4}\cref{it:class:4c}, where $J \cap \SS = k^3 \otimes 1$ and $\times$ is the cross product (see \cref{ex:C1xC2}); these are exactly the examples corresponding to the large exceptional Tits hexagons from \cite[\S 13]{Muhlherr26}.
	\end{enumerate}
\end{remark}

\section{Tits hexagons}\label{se:tits hex}

In this final section, we make the connection with the theory of Tits hexagons in \cite{Muhlherr26}, which was the original inspiration for our work.
As we will see, the connection with $(AG_2, \Delta)$-graded Lie algebras is not so straightforward, due to the specific nature of Tits polygons.
The $AG_2$-grading of a Lie algebra gives rise to root groups and commutator relations between non-opposite root groups, but the existence of a Tits hexagon also depends on the interaction between opposite root groups.


\medskip

The following way to pass from graded Lie algebras to groups is standard.
\begin{definition}\label{def:to groups}
	Let $L$ be an $(AG_2, A_2)$-graded Lie algebra.
	For each non-zero root $\alpha \in AG_2$, we define
	\[ \fu_\alpha := \begin{cases}
		L_\alpha & \text{if $\alpha$ is long or very long}, \\
		L_\alpha \oplus L_{2\alpha} & \text{if $\alpha$ is short},
	\end{cases} \]
	so $\fu_\alpha$ is at most two-step nilpotent. Now let
	\[ U_\alpha := \exp(\fu_\alpha) \leq \Aut(L) , \]
	where $\exp$ is the usual truncated exponential defined by
	\[ \exp(x) = \sum_{i=0}^4 \frac{1}{i!} (\ad x)^i \]
	for all $x \in \fu_\alpha$. (Recall that $\Char(k) \neq 2,3$.)
	
	In exponential coordinates, a short-root group $U_\alpha$ has multiplication
	\[ (x,s)(y,r) = \bigl( x + y,\, s + r + \tfrac{1}{2}[x,y] \bigr) \]
	for all $x, y \in L_\alpha$ and $s,r \in L_{2\alpha}$.
	The doubled-root subgroup is
	\[ V_\alpha := U_{2\alpha} \leq U_\alpha , \]
	where we set $V_\alpha = 1$ if $\alpha$ is long or very long.
	Now choose a cyclic ordering $\alpha_i$ (with $i \in \Z/12\Z$) of the roots of $\Phi \cong G_2 \subset AG_2$ such that even $i$ correspond to short roots, and write $U_i := U_{\alpha_i}$.
	This gives a canonical cyclic root-group sequence
	\[ \mathcal{U}(L) = \bigl( U_i \bigr)_{i \in \Z/12\Z}, \]
	together with the distinguished subgroups $V_i := V_{\alpha_i}$ for even $i$.
	Write
	\[ G := \langle U_i \mid i \in \Z/12\Z \rangle . \]
	We can then define the ``torus''
	\[ H := \{ g \in G \mid g(L_\alpha) = L_\alpha \text{ for all } \alpha \in AG_2^o \} . \]
	Notice that, in particular, $H$ normalizes each $U_\alpha$.
\end{definition}
We can now make the connection with Tits hexagons, which we will do using the machinery developed in \cite[\S 2]{Muhlherr19}.
In particular, we will use the notations from \cite[Notations 2.1 and 2.2]{Muhlherr19} and the (long) definition from \cite[Definition~2.3]{Muhlherr19} without repeating them here.
\begin{definition}\label{def:tits-hex}
	Let $L$ be an $(AG_2, A_2)$-graded Lie algebra and let $\mathcal{U}(L)$, $G$ and~$H$ be as in \cref{def:to groups}.
	Then we call $L$ \emph{Tits-hexagonalizable} if the data $\mathcal{U}(L)$ with torus $H$ forms a $\Phi_6$-grading as defined in \cite[Definition~2.3]{Muhlherr19} that is stable as defined in \cite[(2.4)]{Muhlherr19}.
	By \cite[Construction 5.2 and Theorem 5.3]{Muhlherr19}, there then exists a corresponding Tits hexagon (with root groups $U_i$).
	
	We say that this hexagonalization is \emph{of type $AG_2$} if, in addition, the assignments
	\[ \alpha_i \mapsto U_i, \qquad 2\alpha_i \mapsto V_i \quad \text{for even } i \]
	define an $AG_2$-commutator map in the sense of \cite[Definition 5.1]{Muhlherr26}.
	Thus the resulting Tits hexagon is \emph{$AG_2$-shaped} in the sense of \cite[Definition~5.6]{Muhlherr26}.
	(This is not automatic: a stable $\Phi_6$-grading can also give rise to a Tits hexagon of ordinary $G_2$\dash shape.)
	
	We also call a pair $(\A, t)$ (consisting of a structurable algebra $\A$ and a conjugate-supplementary idempotent $t \in \A$, as usual) \emph{Tits-hexagonalizable} if the corresponding $(AG_2, A_2)$-graded Lie algebra $L = K(\A)$ (as in \cref{th:AG2-A2-conv}) is Tits\dash hexagonalizable.
\end{definition}

Notice that if $L$ is Tits-hexagonalizable of type $AG_2$, then the corresponding $AG_2$-algebra is strong.

\begin{remark}\label{rem:stability}
	One of the conditions in \cref{def:tits-hex} is the \emph{stability condition} \cite[(2.4)]{Muhlherr19}.
	For algebraic structures over arbitrary base rings, this condition is rather delicate and certainly not automatically satisfied.
	(See, for instance, \cite[Observation~12.13]{Muhlherr26}.)
	
	For finite-dimensional central simple structurable algebras, however, the stability condition is automatic.
	Indeed, by \cite[Theorem 4.8]{Stavrova22}, the algebraic group $\mathbf G := \Aut(L)^\circ$ is adjoint absolutely simple.
	Let $\mathbf S$ be the two-dimensional split torus defining the $\mathbb Z^2$-grading of $L$.
	For each root $\alpha$, put
	\[ \mathbf S_\alpha:=(\ker\alpha)^\circ, \qquad \mathbf R_\alpha:=C_{\mathbf G}(\mathbf S_\alpha), \qquad \mathbf M:=C_{\mathbf G}(\mathbf S). \]
	The group $\mathbf R_\alpha$ is reductive.
	Choosing a complementary cocharacter in $\mathbf S$ gives opposite parabolic subgroups of $\mathbf R_\alpha$ with common Levi subgroup $\mathbf M$ and unipotent radicals whose groups of $k$-rational points are $U_\alpha$ and $U_{-\alpha}$.
	In particular, when $\alpha$ is short, these radicals include the doubled-root subgroups.

	Gauss decomposition \cite[Exp.~XXVI, Th\'eor\`eme~5.1]{SGA3-III} applied to the reductive group $\mathbf R_\alpha$ now gives
	\[ \mathbf R_\alpha(k) = U_{-\alpha}U_\alpha\mathbf M(k)U_{-\alpha} = U_{-\alpha}U_\alpha U_{-\alpha}\mathbf M(k), \]
	where the second equality holds because $\mathbf M$ normalizes $U_{-\alpha}$.
	Put $G_\alpha := \langle U_\alpha, U_{-\alpha} \rangle$ and $H_\alpha := G_\alpha\cap H$.
	Since preserving every homogeneous component is equivalent to centralizing $\mathbf S$, we have $H = G \cap \mathbf M(k)$.
	Moreover, if an element of $G_\alpha$ is written using the preceding decomposition, its Levi factor also belongs to $G_\alpha$.
	Consequently,
	\[ G_\alpha = U_{-\alpha} U_\alpha U_{-\alpha} H_\alpha , \]
	which is precisely the required stability condition.
\end{remark}

\medskip

The key ingredient of \cite[Definition 2.3]{Muhlherr19}, in addition to the setup of our paper and the stability condition, is its condition (iii), namely the existence of ``Weyl elements'':
for each root~$\alpha$, there exists some $m \in U_{-\alpha} U_\alpha U_{-\alpha} \subset G$ conjugating each root group $U_\beta$ to $U_{s_\alpha(\beta)}$ and normalizing $H$.
In particular, $U_\beta$ and $U_{-\beta}$ are conjugate inside $G$.
By \cref{rem:stability}, if $\A$ is finite-dimensional and central simple, then $(\A, t)$ is Tits-hexagonalizable if and only if every root has such a Weyl element.
It turns out that this forces the idempotent $t \in \A$ to be ``balanced'' in a certain sense, which will impose restrictions on the cases \cref{it:class:3}\cref{it:class:3a} and \cref{it:class:4}\cref{it:class:4a} from \cref{th:classification};
see \cref{pr:balanced} below.
In fact, we will see that among the strong examples from \cref{th:classification}, those in cases \cref{it:class:3}\cref{it:class:3b}, \cref{it:class:4}\cref{it:class:4b} and \cref{it:class:4}\cref{it:class:4c} are Tits-hexagonalizable, with one family of exceptions in case \cref{it:class:3}\cref{it:class:3b}, whereas in cases \cref{it:class:3}\cref{it:class:3a} and \cref{it:class:4}\cref{it:class:4a}, this happens precisely when $t$ is balanced; see \cref{co:hexagonalizable} at the end of this section.

\medskip

To compare our examples with the Tits hexagons from \cite{Muhlherr26}, we fix the following conventions.
\begin{notation}\label{not:coordinates}
	Let $\A$ be a structurable algebra with a conjugate-supplementary idempotent $t$, and let $L = \TKK$.
	As in the proof of \cref{th:AG2-A2-conv}, we identify the roots of $AG_2$ with the pairs $(i,j)$ for which $L_{ij} = L_i \cap L'_j$ is a root space (see \cref{ta:Lij}).
	As the cyclic ordering of the roots of the $G_2$-subsystem in \cref{def:to groups}, we always take
	\[ (\alpha_1, \alpha_2, \alpha_3, \alpha_4, \alpha_5, \alpha_6) = \bigl( (1,2), (1,1), (2,1), (1,0), (1,-1), (0,-1) \bigr) \]
	and $\alpha_{i+6} = -\alpha_i$, with indices modulo $12$, so that $U_i = U_{\alpha_i}$.
	As in \cref{def:TKK}, we write $x_+$ both for elements of $\A_+ = L_1$ and of $\SS_+ = L_2$, and similarly for $x_-$.
	We view $\Aut(L)$ as a group under composition, and we use the conventions $[g,h] = g^{-1}h^{-1}gh$ and $g^h = h^{-1}gh$ from \cite[1.5]{Muhlherr26}.
\end{notation}

Our description of the examples of hermitian type differs from the description in \cite[\S 12]{Muhlherr26}, so we will first translate their description in terms of $\epsilon$-pseudo-quadratic modules to our setting.
Notice that \cite{Muhlherr26} allows for arbitrary alternative rings whereas we only deal with alternative $k$-algebras over fields $k$ (of characteristic not $2$ or $3$).

\begin{construction}\label{con:MW}
	Let $\Lambda = (K, K_0, \sigma, \epsilon, V, q, f)$ be an $\epsilon$-pseudo-quadratic module as in \cite[Notation~12.4]{Muhlherr26} with $\epsilon \in \{1, -1\}$, where we allow $V = 0$; the module $V$ is denoted by~$L$ in \cite{Muhlherr26}.
	If $V = 0$, we also allow $K$ to be an octonion algebra with its standard involution~$\sigma$, with $\epsilon = -1$ and $K_0 = \Fix_K(\sigma)$, as in \cite[Observations~12.2 and~12.7]{Muhlherr26}.
	In particular, $V$ is a right $K$-module, $f \colon V \times V \to K$ is an $\epsilon$-hermitian form, i.e., $f(b,a) = \epsilon f(a,b)^\sigma$ for all $a, b \in V$, and $K_0 = \{ u \in K \mid u^\sigma = -\epsilon u \}$ because $\Char(k) \neq 2$.
	We assume that $f$ is non-degenerate, that $K$ is a $k$-algebra and $\sigma$ is $k$-linear, and that $(K, \sigma)$ is simple, i.e., $K$ has no non-trivial $\sigma$-stable ideals.
	Then, as explained in \cite[Observation~12.2]{Muhlherr26}, either
	\begin{enumerate}[\rm (a)]
		\item\label{con:MW:a} $K$ is a simple alternative $k$-algebra, or
		\item\label{con:MW:b} $K \cong C \oplus C^\op$ for some simple associative $k$-algebra $C$, and $\sigma$ is the exchange involution.
	\end{enumerate}
	(See also \cref{le:csa assoc}.)
	As in \cite[Proposition~12.5]{Muhlherr26}, let $U_\Lambda$ be the group with underlying set
	\[ U_\Lambda = \{ (a,r) \in V \times K \mid r - \tfrac{1}{2} f(a,a) \in K_0 \} \]
	and multiplication $(a,r)(b,s) = (a+b,r+s+f(a,b))$.

	Now let $E := M_2(K)$, equipped with the involution
	\[
	\overline{\begin{pmatrix}a&b\\c&d\end{pmatrix}}
	=
	\begin{pmatrix}
	d^\sigma & \epsilon b^\sigma \\
	\epsilon c^\sigma & a^\sigma
	\end{pmatrix}
	\]
	as in \cref{ex:M2B}.
	If $V = 0$, we set $\A := E$.
	If $V \neq 0$, we let $M := V \oplus V$ consist of two copies of~$V$, equipped with the usual right $E$-module structure and with the non-degenerate hermitian form
	\[ h\bigl((a,b),(c,d)\bigr) = \begin{pmatrix}
			f(b,c) & f(b,d) \\
			\epsilon f(a,c) & \epsilon f(a,d)
		\end{pmatrix},
	\]
	and we set $\A := E \oplus M$ as in \cref{def:stralg hermitian}.
	In both cases, we let $t := \idempotentmatrix \in E \subseteq \A$.
\end{construction}

\begin{proposition}\label{pr:MW}
	Let $\Lambda$ and $(\A, t)$ be as in \cref{con:MW}.
	Then the root-group sequence in \cref{def:to groups} has odd root groups isomorphic to $(K, +)$ and even root groups isomorphic to $K \times U_\Lambda$, and it satisfies the commutator relations in \cite[Proposition~12.5]{Muhlherr26}.
	If, in addition, $\A$ is finite-dimensional and central simple and $U_\Lambda$ contains an element $(a_0, r_0)$ with $r_0 \in K^\times$ (for instance, $(a_0, r_0) = (0,1)$ if $\epsilon = -1$), then $(\A, t)$ is Tits-hexagonalizable.
\end{proposition}
\begin{proof}
	This follows from \cref{ex:M2B,ex:herm:B+B,ex:herm:M2B} together with an explicit description of the root groups $U_i$ of $L := \TKK$, labeled as in \cref{not:coordinates}, which we now give.
	For $c \in K$, we write $c e_{ij} \in E = M_2(K)$ for the matrix with $c$ in position $(i,j)$ and zeros elsewhere, and for $a \in V$, we write $(a,0) \in Mt$ and $(0,a) \in M\ov{t}$ for the corresponding elements of $M = V \oplus V$ (if $V = 0$, these terms are simply absent).
	For $(a,r) \in U_\Lambda$, we put $r' := r - \tfrac{1}{2} f(a,a) \in K_0$.
	Now define, for all $\lambda \in K$ and $(a,r) \in U_\Lambda$,
	\begin{alignat*}{2}
		x_1(\lambda) &:= \exp\bigl( (\lambda^\sigma e_{11})_+ \bigr) && \in U_1, \\
		x_2(\lambda, a, r) &:= \exp\bigl( (\lambda^\sigma e_{21} + (0,a))_+ - (r' e_{12})_+ \bigr) && \in U_2, \\
		x_3(\lambda) &:= \exp\bigl( (\lambda^\sigma e_{11} - \lambda e_{22})_+ \bigr) && \in U_3, \\
		x_4(\lambda, a, r) &:= \exp\bigl( (\lambda e_{12} - \epsilon (a,0))_+ - \epsilon (r' e_{21})_+ \bigr) && \in U_4, \\
		x_5(\lambda) &:= \exp\bigl( (-\lambda e_{22})_+ \bigr) && \in U_5, \\
		x_6(\lambda, a, r) &:= \exp\bigl( -V_{\lambda^\sigma e_{21}, \ov{t}} + V_{(0,a), \ov{t}} - \tfrac{1}{2} V_{r' e_{12}, \ov{t}} \bigr) \ && \in U_6, \\
		x_7(\lambda) &:= \exp\bigl( (-\lambda e_{22})_- \bigr) && \in U_7, \\
		x_0(\lambda, a, r) &:= \exp\bigl( V_{\lambda e_{12} - \epsilon (a,0), t} - \tfrac{\epsilon}{2} V_{r' e_{21}, t} \bigr) && \in U_0,
	\end{alignat*}
	where $\exp$ is as in \cref{def:to groups}.
	It can now be verified that $x_i$ is an isomorphism from $(K, +)$ onto $U_i$ for $i \in \{1,3,5\}$ and from $K \times U_\Lambda$ onto $U_i$ for $i \in \{2,4,6\}$, that the root groups $U_1, \dots, U_6$ satisfy precisely the commutator relations in \cite[Proposition~12.5]{Muhlherr26} (with our value of~$\epsilon$; see also \cite[Observation~12.11]{Muhlherr26}), and that $m_1 = x_7(1)\, x_1(1)\, x_7(1)$ acts on these root groups as in \cite[Observation~12.12]{Muhlherr26}.
	Moreover, if $(a_0, r_0) \in U_\Lambda$ with $r_0 \in K^\times$, then, writing $r_0^{-\sigma} := (r_0^{-1})^\sigma$, the element
	\[ m_6 = x_0(1, \epsilon a_0 r_0^{-1}, r_0^{-\sigma})\, x_6(1, a_0, r_0)\, x_0(1, a_0 r_0^{-\sigma}, r_0^{-\sigma}) \]
	acts on these root groups as in \cite[Observation~12.12]{Muhlherr26}; for $\epsilon = -1$ and $(a_0, r_0) = (0, 1)$, this is the element $m_6 = x_0(1,0,1)\, x_6(1,0,1)\, x_0(1,0,1)$ from \cite[Proposition~12.6]{Muhlherr26}.
	In particular, $m_1$ and $m_6$ are Weyl elements for the roots $\alpha_1$ and~$\alpha_6$, respectively.
	Since $s_{\alpha_1}$ and $s_{\alpha_6}$ generate the Weyl group of~$G_2$, every root then has a Weyl element, so if $\A$ is finite-dimensional and central simple, then $(\A, t)$ is Tits-hexagonalizable by \cref{rem:stability}.

	In case \cref{con:MW:b} of \cref{con:MW} with $\epsilon = -1$, to compare with \cref{ex:herm:B+B}, we first identify $M_2(C^\op)$ with $M_2(C)^\op$ using transposition. Then set $j := \begin{psmallmatrix} 0 & 1 \\ -1 & 0 \end{psmallmatrix}$ and apply a further change of coordinates
	\[ \Psi \colon (x, y) \mapsto (x, jyj^{-1}) \]
	on $M_2(C) \oplus M_2(C)^\op$; this turns the involution on $M_2(K)$ into the ordinary exchange involution.
	Notice that $\Psi(t) = \bigl( \idempotentmatrix, \begin{psmallmatrix} 0 & 0 \\ 0 & 1 \end{psmallmatrix} \bigr)$, in agreement with \cref{pr:assoc idemp}\cref{it:assoc idemp:2}.
\end{proof}

Observe that in case \cref{con:MW:b} of \cref{con:MW}, the conjugate-supplementary idempotent $t$ is \emph{not} central. However, it is ``balanced'' in the sense that the two components of $t = (t_1, t_2) \in B \oplus B^\op$ (where $B = M_2(C)$) are isomorphic (as defined in \cref{rem:isom-idemp}), or equivalently, that the idempotents $t$ and $\overline{t}$ are isomorphic inside the associative algebra~$E$.

We now show that conversely, this balancing condition is forced by hexagonalizability.
\begin{proposition}\label{pr:balanced}
	Let $\A = E \oplus M$ and $t \in E$ (where we allow $M=0$) be as in \cref{th:classification}\cref{it:class:3}\cref{it:class:3a} or \cref{it:class:4}\cref{it:class:4a}, and assume that $(\A, t)$ is Tits-hexagonalizable.
	Then:
	\begin{enumerate}
		\item\label{pr:balanced:1}
			The idempotent $t$ is balanced, i.e., $t$~and~$\overline{t}$ are isomorphic idempotents in $E$.
		\item\label{pr:balanced:2}
			Write $E = B \oplus B^\op$, where $B = M_n(D)$ for some central division algebra~$D$.
			Let $t = (u, \overline{1-u}) \in E$, with $u \in B$.
			Then $u$ and $1-u$ have the same $D$-rank.
			In particular, $n$ is even, and we can write $B = M_2(C)$ for some central simple algebra $C$, such that $u = \idempotentmatrix \in M_2(C)$.
	\end{enumerate}
\end{proposition}
\begin{proof}
	To show \cref{pr:balanced:1}, we will use the existence of Weyl elements to explicitly construct elements $r \in tE\overline{t}$ and $s \in \overline{t}Et$ such that $rs = t$ and $sr = \overline{t}$.
	Choose the short root $\alpha = (0,1)$, where we use the $\Z^2$-coordinates of our $5 \times 5$-grading as in \cref{ta:Lij}.
	By assumption, there exists a Weyl element
	\[ w := v_1 u v_2, \quad v_1, v_2 \in U_{-\alpha}, \ u \in U_\alpha , \]
	acting on the root groups as the reflection $s_\alpha$; in particular, it exchanges the long-root spaces $L_{1, -1}$ and $L_{1, 2}$.
	Notice that both spaces are contained in $L_1 = \A_+$, so we can compute the action of $w$ by performing exponentials of $V$-operators on~$\A$.
	
	Explicitly, write $u \in U_\alpha$ as $u = \exp(\ad_X + \ad_Y)$ where $X := V_{x,t} \in L_{0,1}$ and $Y := \tfrac{1}{2}V_{z,t} \in L_{0,2}$ for some $x = (r,m) \in J = tE\overline{t} \oplus Mt$ and $z \in J' \cap \SS = \overline{t}Et \cap \SS$.
	Let $q := z + \tfrac{1}{2} h(m,m)$, so $q \in \overline{t}Et$.
	
	We first record how $X$ and $Y$ operate on $\A = E \oplus M$, using \cref{def:struct algebra,def:stralg hermitian}.
	So let $(a, n) \in E \oplus M$ be arbitrary. Then using $\overline{t} r = 0 = \overline{r} t$, $r \overline{t} = r$, $\overline{z} = -z$, $\overline{t} z = z = zt$ and $z \overline{t} = 0 = tz$, we obtain
	\begin{align*}
		X(a,n) &= \bigl( ra + h(n,m),\ m\overline{a} + n\overline{r} \bigr) , \\
		Y(a,n) &= \bigl( az,\ 0 \bigr) .
	\end{align*}
	Observe that $X$ and $Y$ commute: since $mz = 0$, we have
	\[ XY(a,n) = \bigl( raz, \ 0 \bigr) = YX(a,n) . \]
	Repeated applications on an element $(b,0)$ now give
	\begin{align*}
		X(b,0) &= \bigl( rb,\  m\overline{b} \bigr), \\
		X^2(b,0) &= \bigl( b h(m,m), \ 2m\overline{b}\overline{r} \bigr), \\
		X^3(b,0) &= \bigl( 3rb h(m,m), \ 0 \bigr) ,
	\end{align*}
	and also $X^4 = Y^2 = 0$. Combining all this, we get
	\[ u(b,0) = \bigl( b + rb + bq + rbq, \ m\overline{b} + m\overline{b}\overline{r} \bigr) . \]
	In particular, if $b \in B^\op = \overline{t}E\overline{t}$, so $b_+ \in L_{1, -1} = B^\op_+ \subseteq \A_+$,
	then the $L_{1,2}$-component of $u(b)$ is equal to
	\[ \pr_{tEt} u(b,0) = rbq . \]
	Now notice that $v_2 \in U_{-\alpha}$ fixes $\overline{t}E\overline{t}$ pointwise and that $v_1 \in U_{-\alpha}$ does not change the $tEt$-component. Since we know that $w(b) \in L_{1,2} = tEt$, we conclude that the Weyl element $w$ maps each $b \in L_{1,-1}$ to
	\[ w(b) = rbq \in L_{1,2} . \]
	Since $w$ exchanges $L_{1,-1}$ and $L_{1,2}$ bijectively, this map $b \mapsto rbq$ is surjective.
	Let $b_0 \in \overline{t}E\overline{t}$ be such that $rb_0q = t$ and set $s := b_0q \in \overline{t}Et$; then
	\[ rs = t, \quad \text{with} \quad r \in tE\overline{t} \quad \text{and} \quad s \in \overline{t}Et . \]
	
	We claim that since $E$ is a finite-dimensional algebra, this already implies that also $sr = \overline{t}$. The key point is that such an algebra is ``Dedekind-finite'', i.e., $ab = 1$ implies $ba = 1$.
	Now $r\overline{r} = 0$ and $\overline{s} s = 0$ imply
	\[ (r + \overline{s})(s + \overline{r}) = rs + \overline{rs} = t + \overline{t} = 1 , \]
	and therefore (using also $\overline{r} r = 0$ and $s \overline{s} = 0$)
	\[ (s + \overline{r})(r + \overline{s}) = sr + \overline{sr} = 1 . \]
	Multiplying this identity on both sides by $\overline{t}$ then gives $sr = \overline{t}$ as claimed, and we conclude that $t$ and $\overline{t}$ are isomorphic idempotents as in \cref{rem:isom-idemp}.
	
	\medskip
	
	We now show \cref{pr:balanced:2}. Recall that an idempotent $u \in B = M_n(D)$ has $D$-rank $\ell$ if and only if $u$ is conjugate in $M_n(D)$ to $\begin{psmallmatrix} I_\ell & 0 \\ 0 & 0 \end{psmallmatrix}$.
	Identifying $M_n(D)$ with $\End_D(V)$ where $V = D^n$ as a right $D$-module, the $D$-rank of $u$ is precisely the $D$-rank of the submodule $uV$ of $V$.
	
	By \cref{pr:balanced:1}, we have $r \in tE\overline{t}$ and $s \in \overline{t}Et$ such that $rs = t$ and $sr = \overline{t}$ in $E = B \oplus B^\op$.
	Now $t = (u, \overline{v})$ with $u+v = 1$, so when we apply the projection $\pi \colon B \oplus B^\op \to B$, we find elements $r_1 := \pi(r) \in uBv$ and $s_1 := \pi(s) \in vBu$ with $r_1 s_1 = u$ and $s_1 r_1 = v$, so $u$ and $v$ are isomorphic idempotents in $B$. In particular, they have the same $D$-rank $\ell$.  Since $V = uV \oplus vV$ and both $uV$ and $vV$ have $D$-rank $\ell$, we see that $n = 2\ell$~is even.
	
	Now set $C := M_\ell(D)$ and choose a basis of $V = uV \oplus vV$ such that $u = \begin{psmallmatrix} I_\ell & 0 \\ 0 & 0 \end{psmallmatrix}$. Then $B \cong M_2(C)$, where the isomorphism maps $u \in B$ to $\idempotentmatrix \in M_2(C)$.
\end{proof}

\begin{example}\label{ex:M3k}
	Let $B := M_3(k)$, let $\A = B \oplus B^\op$ with the exchange involution, and let $t = (e_{11}, \overline{1 - e_{11}})$.
	Then $t$ is not balanced, because $e_{11}$ has $k$-rank $1$ while $1 - e_{11}$ has $k$-rank $2$.
	By \cref{pr:balanced}, $(\A, t)$ is not Tits-hexagonalizable.
	On the other hand, by \cref{ex:B+B}\cref{it:B+B:non-central}, the corresponding $AG_2$-algebra is strong (but not mighty).
	
	Similarly, we can consider a non-degenerate hermitian form over $E = B \oplus B^\op$ and produce an example of a mighty $AG_2$-algebra (where $\A$ is of hermitian type over $E$), where again $(\A, t)$ is not Tits-hexagonalizable.
\end{example}

Finally, we turn to the exceptional examples in \cref{th:classification}\cref{it:class:4}\cref{it:class:4c}.
They correspond to the large exceptional Tits hexagons of type $AG_2$ described in \cite[\S 13]{Muhlherr26}.

\begin{proposition}\label{pr:MW-exc}
	Let $C$ be a composition algebra over $k$, let $\A = \Zorn(k) \otimes C$ and let $t = \idempotentmatrix \otimes 1$ as in \cref{ex:C1xC2}.
	Then $(\A, t)$ is Tits-hexagonalizable of type $AG_2$, and the corresponding root-group sequence in \cref{def:to groups} is that of the large exceptional Tits hexagon of type $AG_2$ from \cite[\S 13]{Muhlherr26}, with $K = C$ in their notation.
	In particular, the odd root groups are isomorphic to $(C, +)$ and the even root groups are isomorphic to the group $\Theta = C^3 \times k^3$ from \cite[13.2]{Muhlherr26}.
\end{proposition}
\begin{proof}
	We give an explicit description of the root groups $U_i$ of $L := \TKK$, labeled as in \cref{not:coordinates}.
	We write the elements of~$\A$ as matrices with entries in $C$ and $C^3$ as in \cref{ex:C1xC2}, and for $c \in C$ and $u \in C^3$, we write
	\[ c e_{11} := \begin{psmallmatrix} c & 0 \\ 0 & 0 \end{psmallmatrix}, \quad c e_{22} := \begin{psmallmatrix} 0 & 0 \\ 0 & c \end{psmallmatrix}, \quad u e_{12} := \begin{psmallmatrix} 0 & u \\ 0 & 0 \end{psmallmatrix}, \quad u e_{21} := \begin{psmallmatrix} 0 & 0 \\ u & 0 \end{psmallmatrix} ; \]
	we view $k^3$ as a subset of~$C^3$.
	Let $g \colon C^3 \times C^3 \to k^3$ be as in \cite[13.2]{Muhlherr26}.
	For $b = (b_1, b_2, b_3) \in C^3$, we write $b^\sigma := (b_1^\sigma, b_2^\sigma, b_3^\sigma)$, and for $(b,s) \in \Theta$, we put $s' := s - \tfrac{1}{2} g(b,b) \in k^3$.
	Now define, for all $d \in C$ and $(b,s) \in \Theta$,
	\begin{alignat*}{2}
		x_1(d) &:= \exp\bigl( (d^\sigma e_{11})_+ \bigr) && \in U_1, \\
		x_2(b,s) &:= \exp\bigl( (b^\sigma e_{21})_+ + ((s'_3, -s'_2, s'_1) e_{12})_+ \bigr) && \in U_2, \\
		x_3(d) &:= \exp\bigl( (d^\sigma e_{11} - d e_{22})_+ \bigr) && \in U_3, \\
		x_4(b,s) &:= \exp\bigl( (b e_{12})_+ + ((-s'_3, s'_2, -s'_1) e_{21})_+ \bigr) && \in U_4, \\
		x_5(d) &:= \exp\bigl( (-d e_{22})_+ \bigr) && \in U_5, \\
		x_6(b,s) &:= \exp\bigl( V_{(b_1^\sigma, -b_2^\sigma, b_3^\sigma) e_{21}, \ov{t}} - \tfrac{1}{2} V_{(s'_3, s'_2, s'_1) e_{12}, \ov{t}} \bigr) \ && \in U_6, \\
		x_7(d) &:= \exp\bigl( (-d e_{22})_- \bigr) && \in U_7, \\
		x_0(b,s) &:= \exp\bigl( V_{(b_3^\sigma, b_2^\sigma, -b_1^\sigma) e_{12}, t} + \tfrac{1}{2} V_{(-s'_1, s'_2, s'_3) e_{21}, t} \bigr) \ && \in U_0,
	\end{alignat*}
	where $\exp$ is as in \cref{def:to groups}.
	It can now be verified that $x_i$ is an isomorphism from $(C, +)$ onto $U_i$ for $i \in \{1,3,5\}$ and from $\Theta$ onto $U_i$ for $i \in \{2,4,6\}$, that the root groups $U_1, \dots, U_6$ satisfy precisely the commutator relations in \cite[13.4]{Muhlherr26}, and that the Weyl elements
	\[ m_1 = x_7(1)\, x_1(1)\, x_7(1) \quad \text{and} \quad m_6 = x_0(e)\, x_6(e)\, x_0(e), \]
	where $e := (0,1,0,0,1,0) \in \Theta$, act on these root groups as in \cite[13.5]{Muhlherr26}.

	Once the root groups in both settings have been identified, the fact that $(\A, t)$ is Tits-hexagonalizable of type $AG_2$ is automatic, because the large exceptional Tits hexagons of type $AG_2$ are $AG_2$-shaped \cite[Definition~6.9]{Muhlherr26}.
\end{proof}

The following elementary lemma will allow us to find Weyl elements in the situation of \cref{pr:MW} when $\epsilon = 1$.
\begin{lemma}\label{le:invertible r0}
	Let $\sigma$ be an orthogonal involution on $K = M_n(k)$ and let $h \in K$ be non-zero with $h^\sigma = h$.
	Then there is some $z \in K$ with $z^\sigma = -z$ such that $h + z$ is invertible.
\end{lemma}
\begin{proof}
	We have $x^\sigma = g^{-1} x^{\mathsf{T}} g$ for some invertible symmetric matrix~$g$ (see, e.g., \cite[Chapter~I]{Knus1998}), and $x \mapsto gx$ maps the elements fixed (resp.\ negated) by~$\sigma$ onto the symmetric (resp.\ skew-symmetric) matrices.
	Hence we may assume that $\sigma$ is the transpose map.
	Replacing $h$ by $P^{\mathsf{T}} h P$ for a suitable invertible matrix~$P$ (see \cite[Chapter~1, Corollary~3.6]{Scharlau1985}), we may then assume that $h = \operatorname{diag}(\lambda_1, \dots, \lambda_n)$, with $\lambda_n \neq 0$ if $n$ is odd.
	Now let $z$ be block-diagonal with $2 \times 2$ blocks $z_i$ in rows and columns $2i-1, 2i$ (and a final entry~$0$ if $n$ is odd), where $z_i = \begin{psmallmatrix} 0 & 1 \\ -1 & 0 \end{psmallmatrix}$ if $\lambda_{2i-1} \lambda_{2i} = 0$ and $z_i = 0$ otherwise.
	Then every diagonal block of $h + z$ has non-zero determinant, so $h + z$ is invertible.
\end{proof}

\medskip

We can now summarize what we have obtained for the strong examples from \cref{th:classification}.
\begin{corollary}\label{co:hexagonalizable}
	Let $\A$ and $t \in \A$ be as in \cref{th:classification}, and assume that the associated $AG_2$-algebra is strong.
	\begin{enumerate}
		\item\label{co:hex:1}
			In cases \cref{it:class:3}\cref{it:class:3b} and \cref{it:class:4}\cref{it:class:4b}, $(\A, t)$ is Tits-hexagonalizable, except in case \cref{it:class:3}\cref{it:class:3b} when $B \cong M_n(k)$ for some odd $n \geq 3$, the involution of~$B$ is orthogonal and $\epsilon = 1$ in \cref{ex:M2B}.
		\item\label{co:hex:2}
			In cases \cref{it:class:3}\cref{it:class:3a} and \cref{it:class:4}\cref{it:class:4a}, $(\A, t)$ is Tits-hexagonalizable if and only if $t$ is balanced.
		\item\label{co:hex:3}
			In case \cref{it:class:4}\cref{it:class:4c}, $(\A, t)$ is Tits-hexagonalizable of type $AG_2$.
	\end{enumerate}
\end{corollary}
\begin{proof}
	\begin{enumerate}[leftmargin=1.5em]
		\item
			Let $K := B$, equipped with its involution~$\sigma$, and let $\epsilon$ be as in \cref{ex:M2B}.
			In case \cref{it:class:3}\cref{it:class:3b}, the pair $(\A, t)$ is precisely as in \cref{con:MW}\cref{con:MW:a} with $V = 0$.
			In case \cref{it:class:4}\cref{it:class:4b}, we have $\A = E \oplus M$ with $E = M_2(K)$, and a standard hermitian Morita argument shows that $M \cong V \oplus V$ for the right $K$-module $V := Mt$, where the hermitian form on $M$ is obtained as in \cref{con:MW} from a non-degenerate $\epsilon$-hermitian form $f \colon V \times V \to K$.
			By \cref{pr:MW}, it remains to find some $(a_0, r_0) \in U_\Lambda$ with $r_0 \in K^\times$.
			If $K$ contains an invertible element~$z$ with $z^\sigma = -\epsilon z$, we can take $(a_0, r_0) = (0, z)$.
			Such an element exists if $\epsilon = -1$ (take $z = 1$), and also if $\epsilon = 1$, unless $K \cong M_n(k)$ for some odd~$n$ and $\sigma$ is orthogonal (compare with \cite[Observation~12.9]{Muhlherr26}); in this exceptional situation, $n \geq 3$ because $\epsilon = -1$ if $B = k$.

			So assume that we are in this exceptional situation.
			In case \cref{it:class:4}\cref{it:class:4b}, we first observe that $f(a_0, a_0) \neq 0$ for some $a_0 \in V$.
			Otherwise, polarizing gives $f(a,c) = -f(c,a) = -f(a,c)^\sigma$, hence $\Trd f(a,c) = 0$ (as $\Trd(x^\sigma) = \Trd(x)$), for all $a, c \in V$; since $f(a, c\lambda) = f(a,c)\lambda$ for all $\lambda \in K$ and the trace form of~$K$ is non-degenerate, this would imply $f = 0$.
			By \cref{le:invertible r0}, applied to $h = \tfrac{1}{2} f(a_0, a_0)$, there is a skew element $z_0 \in K$ such that $r_0 := \tfrac{1}{2} f(a_0, a_0) + z_0$ is invertible; notice that $(a_0, r_0) \in U_\Lambda$.
			In case \cref{it:class:3}\cref{it:class:3b}, on the other hand, the computation in the proof of \cref{pr:balanced}\cref{pr:balanced:1}, which does not depend on the specific form of~$E$, shows that a Weyl element for the short root $(0,1)$ would map each $b \in L_{1,-1} = \overline{t}E\overline{t}$ to $rbq \in L_{1,2} = tEt$, for certain $r \in tE\overline{t}$ and $q \in \overline{t}Et \cap \SS$ (because $M = 0$).
			Writing $q = q_0 e_{21}$, the element $q_0 \in K$ is skew, and it would be invertible because $b \mapsto rbq$ is bijective, which is impossible.
			Hence $(\A, t)$ is not Tits-hexagonalizable in this case.
		\item
			If $(\A, t)$ is Tits-hexagonalizable, then $t$ is balanced by \cref{pr:balanced}.
			Conversely, if $t$ is balanced, then the proof of \cref{pr:balanced}\cref{pr:balanced:2} shows that $B \cong M_2(C)$ for some central simple algebra~$C$, with $u = \idempotentmatrix$.
			Using the change of coordinates $\Psi$ from the proof of \cref{pr:MW} (and the same Morita argument as in~\cref{co:hex:1} if $M \neq 0$), we see that $(\A, t)$ is as in \cref{con:MW}\cref{con:MW:b} with $K = C \oplus C^\op$ and $\epsilon = -1$, so $(\A, t)$ is Tits-hexagonalizable by \cref{pr:MW}.
		\item
			This is \cref{pr:MW-exc}.
		\qedhere
	\end{enumerate}
\end{proof}

\begin{remark}\label{rem:epsilon}
	\begin{enumerate}
		\item\label{rem:epsilon:1}
			If $B$ is associative and contains an invertible element~$s$ with $\overline{s} = -s$, then conjugation by $\begin{psmallmatrix} 1 & 0 \\ 0 & s \end{psmallmatrix}$ is an isomorphism from $M_2(B)$ with $\epsilon = 1$ onto $M_2(B)$ with $\epsilon = -1$ and with the involution $b \mapsto s \overline{b} s^{-1}$ on~$B$, and it fixes $t = \idempotentmatrix$; by transport of structure, it also induces an isomorphism in case \cref{it:class:4}\cref{it:class:4b} of \cref{th:classification}.
			This corresponds to the normalization $\epsilon = -1$ in \cite[Observation~12.9]{Muhlherr26}.
			Consequently, $\epsilon = 1$ can be normalized to $\epsilon = -1$ unless $B \cong M_n(k)$ for some odd $n \geq 3$ and the involution of~$B$ is orthogonal.
		\item\label{rem:epsilon:2}
			Now assume that we are in this case.
			In case \cref{it:class:3}\cref{it:class:3b}, we have $K(\A) \cong \mathfrak{so}(X, Q)$ for a $6n$-dimensional hyperbolic quadratic space $(X, Q)$, and the reason why $(\A, t)$ is not Tits-hexagonalizable is that a Weyl element inducing the longest element of the Weyl group of~$G_2$ would have to come from an isometry of determinant $(-1)^{3n} = -1$; compare with the condition that $d$ be even in \cite[Proposition~7.8(ii)]{Muhlherr26}.
			In case \cref{it:class:4}\cref{it:class:4b}, the quadratic space $(X,Q)$ has an additional non-zero orthogonal summand corresponding to~$M$, on which this determinant can be compensated.
			The resulting Tits hexagons correspond to the classical Tits hexagons of type~$+1$ in the sense of \cite[Observation~12.9]{Muhlherr26}, i.e., those arising from \cite[Proposition~7.8(i)--(iii)]{Muhlherr26} with $d = n$ odd and with $\Delta$ the building of a quadratic form; the element $m_6$ in the proof of \cref{pr:MW} is the one from \cite[Observation~12.12]{Muhlherr26}.
			By \cite[Observation~12.11]{Muhlherr26}, these Tits hexagons do not satisfy hypothesis~(iv) of \cite[Theorem~11.1]{Muhlherr26}; see also \cref{rem:BC3}.
	\end{enumerate}
\end{remark}

\begin{remark}\label{rem:BC3}
	The results in \cref{pr:balanced} and \cref{co:hexagonalizable} have a natural interpretation in terms of the root system~$BC_3$; compare with \cite[11.46]{Muhlherr26}.
	Let $\A = E \oplus M$ be as in case \cref{it:class:3}\cref{it:class:3b} or \cref{it:class:4}\cref{it:class:4b} of \cref{th:classification}, with $B$ associative, $E = M_2(B)$, $t = \idempotentmatrix$, and $M = 0$ in the first case.
	The operator on~$\A$ given by right multiplication with $\overline{t} - t$ on~$E$ and by~$0$ on~$M$ induces a derivation of $L = \TKK$ commuting with the grading derivations of the $5 \times 5$-grading in \cref{ta:Lij}.
	The resulting $\Z^3$-grading of~$L$ is a grading by the root system~$BC_3$ (or by~$C_3$ if $M = 0$), and in suitable coordinates, the $AG_2$-grading is obtained from it via the projection $\Z^3 \to \Z^3/\Z(1,1,1)$.
	This projection maps the roots $\epsilon_i - \epsilon_j$ to the long roots of~$AG_2$, it maps both $\epsilon_i + \epsilon_j$ and $-\epsilon_k$ (where $\{i,j,k\} = \{1,2,3\}$) to the same short root~$\alpha$, and it maps $-2\epsilon_k$ to~$2\alpha$.
	Correspondingly, every short root space splits into an $E$-part and an $M$-part; for instance, $L_{1,0} = (tE\overline{t})_+ \oplus (Mt)_+$.

	The Weyl group of~$G_2$ lifts to the subgroup of the Weyl group of~$BC_3$ generated by the permutations of the~$\epsilon_i$ and by~$-1$.
	The reflection in a short root of~$G_2$ lifts to an element of the form $-(i\ j)$, which involves an odd number of sign changes.
	Accordingly, a Weyl element for a short root of~$AG_2$ is essentially the product of a Weyl element for a root $\pm(\epsilon_i + \epsilon_j)$ and a Weyl element realizing the sign change of~$\epsilon_k$, i.e., a Weyl element for $\pm\epsilon_k$ or $\pm 2\epsilon_k$ (where the root group of $\pm\epsilon_k$ contains that of $\pm 2\epsilon_k$, as in \cref{def:to groups}).
	This is exactly the factorization $b \mapsto rbq$ in the proof of \cref{pr:balanced}.
	The invertibility of $r \in tE\overline{t}$ concerns the roots $\pm(\epsilon_i + \epsilon_j)$; in cases \cref{it:class:3}\cref{it:class:3a} and \cref{it:class:4}\cref{it:class:4a}, this is what forces $t$ to be balanced.
	The invertibility of $q = z + \tfrac{1}{2} h(m,m)$ concerns the sign change, where $z$ comes from the root space of $\pm 2\epsilon_k$ and $m$ from that of~$\pm\epsilon_k$.

	Now let $\Delta$ be the set of roots of the $BC_3$-grading for which Weyl elements exist.
	Then $\Delta$ always contains the roots $\pm\epsilon_i \pm \epsilon_j$, which form a root system of type~$D_3$.
	If $\epsilon = -1$, or if $\epsilon = 1$ and $B$ contains an invertible skew element, then $\Delta$ consists of all roots, so $\Delta$ is of type $C_3$ or $BC_3$, according to whether $M = 0$ or not.
	If $\epsilon = 1$, $B \cong M_n(k)$ for some odd~$n$ and the involution of~$B$ is orthogonal, on the other hand, the roots $\pm 2\epsilon_k$ never belong to~$\Delta$, so $\Delta$ is of type~$B_3$ if $M \neq 0$ (by the proof of \cref{co:hexagonalizable}\cref{co:hex:1}) and of type~$D_3$ if $M = 0$.
	Since $-1$ belongs to the Weyl groups of type $B_3$ and~$C_3$, but not to the Weyl group of type~$D_3$, this explains the exception in \cref{co:hexagonalizable}\cref{co:hex:1}.
	Notice that the case where $\Delta$ is of type~$B_3$ is excluded by hypothesis~(iv) of \cite[Theorem~11.1]{Muhlherr26}, which amounts to the existence of Weyl elements for the roots~$\pm 2\epsilon_k$; compare with \cite[Observation~12.11]{Muhlherr26}.
\end{remark}

\clearpage

\appendix
\renewcommand{\thesubsection}{\thesection.\arabic{subsection}}

\newcommand{\JH}{J_{\mathrm h}}
\newcommand{\JS}{J_{\mathrm s}}
\newcommand{\JpH}{J'_{\mathrm h}}
\newcommand{\JpS}{J'_{\mathrm s}}
\newcommand{\ax}[1]{\textup{(#1)}}

\section{An axiom list for \texorpdfstring{$AG_2$}{AG2}-algebras}\label{se:appendix}

In this appendix, we give a self-contained list of axioms for $AG_2$-algebras that does not refer to the corresponding structurable algebra. We have produced this list with the help of AI (Claude Opus 5.5 by Anthropic), although it could certainly also have been produced by hand. It is not clear yet how useful this list is, which is precisely why it appears as an appendix only.

This axiom system should be compared with the axioms of a cubic norm structure (or even better, of a cubic norm \emph{pair} as in \cite{Smet2025,DeMedts2026}). See also \cref{rem:axioms}\cref{rem:axioms:H} below.

\subsection{The data}

Let $B$ be a unital $k$-algebra and let $\overline{B}=\{\overline b\mid b\in B\}$ be a copy of $B$ with multiplication $\overline a\,\overline b:=\overline{ba}$, so that $\overline{B}\cong B^{\op}$. As in \cref{def:bimodules} and \cref{def:AG2-algebra}, we are given:
\begin{itemize}[leftmargin=2em]
\item a vector space $J$ with an involutory linear map $v\mapsto\overline v$, a left action $B\times J\to J$, $(a,v)\mapsto av$, and a right action $J\times\overline{B}\to J$, $(v,\overline a)\mapsto v\overline a$, related by $\overline{av}=\overline v\,\overline a$;
\item a vector space $J'$ with an involutory linear map $v'\mapsto\overline v'$, a right action $J'\times B\to J'$, $(v',a)\mapsto v'a$, and a left action $\overline{B}\times J'\to J'$, $(\overline a,v')\mapsto \overline a v'$, related by $\overline{v'a}=\overline a\,\overline v'$;
\item bilinear maps $\times\colon J\times J\to J'$, $\times'\colon J'\times J'\to J$ and $T\colon J\times J'\to B$.
\end{itemize}
We put $T'(v',v):=\overline{T(\overline v,\overline v')}\in\overline{B}$ (see \cref{le:T'}), and we always assume
\begin{equation*}\label{eq:Inv}\tag{Inv}
\overline{v\times w}=\overline w\times\overline v,\qquad \overline{v'\times' w'}=\overline w'\times'\overline v' \qquad (v,w\in J,\ v',w'\in J').
\end{equation*}
Condition~\eqref{eq:Inv} is exactly what is needed for the map $x\mapsto\overline x$ of \cref{pr:Peirce multiplication new} to be an involution of
\[ \A = \begin{pmatrix} B & J\\ J' & \overline{B}\end{pmatrix}; \]
so \eqref{eq:Inv} is implicit in \cref{def:AG2-algebra}. Note that the multiplication of $\A$ is then simply
\bgroup
\renewcommand{\arraystretch}{1.15}
\[
\begin{array}{c|cccc}
\cdot & B & J & J' & \overline{B}\\\hline
B & ab & av & 0 & 0\\
J & 0 & v\times w & T(v,v') & v\overline a\\
J' & v'a & T'(v',v) & v'\times' w' & 0\\
\overline{B} & 0 & 0 & \overline a v' & \overline a\overline b
\end{array}
\]
\egroup
We write $J=\JH\oplus\JS$ and $J'=\JpH\oplus\JpS$ for the decompositions into hermitian ($\overline v=v$) and skew ($\overline v=-v$) elements.

For elements $x,y,z\in B\cup J\cup J'\cup\overline{B}$, we write $[x,y,z]=(xy)z-x(yz)$ for their associator \emph{computed in $\A$}, i.e., with the table above. This is only a shorthand; the associators that occur below are
\begin{align*}
[a,b,v]&=(ab)v-a(bv), & [a,v,\overline b]&=(av)\overline b-a(v\overline b),\\
[v',a,b]&=(v'a)b-v'(ab), & [\overline b,v',a]&=(\overline bv')a-\overline b(v'a),\\
[a,v,v']&=T(av,v')-aT(v,v'), & [v,\overline a,v']&=T(v\overline a,v')-T(v,\overline av'),\\
[v,v',a]&=T(v,v')a-T(v,v'a).
\end{align*}

\subsection{The axioms}

\begin{theorem}\label{th:axioms}
Let $\Xi = (B,J,J',\times,\times',T)$ be as above, satisfying \eqref{eq:Inv}. Then $\Xi$ is an $AG_2$-algebra \textup(i.e., $\A$ is a structurable algebra\textup) if and only if the following identities hold for all
$a,b,c\in B$, $v,w\in J$, $v',w'\in J'$, $s\in\JS$, $s'\in\JpS$, $x,y,z\in\JH$ and $x',y',z'\in\JpH$.

\medskip
\noindent\textbf{Group I} \textup{(}skew-alternativity\textup{):}
\begin{enumerate}[leftmargin=3.2em,itemsep=3pt]
\item[\ax{A}] $B$ is alternative.
\item[\ax{S}] $s\times w=-\,w\times s$ \ and \ $s'\times' w'=-\,w'\times' s'$.
\item[\ax{M}] $[a,b,v]=-[b,a,v]=[a,v,\overline b]$ \ and \ $[a,b,s]=0$.
\item[\ax{M$'$}] $[v',b,a]=-[v',a,b]=[\overline b,v',a]$ \ and \ $[s',a,b]=0$.
\item[\ax{T}] $T(av,v')-aT(v,v')=T(v\overline a,v')-T(v,\overline av')=T(v,v')a-T(v,v'a)$, and this element vanishes whenever $v\in\JS$ or $v'\in\JpS$.
\item[\ax{X}] $(av)\times w-v\times(aw)=(v\overline a)\times w-\overline a\,(v\times w)$.
\item[\ax{X$'$}] $(\overline av')\times' w'-v'\times'(\overline aw')=(v'a)\times' w'-a\,(v'\times' w')$.
\item[\ax{N}] $T(s,v\times w)=-\,T(v,s\times w)$.
\item[\ax{N$'$}] $T'(s',v'\times' w')=-\,T'(v',s'\times' w')$.
\item[\ax{J}] $(v\times s)\times' v'=v\,T'(v',s)-T(v,v')\,s$ \ and \\ $s'\times'(v\times w)=T(v,s')\,w-v\,T'(s',w)$.
\item[\ax{J$'$}] $(v'\times' s')\times v=v'\,T(v,s')-T'(v',v)\,s'$ \ and \\ $s\times(v'\times' w')=T'(v',s)\,w'-v'\,T(s,w')$.
\end{enumerate}

\medskip
\noindent\textbf{Group II} \textup(on hermitian elements\textup{):}

Put $N(x,y,z):=T(x,y\times z)\in B$ and $N'(x',y',z'):=T'(x',y'\times' z')\in\overline{B}$, and
\begin{align*}
	x^\sharp &:= \tfrac12\, x\times x\in\JpH, & N(x) &:= \tfrac13\, T(x,x^\sharp)\in B, \\
	x'^\sharp &:= \tfrac12\, x'\times' x'\in\JH, & N'(x') &:= \tfrac13\, T'(x',x'^\sharp)\in\overline{B},
\end{align*}
so that $N(x,x,x)=6N(x)$ and $N'(x',x',x')=6N'(x')$.
\begin{enumerate}[leftmargin=3.2em,itemsep=3pt]
\item[\ax{H1}] $N(x,y,z)-N(y,x,z)=N(z,x,y)-N(z,y,x)$.
\item[\ax{H1$'$}] $N'(x',y',z')-N'(y',x',z')=N'(z',x',y')-N'(z',y',x')$.
\item[\ax{H2}] $(z\times x)\times' y'-y'\times'(x\times z)=T(x,y')z+T(z,y')x-x\,T'(y',z)-z\,T'(y',x)$.
\item[\ax{H2$'$}] $(z'\times' x')\times y-y\times(x'\times' z')=T'(x',y)z'+T'(z',y)x'-x'\,T(y,z')-z'\,T(y,x')$.
\item[\ax{H3}] $(y\times x)\times' x^\sharp=T(y,x^\sharp)\,x+2N(x)\,y-y\,\overline{N(x)}$.
\item[\ax{H4}] $(y'\times' x^\sharp)\times x=T'(y',x)\,x^\sharp+2\,\overline{N(x)}\,y'-y'N(x)$.
\item[\ax{H3$'$}] $(y'\times' x')\times x'^\sharp=T'(y',x'^\sharp)\,x'+2N'(x')\,y'-y'\,\overline{N'(x')}$.
\item[\ax{H4$'$}] $(y\times x'^\sharp)\times' x'=T(y,x')\,x'^\sharp+2\,\overline{N'(x')}\,y-y\,N'(x')$.
\end{enumerate}
\end{theorem}

\begin{remark}
Replacing $t$ by $\overline t$ replaces $(B,J,J',\times,\times',T)$ by $(\overline{B},J',J,\times',\times,T')$, and the list is invariant under this substitution: it interchanges each unprimed family with its primed version, while axioms \ax{A}, \ax{S} and \ax{T} are mapped to equivalent forms of themselves. So, up to this symmetry, there are $11$ families of axioms.
\end{remark}

We collect some useful remarks on the meaning of these axioms.

\begin{remark}\label{rem:axioms}
\begin{enumerate}[leftmargin=2em]
\item Group~I is \emph{equivalent} to the skew-alternativity of $\A$ (\cref{pr:skew alternative}): $[s,x,y]=-[x,s,y]=[x,y,s]$ for all $s\in \SS$, $x,y\in \A$. More precisely, each family in Group~I is the set of Peirce components of this identity of one fixed Peirce type, see \cref{tab:SA}.
\item Group~II is what remains of Allison's two conditions on hermitian elements (see \S\ref{se:proof}). Only hermitian elements of $J$ and $J'$ occur; the hermitian elements $b+\overline b$ ($b\in B$) give nothing new.
\item Axiom \ax{M} says that the associator $[a,b,v]$ is alternating in $a,b$, only depends on the hermitian part of $v$, and can be moved to the ``middle'' position. In combination with the other axioms it also forces $\overline{[a,b,v]}=[a,b,v]$. Axiom~\ax{T} says that $T$ is ``$B$-balanced up to one associator'', which only sees hermitian parts.
\item\label{rem:axioms:H} Axioms \ax{H3} and \ax{H4} (and their primed versions) are the analogues of the adjoint identity $x^{\sharp\sharp}=N(x)x$ for cubic norm structures. Indeed, for $\A=M(J,\eta)$ as in \cref{ex:skewdim1}, we have $B=k$, the involutions on $J$ and $J'$ are trivial, $\times'=\eta\,\times$, and the trace map of the $AG_2$-algebra is $\eta T_J$, where $T_J$ denotes the trace form of the cubic norm structure~$J$. Hence $x^\sharp$ is the adjoint of $x$ in $J$ and $N(x)=\eta N_J(x)$, where $N_J$ denotes the norm of~$J$, and after dividing by~$\eta$, axiom \ax{H3} becomes
\[ x^\sharp\times(x\times y)=T_J(x^\sharp,y)\,x+N_J(x)\,y, \]
which is the linearization of $x^{\sharp\sharp}=N_J(x)\,x$; this is also the identity used in the proof of \cref{co:G2-frame}.
Similarly, axiom \ax{H4} becomes $x\times(x^\sharp\times y')=T_J(x,y')\,x^\sharp+N_J(x)\,y'$.
\item In the weak case ($\JS=\JpS=0$) all axioms containing $s$ or $s'$ are void, the involutions on $J,J'$ are trivial, so $v\overline a=av$ and $\overline av'=v'a$, and axiom~\ax{M} reduces to $(ab)v=b(av)$: $J$ is a right $B$-module in the associative sense via $v\cdot a:=av$ (and similarly for $J'$). This is consistent with \cref{th:classification}\cref{it:class:2}, where $B$ is associative.
\item In the feeble case, only axiom \ax{A} remains.
\end{enumerate}
\end{remark}

\subsection{Sketch of the proof}\label{se:proof}

The crucial ingredient is Allison's criterion \cite[Theorem~13]{Allison1978}:
A unital algebra $\A$ with involution is structurable if and only if
\begin{enumerate}[label=(\roman*)]
\item\label{it:pf:1} $[s,x,y]=-[x,s,y]$ for all $s\in \SS$ and $x,y\in \A$ \ (this implies $=[x,y,s]$ by applying the involution);
\item\label{it:pf:2} $[x,y,z]-[y,x,z]=[z,x,y]-[z,y,x]$ for all $x,y,z\in \HH$;
\item\label{it:pf:3} $D_{x^2,x}(y)=0$ for all $x,y\in \HH$, where
	\[ D_{x,y}(z) = \tfrac13 [[x,y] + [\overline x, \overline y],z] + [z,y,x] - [z,\overline x,\overline y] . \]
\end{enumerate}
All of these have degree $\leq 4$, which is the key reduction compared with the $5$-linear defining identity. Since $\Char(k) \neq 2,3$, \ref{it:pf:3} may be replaced by its full linearization in $x$. As all conditions are multilinear, it suffices to check them on the homogeneous spanning sets
\[ \SS = \{b-\overline b\}\oplus\JS\oplus\JpS,\qquad \HH = \{b+\overline b\}\oplus\JH\oplus\JpH, \]
with $x,y$ in \ref{it:pf:1} running over $B\cup J\cup J'\cup\overline{B}$, and then to project onto the four Peirce components.

\medskip

The first step consists of showing that conditions \ref{it:pf:1} translate to the axioms in Group I.
Table~\ref{tab:SA} lists, for each family of Group~I, the choices of $(s;x,y)$ in~\ref{it:pf:1} which produce it; conversely, every component of \ref{it:pf:1} follows from axioms \ax{A} and~\ax{S} together with the single family of the same Peirce type.
For instance, for $s=b-\overline b$ and $(x,y)=(v,v')$ one gets
$[s,v,v']=[b,v,v']$, $[v,s,v']=-[v,\overline b,v']$, $[v,v',s]=[v,v',b]$, which is the first part of \ax{T}; for $s\in\JS$ and $(x,y)=(v,t)$ one gets $[s,v,t]=s\times v$ and $[v,s,t]=v\times s$, which is \ax{S}.

\begin{table}[ht]
\centering
\begin{tabular}{@{}lll@{}}
\toprule
family & Peirce type & produced by $(s;\,x,y)$ in \ref{it:pf:1} \\\midrule
\ax{A} & $BBB \to B$ & $(b-\overline b;\ a,c)$\\
\ax{S} & $JJB \to J'$ & $(s;\ v,t)$\\
\ax{M} & $BBJ \to J$ & $(b-\overline b;\ a,v)$, $(s;\ a,b)$\\
\ax{T} & $BJJ' \to B$ & $(b-\overline b;\ v,v')$, $(s;\ a,v')$, $(s';\ v,a)$\\
\ax{X} & $BJJ \to J'$ & $(b-\overline b;\ v,w)$\\
\ax{N} & $JJJ \to B$ & $(s;\ v,w)$\\
\ax{J} & $JJJ' \to J$ & $(s;\ v,v')$, $(s';\ v,w)$\\
\ax{J$'$} & $JJ'J' \to J'$ & $(s';\ v',v)$, $(s;\ v',w')$\\
\bottomrule
\end{tabular}
\caption{Skew-alternativity, sorted by Peirce type (the primed families are the mirror images).}\label{tab:SA}
\end{table}

Once this is done, the next step is to show that \ref{it:pf:2} and \ref{it:pf:3} together are equivalent with the axioms in Group II, \emph{assuming the axioms in Group I}.

For $x,y,z\in\JH$, the $B$-component of \ref{it:pf:2} is exactly \ax{H1} (the $\overline{B}$-component is its conjugate). For $x,z\in\JH$ and $y'\in\JpH$, condition \ref{it:pf:2} for $(x,y',z)$ is \ax{H2}. All other types in \ref{it:pf:2} follow from Group~I.
For \ref{it:pf:3}, take $x,y\in\JH$. Then $x^2=x\times x=2x^\sharp\in\JpH$ and $xx^2=T(x,x\times x)=6N(x)$, so $[x^2,x]=6\bigl(\overline{N(x)}-N(x)\bigr)$, hence
\begin{multline*}
	[[x^2,x],y]=-6N(x)\,y-6y\,\overline{N(x)},\qquad [y,x,x^2]=2(y\times x)\times' x^\sharp, \\ [y,x^2,x]=2T(y,x^\sharp)\,x-6y\,\overline{N(x)},
\end{multline*}
so that
\[ D_{x^2,x}(y)=2\bigl((y\times x)\times' x^\sharp-T(y,x^\sharp)\,x-2N(x)\,y+y\,\overline{N(x)}\bigr), \]
and therefore $D_{x^2,x}(y)=0$ if and only if \ax{H3} holds. In the same way $D_{x^2,x}(y')=0$ for $x\in\JH$, $y'\in\JpH$ is \ax{H4}. The components of \ref{it:pf:3} with $x$ in $\JH\oplus\JpH$ ``mixed'' (or involving $b+\overline b$) follow from Group~I and \ax{H1}, \ax{H1$'$}, \ax{H2}, \ax{H2$'$}.

\bigskip

\subsection*{AI Declaration}

This paper has been written by the authors, with the following exceptions.
\begin{itemize}\setlength{\itemindent}{-3ex}
	\item We have used ChatGPT Sol 6.1 (OpenAI) to assist with the narrative of the abstract and the introduction.
	\item We have used ChatGPT Sol 5.5 to assist with some parts of the proofs in \cref{se:AG2 to stralg,se:AG2 G2}.
	\item \cref{le:Phi} and its use in the proof of \cref{th:AG2-A2}, which make this proof valid in characteristic~$5$ as well, have been suggested by Claude Opus 5.5 (Anthropic).
	\item The comparison with Yaseen's work \cite{Yaseen22} has been spelled out by Claude Opus 5.5 in \cref{rem:Yaseen}.
	\item We have used ChatGPT Sol 5.5 to assist with the proof of \cref{pr:balanced}. In addition, \cref{rem:stability} has been written by ChatGPT Sol 6.1.
	\item The explicit descriptions of the root groups and Weyl elements in the proofs of \cref{pr:MW,pr:MW-exc} have been worked out and verified by Claude Opus 5.5. The same holds for the treatment of the case $\epsilon = 1$ in \cref{se:tits hex} (in particular in \cref{co:hexagonalizable}\cref{co:hex:1} and \cref{rem:epsilon}). \Cref{rem:BC3} has been written by Claude Opus 5.5, based on a suggestion of the authors.
	\item \cref{se:appendix} has been produced by Claude Opus 5.5 and subsequently polished by us. (We had tried this before with ChatGPT Sol 5.5, but it had failed at producing a useful list of axioms.)
\end{itemize}
When the paper was essentially finished, we asked Claude Opus 5.5 to proofread the manuscript in detail. This resulted in a six-page long list of useful suggestions, which we have carefully addressed; a second round of proofreading of the revised version led to a few more corrections.

\bibliographystyle{alpha}
\bibliography{AG2.bib}

\end{document}